\documentclass[11pt,a4paper]{amsart}
\usepackage[textwidth=158mm,textheight=243mm,centering]{geometry}
\usepackage{amsmath,amssymb,amsthm,mathtools}
\usepackage[T1]{fontenc}
\usepackage{lmodern,microtype,needspace}
\usepackage[hidelinks]{hyperref}
\usepackage{xurl}
\hypersetup{pdftitle={Lorentzian entropy, subgaussian matrix concentration, and sharp functional inequalities for all matroids},pdfauthor={Eren Ercan}}
\numberwithin{equation}{section}
\newtheorem{theorem}{Theorem}[section]
\newtheorem{lemma}[theorem]{Lemma}
\newtheorem{corollary}[theorem]{Corollary}
\newtheorem{proposition}[theorem]{Proposition}
\theoremstyle{definition}
\theoremstyle{remark}\newtheorem{remark}[theorem]{Remark}
\newcommand{\one}{\mathbf1}
\newcommand{\E}{\mathbb E}
\newcommand{\R}{\mathbb R}
\newcommand{\supp}{\operatorname{supp}}
\DeclareMathOperator{\Tr}{Tr}
\DeclareMathOperator{\diag}{diag}
\DeclareMathOperator{\Spec}{Spec}
\DeclareMathOperator{\Var}{Var}
\DeclareMathOperator{\rank}{rank}
\DeclareMathOperator{\Sym}{Sym}
\title[Lorentzian entropy and matrix concentration]{Lorentzian entropy, subgaussian matrix concentration, and sharp functional inequalities for all matroids}
\author{Eren Ercan}
\email{eren321@gmail.com}
\date{31 October 2026}
\microtypesetup{expansion=false}
\allowdisplaybreaks[2]
\begin{document}
\begin{abstract}
We prove complete entropy contraction for homogeneous Lorentzian laws
and derive subgaussian matrix concentration and sharp functional inequalities
for every matroid. Uniform retention of $\ell$ out of $r$ particles
contracts matrix entropy by the optimal universal factor $\ell/r$,
in every finite matrix dimension and positive field, including repeated labels.
The matrix Laplace comparison retains matrix variance and arbitrary Hermitian
offsets. It gives subgaussian bounds for arbitrary Hermitian observables with
bounded exchange differences, and Chernoff and Bernstein bounds for basis sums.

Every unit-edge matroid basis graph satisfies the cube's small-set and Gross
inequalities, with equality classified and optimal universal complete
entropy-decay rate four. Its spectrum lies in $\{0,2\}\cup[3,\infty)$,
with binary direct-factor contrasts spanning the eigenspace at two.
We classify extremal down-up eigenfunctions, resolving Alaniz Pintos's
degree-one conjecture and the equality problem for the Feder--Mihail walk.

For degree at least two, every matrix entropy extremizer extends uniquely
from its two-particle marginal. Equality depends only on the normalized pair
law, with noncommutativity confined to mutually exclusive labels. We classify
power-entropy extremizers for $1<p<2$ and give a fourth-order criterion for
affine equality. Retained-level Gibbs reconstruction has sharp
dimension-independent stability. For multiaffine laws whose generating
polynomial and complement are Lorentzian, we determine the exact simultaneous
all-field quantum marginal and matrix H\"older regions and classify the sharp
half-marginal fields. Connected sharp laws are determinantal through rank
three; stable nondeterminantal examples exist in every higher rank.
\end{abstract}
\maketitle
\section{Introduction}

Uniform retention of $\ell$ out of $r$ particles from a homogeneous
Lorentzian law contracts matrix entropy by the sharp universal factor
$\ell/r$, at every positive field and in every finite matrix dimension.
Repeated labels are allowed. For $r\ge2$, every matrix extremizer extends
uniquely from its two-particle marginal; its noncommutative structure is confined
to labels that never occur together.

For matroids, complete contraction gives nonlinear subgaussian concentration
and matrix Chernoff and Bernstein bounds. The equality theory identifies
all extremal down-up eigenfunctions and, for the unit-edge basis graph,
the spectral separation $\Spec(L_G)\subseteq\{0,2\}\cup[3,\infty)$.
The eigenvalue-two space consists exactly of binary direct-factor
contrasts, and the optimal universal complete entropy-decay rate is four.

Contraction, profiles and trace duality form the common entrance
(Sections~\ref{sec:complete-contraction}--\ref{sec:duality}).
The concentration bounds and all-field marginal region can then be read directly.
For spectral applications, continue through scalar translations to
Sections~\ref{sec:deletion-dynamics} and~\ref{sec:exchange-edges};
the quantum equality route passes through reconstruction and
Sections~\ref{sec:universal-equality} and~\ref{sec:power-entropies}.

\subsection{Complete contraction}
Let $f(z)=\sum_{|\alpha|=r}a_\alpha z^\alpha$ be a nonzero
homogeneous Lorentzian polynomial of degree $r\ge1$, and fix a
positive field $\lambda$. Its coefficient law is
\[
 \mu(\alpha)=\frac{a_\alpha\lambda^\alpha}{f(\lambda)}.
\]
An occupation $\alpha$ contains $\alpha_i$ particles of label $i$.
Retain a uniformly chosen set of $\ell$ particles, write $\mu_\ell$
for the retained law, and let
$F_\ell(\beta)=\E[F(\alpha)\mid\beta]$ under the reference
deletion coupling. For a positive semidefinite matrix function put
\[
 \mathcal H_\mu(F)=\E_\mu\Tr(F\log F)
             -\Tr[(\E_\mu F)\log(\E_\mu F)].
\]
No trace normalization is imposed.

\begin{theorem}[Complete particle-deletion contraction]
\label{thm:intro-complete-deletion}
For every $0\le\ell\le r$ and every positive semidefinite matrix
function of any finite size,
\begin{equation}\label{eq:intro-matrix-entropy}
 \mathcal H_{\mu_\ell}(F_\ell)\le\frac\ell r\mathcal H_\mu(F).
\end{equation}
The coefficient is uniform in the positive field and matrix size,
and is optimal over the class of Lorentzian laws.
\end{theorem}

The word \emph{complete} records the uniformity in matrix size.
For a normalized classical--quantum state, the entropy is relative
entropy from $\mu\otimes\rho_Q$, using the actual quantum marginal
$\rho_Q$. The input classical marginal may differ from $\mu$.
If it equals $\mu$, the entropy is mutual information. Deleting one
particle and refilling conditionally gives complete modified
logarithmic Sobolev constant at least $1/r$ for the generator $I-K$.

Lieb's three-matrix inequality turns the Lorentzian Hessian signature
into a quadratic entropy inequality. Conditioning on two unmarked
particles identifies this defect with the second difference of the
deletion profile. Convexity then gives contraction, while the same
defects control strictness and equality.

Cryan--Guo--Mousa established scalar deletion contraction and scalar
concentration for homogeneous strongly log-concave subset laws
\cite{R4CGM}. The scalar entropy framework also includes
\cite{AOV,AJKPV}; Lee--Seo give the one-particle occupation version
with repeated labels \cite[Lemma 3.5 and Remark 3.6]{LeeSeo2026}.
The quadratic converse is exact: scalar contraction, complete
contraction and one positive eigenvalue are equivalent. Infinitesimal
tests on coordinate derivatives at one positive field therefore
recognize Lorentzianity (Corollary~\ref{cor:one-field-recognition}).

\subsection{Quadratic models for entropy equality}
For $r\ge2$, absorb the positive field and normalize $f(\one)=1$.
Its translation space is
\[
 \mathcal N_f=\{h:D_hf=0\}=\ker\nabla^2f(\one).
\]
It determines all scalar extremizers. For matrices, the mean and the
one-particle conditional densities also determine every extremizer,
by Gibbs reconstruction. The following theorem constructs every
equality density from its two-particle marginal.

\Needspace{17\baselineskip}
\begin{samepage}
\begin{theorem}[Quadratic extension]
\label{thm:intro-quadratic-extension}
Restrict to active labels and put $W=\nabla^2f(\one)/(r(r-1))$ and $f_2(z)=z^{\mathsf T}Wz$.
For each finite matrix dimension, two-particle conditional expectation
is a bijection between the complete entropy equality densities of
$f$ and those of $f_2$, normalized by $\Tr\E F=1$.
Given a quadratic extremizer $E_{ij}$, set
$\sigma=\sum_{ij}W_{ij}E_{ij}$. On $\supp\sigma$, put
\[
 G_i=\frac{\sum_jW_{ij}E_{ij}}{(W\one)_i},\qquad
 B_i=\sigma^{-1/2}G_i\sigma^{-1/2}.
\]
Its unique extension is
\[
 F_\alpha=\sigma^{1/2}\left(\prod_iB_i^{\alpha_i}\right)\sigma^{1/2}.
\]
Factors occurring in the same occupation commute. The extension has
$\mathcal H_\mu(F)/r=\mathcal H_W(E)/2$. More generally, any two Lorentzian laws with
the same normalized $W$ have canonically identified equality families,
even when their degrees differ.
\end{theorem}
\end{samepage}

The finite marginal test and normal form
(Theorem~\ref{thm:quadratic-extension}) combine Petz recovery,
Hayden--Jozsa--Petz--Winter's Markov structure theorem \cite{HJPW},
and Lorentzian translations. For matroid basis laws, the quantum
register is independent of the simplified basis; only parallel-class
choices carry equality information (Corollary~\ref{cor:matroid-quantum-parallel}).
When every two distinct active labels can occur together, logarithmic
extremizers are simultaneously diagonalizable. In this case, fixing the classical
marginal leaves only constant equality encodings, although a singular Hessian
makes $\ell/r$ the sharp supremum, approached by qubit encodings.
In general, informative fixed-marginal equality exists exactly when
two distinct active labels never co-occur, and a qubit suffices.

Power equality for $1<p<2$ places the exclusive classes in orthogonal
reducing subspaces; those inputs have linear profiles for every trace
function (Theorem~\ref{thm:universal-matrix-power-equality}).
For affine matrix translation families, a quartic sum of squares
decides exact equality for the whole family
(Theorem~\ref{thm:complete-quartic-matrix-obstruction}).
This criterion applies to every homogeneous occupation law with
the stated translation identities.

\subsection{Extremal matroid walks and spectral separation}
For the spectral statements, let $M$ be a loopless rank-$r$ matroid
with at least two bases,
$r\ge1$, and let $\mu$ be uniform on $\mathcal B(M)$. The
\emph{down-up walk} $P_{\mathrm{du}}$ removes one uniformly chosen
element and chooses a uniform basis completing the retained face,
allowing the removed element to return. Write $C$ for the parallel
classes. The basis graph $G$ joins bases differing by one exchange,
and its positive unit-edge Laplacian is
\[
 (L_Gg)(B)=\sum_{B'\sim B}\bigl(g(B)-g(B')\bigr).
\]
The constant-edge Feder--Mihail walk is
$P_{\mathrm{FM}}=I-pL_G$, where $0<p\le1/\Delta(G)$ and
$\Delta(G)$ is the maximum degree.
We use the gap $1-\lambda_2(P)$, where $\lambda_2(P)$ is the
second-largest eigenvalue.

\Needspace{18\baselineskip}
\begin{samepage}
\begin{theorem}[Complete classification of spectral equality]
\label{thm:intro-matroid-rigidity}
With these normalizations, the following statements hold.
\begin{enumerate}
\item The down-up gap is at least $1/r$, and
\begin{equation}\label{eq:intro-down-up-equality}
 \ker\bigl(P_{\mathrm{du}}-(1-1/r)I\bigr)
 =\left\{B\mapsto\sum_{i\in B}h_i:
            \sum_{i\in C}h_i=0\text{ for every }C\right\}.
\end{equation}
This space has dimension $\sum_C(|C|-1)$. Consequently the gap
equals $1/r$ if and only if $M$ has a nontrivial parallel class.
\item If $M$ has $k$ direct summands isomorphic to $U_{1,2}$, then
\begin{equation}\label{eq:intro-spectral-separation}
 \Spec(L_G)\subseteq\{0,2\}\cup[3,\infty),\qquad
 \dim\ker(L_G-2I)=k.
\end{equation}
The eigenvalue-two space is spanned by the $k$ binary contrasts
$\mathbf1_{\{a\in B\}}-\mathbf1_{\{b\in B\}}$ on these summands.
The Feder--Mihail gap equals $2p$ exactly when $k>0$; if $k=0$,
it is at least $3p$. The threshold three is sharp.
\end{enumerate}
\end{theorem}
\end{samepage}

Part~(1) proves more than the existence of a degree-one extremizer
asked for in \cite[Conjecture 8.1]{Pintos2026}: it identifies every
extremizer. It gives an exact criterion on the original matroid
for the structural question in \cite[Conjecture 8.2]{Pintos2026}.
Part~(2) answers \cite[Problem 8.4]{Pintos2026}. Parallel elements
alone suffice for down-up equality; the unit-edge walk singles out
binary direct components. The proofs appear in
Corollaries~\ref{cor:pintos-down-up} and~\ref{cor:pintos-feder-mihail}
and Theorem~\ref{thm:discovery-spectral-separation}.

The extremal-space theorem holds at every proper retained level
for all positive-field Lorentzian laws: the modes are precisely the
Hessian-null linear counts. A weighted upper-star inequality and
parallel-class decomposition give the unit-edge separation; triangles
attain the threshold three.

The same decomposition expresses all retained-level matroid walks
through contraction minors and shows that scalar equality tilts preserve
their spectra simultaneously (Corollary~\ref{thm:parallel-blocks}).
The general splitting theorem uses derivative laws, permits repeated
particles, and needs no Lorentzian hypothesis
(Theorem~\ref{thm:multinomial-blocks}).

\subsection{A general matrix Laplace theorem}
Let $L$ be the positive generator of a finite irreducible reversible
chain with invariant law $\mu$ and rates $q(x,y)$. For a Hermitian
observable $A$, define
\[
 \Gamma_L(A)(x)=\frac12\sum_yq(x,y)(A(x)-A(y))^2.
\]
Suppose that $\alpha>0$ satisfies the complete entropy-production
inequality
\[
 \alpha\mathcal H_\mu(F)
 \le \E_\mu\Tr[(LF)\log F]
 \qquad(F>0).
\]

\begin{theorem}[Matrix Laplace comparison]
\label{thm:intro-nonlinear}
If $\Gamma_L(A)(x)\preceq V$ for a fixed matrix $V\succeq0$
and every $x$, then for every Hermitian offset $K$ and every
real $\theta$,
\begin{equation}\label{eq:intro-nonlinear}
 \E_\mu\Tr e^{K+\theta A}
 \le\Tr\exp\!\left(K+\theta\E_\mu A+
                                    \frac{\theta^2}{\alpha}V\right).
\end{equation}
\end{theorem}

The transport estimate retains $\Tr(MV)$ until the matrix mean $M$
is optimized. For a centered observable,
an affine variance bound $\Gamma_L(A)\preceq aA+B$, with $a\ge0$
and $B\succeq0$, yields a Bernstein comparison with
$\theta^2B/(\alpha-a\theta)$ for $\theta<\alpha/a$; when $a=0$,
all real $\theta$ are allowed (Theorem~\ref{thm:nonlinear-affine}).
For a rank-$r$ matroid, the product-field basis law
$\mu_\lambda(B)\propto\prod_{i\in B}\lambda_i$, $\lambda_i>0$,
is Lorentzian \cite{BrandenHuh}.
If a $d\times d$ Hermitian function $A$ changes by at most $a>0$
in operator norm across every basis exchange, then, for $t>0$,
\begin{equation}\label{eq:intro-nonlinear-tail}
 \Pr\{\lambda_{\max}(A-\E A)\ge t\}
 \le d\exp\!\left(-\frac{t^2}{2ra^2}\right).
\end{equation}
Componentwise complementation replaces rank by the sum of the
smaller rank or corank, with componentwise exchange bounds.

Entropy--trace duality also gives a Laplace comparison for
subconfiguration averages (Theorem~\ref{thm:u-laplace}), yielding
matrix Chernoff bounds and Bernstein variance
$\sum_i\Pr(i\in B)A_i^2-(\E Y)^2/r$ for $Y=\sum_{i\in B}A_i$.

Kathuria explicitly proposed extending nonlinear matrix concentration
to strongly log-concave laws and arbitrary matroids, and asked for
a matrix modified logarithmic Sobolev inequality and a matrix Herbst
argument \cite[Section 5]{Kathuria}. Complete replacement and the
finite-chain comparison give this route for every Lorentzian law;
\eqref{eq:intro-nonlinear-tail} gives the proposed all-matroid conclusion.

Dependent matrix concentration was developed for strongly Rayleigh
measures \cite{R5KS,Kathuria}, under $\ell_\infty$-independence
\cite{R5KKS}, and for nonlinear observables under stochastic covering
assumptions \cite{R5AK}. Harvey--Olver proved matrix Chernoff bounds
for pipage-rounding outputs on matroid base polytopes
\cite{HarveyOlver}. The distribution in the present results is
the product-field basis law itself, for every matroid.
These estimates include the Fano basis law, whose generating polynomial
fails real stability \cite[Example~11.1]{COSW}.
Scalar-to-matrix Poincar\'e methods already give exponential
concentration for Lorentzian measures \cite{GargKathuriaSrivastava}.
Transport and concentration from complete logarithmic Sobolev
inequalities are developed in \cite{GJL2020,GJLL,DePalmaRouze};
Huang--Tropp obtain nonlinear subgaussian matrix concentration from
a curvature criterion \cite{HuangTropp}.

\subsection{Sharp basis-graph inequalities}
Complete entropy factors over exchange edges. For the uniform law
on the $N$ bases of a matroid and the unit-edge Laplacian
$(L_Gg)(B)=\sum_{B'\sim B}(g(B)-g(B'))$, we obtain
\begin{align*}
 |\partial U|&\ge a\log_2(N/a),
       &&a=\min\{|U|,N-|U|\}>0,\\
 \operatorname{Ent}_\mu(g^2)&\le\langle g,L_Gg\rangle_\mu,
       &\mathcal H_\mu(e^{-tL_G}F)&\le e^{-4t}\mathcal H_\mu(F).
\end{align*}
All three constants are sharp uniformly over matroids.
The small-set bound extends the cube's classical logarithmic
edge-isoperimetric inequality \cite{HarperCube,BernsteinCube,HartCube}.
Montenegro--Son proved this bound for balanced matroids
\cite{MontenegroSon}. Anari--Liu--Oveis Gharan--Vinzant proved the matroid case of the
Mihail--Vazirani conjecture: every matroid basis graph has
expansion at least one \cite{ALOV}.
Here equality sets are exactly cylinders fixing binary direct factors,
up to complementation, and the optimal Gross constant is one exactly
when such a factor exists. Spectral separation also bounds distance
from the binary eigenspace.

Under the full transverse gap \eqref{eq:hodge-gap}, the complete
unit-edge rate improves to $4r(1+\kappa)/(r+(r-1)\kappa)$,
with corresponding Gross and small-set improvements.
Every simple matroid of rank at least two therefore has a guaranteed
complete rate greater than four.

\subsection{The exact marginal region}
Suppose a homogeneous multiaffine generating polynomial and its
complement are both Lorentzian. The exact region of weights for which
the complete marginal inequality holds at every positive field is
\[
 \mathcal W(M)=\{w\ge0:w_i+w_j\le1
 \text{ for distinct }i,j\text{ in the same nontrivial component of }M\}.
\]
For $X_i=\mathbf1_{i\in B}$, this means
\[
 \sum_iw_i\mathcal H_{\mu_i}(\mathbb E[F\mid X_i])
       \le\mathcal H_\mu(F)
\]
in every matrix dimension. The same region is exact for the matrix
H\"older inequality required at every positive field, with an arbitrary
Hermitian offset
(Theorem~\ref{r4q:thm:weighted} and
Corollary~\ref{cor:matrix-holder-region}). It depends only on the
component partition, and includes every product-field matroid basis law.

For the half vector, fixed-field sharpness is governed by the common
translation directions of the law and its complement.
On connected rank-$m$ support, a sharp field exists exactly when
$n=2m$ and $\mathcal N_f$ contains a vector $h$ with no zero coordinate.
That translation space is then one-dimensional, and every sharp field
rescales $|h|$ by one constant on its positive half and another on its
negative half. A single translation already forces complement reciprocity
and independent Bernoulli coordinates within both halves
(Theorems~\ref{efg:reciprocity} and~\ref{efg:sharp-fields}).
The optimal coefficient of a product law is the maximum of its component
coefficients. Unless a binary direct factor exists, the half bound is
strict on an open dense set of positive fields, uniformly in matrix size
(Theorem~\ref{thm:marginal-component-coefficient}).
The sharp-field locus is unchanged for every trace-power entropy
$1<p\le2$, with the same complete strictness bound
(Theorem~\ref{thm:joint-kernel-coefficient}).

The matrix equality family has scalar parameters on those sharp components,
with noncommuting factors confined to binary direct choices
(Theorem~\ref{thm:half-matrix-form}).
For a connected rank-three sharp law, the Lorentzian Hessian condition becomes the triangle
condition for a unitary matrix, giving a complete determinantal
classification (Theorem~\ref{efg:cubic-classification}).
In every higher rank, stable matching-mixture laws and determinantal
sharp laws are incomparable (Proposition~\ref{thm:balanced-incomparability}).
The corresponding local entropy spectrum is explicit and independent
of the field parameter; its matrix statement excludes the full
eigenvalue-one space of the second variation.

\subsection{Profiles, equality, and reconstruction}
The profiles $H_k=\mathcal H_{\mu_k}(F_k)$ form exactly the cone
with $H_0=0$, nonnegative initial slope and discrete convexity.
Classical encodings on partition matroids realize the whole cone
(Appendix~\ref{app:profile-cone}); second differences decompose
every deletion deficit into nonnegative terms.

If $\E_\mu\Tr F=1$, its classical--quantum state $\rho$ admits a
Gibbs reconstruction $\tau_\ell$ at each level $1\le\ell<r$, with
normalizing mass $0<Z_\ell\le1$ and
\[
 H_r-\frac r\ell H_\ell
 =D(\rho\Vert\tau_\ell)-\log Z_\ell,
 \qquad
 \|\rho-\tau_\ell\|_1^2
 \le\frac{2r}{\ell}\left(\frac\ell rH_r-H_\ell\right).
\]
The coefficient is optimal. Support compression treats singular
inputs, and trace saturation forces self-consistency. A separate
one-particle specialization reconstructs from one conditional
matrix per active label, even when the known deficit is at a
higher retained level.

Section~\ref{sec:joint-entropy-rigidity} combines scalar translations
with the quantum normal form to classify fixed-field marginal equality.

\setcounter{tocdepth}{1}
\tableofcontents
\section{Complete Lorentzian entropy contraction}\label{sec:complete-contraction}\subsection{Matrix entropy and the three-matrix inequality}

Lieb's trace inequality converts the Lorentzian quadratic signature
into a matrix entropy bound. We fix the notation first.

For a probability measure $\pi$ on a finite set and positive semidefinite
matrices $F(x)$ of one fixed size, write
\[
 \mathcal H_\pi(F)
 =\sum_x\pi(x)\Tr[F(x)\log F(x)]
       -\Tr[\overline F\log\overline F],
 \qquad \overline F=\sum_x\pi(x)F(x).
\]
We use $0\log0=0$. No normalization of the trace of $\overline F$ is
needed. When $\Tr\overline F=1$, the classical--quantum state
$\rho_{XQ}=\sum_x\pi(x)|x\rangle\langle x|\otimes F(x)$ satisfies
\begin{equation}\label{r4q:eq:cq-identity}
 \mathcal H_\pi(F)
 =D(\rho_{XQ}\Vert\pi\otimes\rho_Q).
\end{equation}
Its quantum marginal is $\rho_Q=\overline F$. Thus the reference on the
right uses the actual quantum marginal; the classical marginal of
$\rho_{XQ}$ is allowed to differ from $\pi$.

For a positive definite matrix $\sigma$, define the Fr\'echet derivative
of the logarithm by
\[
 \mathcal T_\sigma(A)
 =D\log(\sigma)[A]
 =\int_0^\infty(\sigma+tI)^{-1}A(\sigma+tI)^{-1}\,dt.
\]
It is self-adjoint and strictly positive as an operator on the real
Hilbert space of Hermitian matrices, and
$\mathcal T_\sigma(\sigma)=I$. In particular
$\langle A,B\rangle_\sigma=\Tr[A\mathcal T_\sigma(B)]$ is a real
inner product. We use Lieb's three-matrix inequality~\cite{R4Lieb}:
\begin{equation}\label{r4q:eq:lieb}
 \Tr\exp(\log A-\log B+\log C)
 \leq\Tr[A\mathcal T_B(C)]\qquad(A,B,C>0).
\end{equation}

\subsection{Occupation laws and complete contraction}
\label{r4l:section}

The quadratic inequality propagates by conditioning on two unmarked
particles. Occupation notation includes repeated labels.

For degree $r\geq2$, a homogeneous polynomial with nonnegative
coefficients is \emph{Lorentzian} when its support is $M$-convex
and every order-$(r-2)$ partial derivative has a Hessian with at most
one positive eigenvalue. Nonnegative homogeneous polynomials of
degree at most one are Lorentzian. Here $M$-convexity means that,
for support points $\alpha,\beta$ and $\alpha_i>\beta_i$,
some $j$ with $\alpha_j<\beta_j$ has both
$\alpha-e_i+e_j$ and $\beta+e_i-e_j$ in the support.
For homogeneous nonnegative polynomials this class is equivalent to
strong log-concavity; see \cite[Theorems 2.25 and 2.30]{BrandenHuh}.

Let
\[
 f(z)=\sum_{\alpha\in\mathbb N^n,\,|\alpha|=r}a_\alpha z^\alpha
 \quad (a_\alpha\geq0),\qquad f\ne0,
\]
be Lorentzian and homogeneous of degree $r\geq1$.  For a positive
field $\lambda$, define the occupation law
\[
 \mu_\lambda(\alpha)=\frac{a_\alpha\lambda^\alpha}{f(\lambda)}.
\]
An occupation $\alpha$ describes $r$ particles, of which $\alpha_i$
have label $i$.  Retaining a uniformly chosen set of $\ell$ particles
defines the classical channel
\begin{equation}\label{r4l:downkernel}
 D_{r\to\ell}(\beta\mid\alpha)
  =\frac{\prod_i\binom{\alpha_i}{\beta_i}}{\binom r\ell},
 \qquad |\beta|=\ell,\quad\beta\leq\alpha.
\end{equation}
Write $\mu_\ell=\mu_\lambda D_{r\to\ell}$.  Every classical channel
below acts as the identity on the quantum register.

\begin{theorem}[Complete Lorentzian downsampling]
\label{r4l:main}
For every finite quantum register $Q$, every classical--quantum state
$\rho_{\mathsf{A} Q}$ whose classical support is contained in
$\operatorname{supp}\mu_\lambda$, and every $0\leq\ell\leq r$,
\begin{equation}\label{r4l:mainineq}
 D\bigl((D_{r\to\ell}\otimes\operatorname{id}_Q)\rho_{\mathsf{A} Q}
       \Vert\mu_\ell\otimes\rho_Q\bigr)
 \leq\frac\ell r
 D\bigl(\rho_{\mathsf{A} Q}\Vert\mu_\lambda\otimes\rho_Q\bigr).
\end{equation}
Equivalently, the matrix entropy of every positive semidefinite density
contracts by $\ell/r$ under conditional expectation onto the retained
occupation.  The coefficient has no dependence on $n$ or $\dim Q$.
It is the best
universal coefficient for the class of Lorentzian laws.
\end{theorem}

\begin{proof}
\emph{The retained laws.}
Absorb the field and normalization into $f$, so that $f(\one)=1$ and
its coefficients are the probabilities of $\mu_\lambda$.  Positive
diagonal substitutions preserve the Lorentzian property.  With
$D_\one=\sum_i\partial_i$, the generating polynomial of the level-$k$
law is
\begin{equation}\label{r4l:levelpoly}
 f_k=\frac{k!}{r!}D_\one^{r-k}f.
\end{equation}
Indeed the coefficient of $z^\beta$ in this polynomial is
$\sum_{\alpha\geq\beta}\mu_\lambda(\alpha)
  \prod_i\binom{\alpha_i}{\beta_i}/\binom rk$.
Every $f_k$ is normalized and Lorentzian.  The closure under
nonnegative directional derivatives used here is
\cite[Corollary 2.11]{BrandenHuh}.

For a positive semidefinite matrix function $F$ on level $r$, let
$F_k(\beta)=\mathbb E[F(\mathsf{A})\mid\mathsf{B}_k=\beta]$ under the base
law and uniform deletion.  All its means equal a fixed matrix $M$.
Put
\[
 H_k(F_k)=\mathbb E_{\mu_k}\operatorname{Tr}(F_k\log F_k)
          -\operatorname{Tr}(M\log M).
\]
We prove that $(H_k)_{k=0}^r$ is convex. Since $H_0=0$, its
endpoint chord is the required contraction inequality.

\emph{The quadratic step.}
At level two, give the two remaining particles a uniform random order
$(U,V)$.  The symmetric ordered-pair law $W$ has entries
\[
 W_{ii}=\mu_2(2e_i),\qquad
 W_{ij}=\tfrac12\mu_2(e_i+e_j)\quad(i\ne j).
\]
Thus $2W=\nabla^2 f_2$ and $W$ has at most one positive eigenvalue.
Its row sums $p_i=\sum_jW_{ij}$ form the one-particle law.
Write $F_{ij}=F_2(e_i+e_j)$. Set $G_i=\mathbb E[F_2\mid U=i]$
and $M=\sum_i p_iG_i$, and
regularize to make the matrix values positive definite.
Conditional expectation rewrites the entropy difference as
\[
 H_2(F_2)-2H_1(G)=\sum_{i,j}W_{ij}\Tr F_{ij}
       (\log F_{ij}-\log G_i-\log G_j+\log M).
\]
Apply Klein's inequality with reference
$\exp(\log G_i+\log G_j-\log M)$, then Lieb's inequality
\eqref{r4q:eq:lieb}. These two steps replace the noncommuting
logarithms by the quadratic form determined by $W$:
\begin{align*}
 H_2(F_2)-2H_1(G)
 &\geq\operatorname{Tr}M-
       \sum_{i,j}W_{ij}\operatorname{Tr}
       \exp(\log G_i+\log G_j-\log M)\\
 &\geq\operatorname{Tr}M-
       \sum_{i,j}W_{ij}\operatorname{Tr}
       G_i\,\mathcal T_M(G_j)\geq0.
\end{align*}
Here $\mathcal T_M$ is the logarithmic derivative defined above.
To justify the last inequality, the form
$\langle A,B\rangle_M=\operatorname{Tr}A\mathcal T_M(B)$ is a
positive definite real inner product on Hermitian matrices.
The one-positive-eigenvalue condition and $\one^\mathsf TW\one=1$
give $c^\mathsf TWc\leq0$ whenever $p^\mathsf Tc=0$: otherwise the
span of $\one$ and $c$ would be a two-dimensional positive subspace.
Expand $G_i-M$ in an orthonormal basis for this inner product.
Every coefficient vector is $p$-centered. All cross terms with $M$
vanish, while $\mathcal T_M(M)=I$, proving the assertion.
Positive regularization followed by continuity covers singular matrix
values.  This proves the base case with repeated labels included.

\emph{Two-particle links and convexity.}
Fix $1\le j<r$ and mark $j-1$ uniformly chosen particles in a
level-$(j+1)$ occupation. If their occupation is $\gamma$, the two
unmarked particles have normalized generating polynomial
\[
 f^\gamma(z)=\frac{\partial^\gamma f_{j+1}(z)}
                        {\partial^\gamma f_{j+1}(\one)}.
\]
It is Lorentzian. Give this link the matrix value
$F_{j+1}(\gamma+\eta)$. Its mean is $F_{j-1}(\gamma)$; after retaining
one unmarked particle, its conditional mean is $F_j(\gamma+e_i)$.
The latter assertion includes repeated labels: the identity
\begin{equation}\label{r4l:markedidentity}
 \binom{\alpha}{\gamma}(\alpha_i-\gamma_i)
   =(\gamma_i+1)\binom{\alpha}{\gamma+e_i}
\end{equation}
shows that, once $\gamma+e_i$ is fixed, remembering the mark
introduces a factor independent of the original occupation.

Write $H_2^\gamma,H_1^\gamma$ for the link entropies, relative to
their common mean. The matrix entropy chain rule gives
\[
 \E_\gamma H_2^\gamma=H_{j+1}-H_{j-1},\qquad
 \E_\gamma H_1^\gamma=H_j-H_{j-1}.
\]
The quadratic step therefore proves the exact identity
\begin{equation}\label{eq:marked-quadratic-defect}
 H_{j+1}-2H_j+H_{j-1}
     =\E_\gamma(H_2^\gamma-2H_1^\gamma)\ge0.
\end{equation}
Thus the profile is convex and $H_\ell\le(\ell/r)H_r$.

Finally, write
$\rho_{\mathsf{A} Q}=\sum_\alpha\mu_\lambda(\alpha)
 |\alpha\rangle\langle\alpha|\otimes F(\alpha)$.
Then $M=\rho_Q$ and $H_r$ is its relative entropy from
$\mu_\lambda\otimes\rho_Q$.  The output density is $F_\ell$, so
the matrix entropy inequality is exactly~\eqref{r4l:mainineq}.

\emph{Sharpness.}
Take the product polynomial
$\prod_{j=1}^r(z_{j,0}+z_{j,1})$, with its uniform product law,
and a nonconstant positive scalar density of mean one depending only
on the choice in the first pair. The retained occupation reveals that choice
with probability $\ell/r$; otherwise its conditional density is one.
Thus its entropy is exactly $\ell/r$ times the input entropy.
\end{proof}

\begin{corollary}[Complete replacement contraction]
\label{r4l:replacement}
Delete one uniformly chosen particle and then draw its replacement
from the conditional law induced by $\mu_\lambda$ and that deletion
kernel.  If $K$ is this reversible heat-bath kernel, then
\[
 D\bigl((K_*\otimes\operatorname{id}_Q)\rho
          \Vert\mu_\lambda\otimes\rho_Q\bigr)
 \leq(1-r^{-1})D(\rho\Vert\mu_\lambda\otimes\rho_Q).
\]
For $r=1$ the replacement draws a fresh particle and the left side
is zero.  With generator $I-K$, continuous time gives decay at least
$e^{-t/r}$ in this relative entropy.  Equivalently, with the convention
$\operatorname{Ent}_\mu F\leq c^{-1}
 \mathbb E_\mu\operatorname{Tr}[(I-K)F\log F]$, the complete
modified logarithmic Sobolev constant is at least $1/r$.
\end{corollary}

\begin{proof}
Apply Theorem~\ref{r4l:main} with $\ell=r-1$, then quantum data
processing for the conditional replacement channel.  Its reference
output is $\mu_\lambda$, and its quantum marginal remains $\rho_Q$.
Convexity gives
$H(F)-H(KF)\leq\mathbb E\operatorname{Tr}[(I-K)F\log F]$.
Differentiation along the semigroup and this inequality prove the
continuous-time assertion.
\end{proof}

\section{Entropy profiles and strict contraction}
\label{sec:entropy-profiles}

The marked-link identity controls every chord of the entropy profile.
Its second differences give the exact deficits and equality conditions.
Appendix~\ref{app:profile-cone} proves that their nonnegativity gives
all universal profile constraints.

Put $H_k=\mathcal H_{\mu_k}(F_k)$ and $H_0=0$.

\begin{theorem}[Convexity and the exact deficit decomposition]
\label{thm:entropy-profile}
For every positive semidefinite matrix input,
\[
 d_j:=H_{j+1}-2H_j+H_{j-1}\ge0\qquad(1\le j<r).
\]
Each $d_j$ is the average complete quadratic entropy defect after
marking $j-1$ particles in the level-$(j+1)$ reference law.
Consequently
\begin{equation}\label{eq:entropy-profile-chord}
 H_\ell-H_j\le\frac{\ell-j}{r-j}(H_r-H_j)
                  \qquad(0\le j<\ell<r),
\end{equation}
and
\begin{equation}\label{eq:profile-green-deficit}
 \frac\ell r H_r-H_\ell
 =\sum_{j=1}^{r-1}
   \frac{\min(\ell,j)\,[r-\max(\ell,j)]}{r}\,d_j.
\end{equation}
Equality in deletion at one proper retained level is equivalent to
$H_k=kH_1$ for every $k$. More generally, equality in any chord
\eqref{eq:entropy-profile-chord} is equivalent to the profile being
affine on the whole interval $\{j,\ldots,r\}$. It is also equivalent
to equality in deletion from $r-j$ to $\ell-j$ in every positive-mass
link obtained by marking $j$ particles.
\end{theorem}
\begin{proof}
The proof of Theorem~\ref{r4l:main} established convexity and its
link interpretation in \eqref{eq:marked-quadratic-defect}.
Summing the second differences gives
\[
 H_k=kH_1+\sum_{j=1}^{k-1}(k-j)d_j.
\]
Substitution proves \eqref{eq:profile-green-deficit}. Applying the
same identity to the shifted profile $(H_{j+k}-H_j)_{k=0}^{r-j}$
proves the chord bound; its deficit is
\[
 \sum_{t=j+1}^{r-1}
 \frac{[\min(\ell,t)-j]\,[r-\max(\ell,t)]}{r-j}\,d_t.
\]
All coefficients are positive. This proves both affine-profile
characterizations.

Finally, mark $j$ particles at level $r$. The chain rule expresses
$H_r-H_j$ and $H_\ell-H_j$ as the averages of the input and output
entropies in the resulting links. The required conditional means
follow from the general identity
\[
 \binom{\alpha}{\gamma}\binom{\alpha-\gamma}{\eta}
 =\binom{\alpha}{\gamma+\eta}\binom{\gamma+\eta}{\gamma},
\]
where $\gamma$ is marked and $\eta$ is retained outside the mark.
Their deletion deficits are
nonnegative and average to the chord deficit. The average vanishes
exactly when every positive-mass link has zero deficit.
\end{proof}

\begin{corollary}[Information after partial disclosure]
\label{cor:conditional-quantum-information}
Suppose the classical input marginal equals the reference law $\mu$.
For nested uniform retentions
$\mathsf A_r\to\mathsf A_\ell\to\mathsf A_j$,
\[
 I(\mathsf A_\ell:Q\mid\mathsf A_j)
 \le\frac{\ell-j}{r-j}I(\mathsf A_r:Q\mid\mathsf A_j)
 \qquad(0\le j<\ell<r).
\]
In particular, $H_j=0$ implies
$H_\ell\le(\ell-j)H_r/(r-j)$.
\end{corollary}
\begin{proof}
Here $H_k=I(\mathsf A_k:Q)$. The Markov property gives
$I(\mathsf A_k:Q\mid\mathsf A_j)=H_k-H_j$ for $k\ge j$;
apply \eqref{eq:entropy-profile-chord}.
\end{proof}

\subsection{A quadratic criterion and quantitative strictness}
\label{sec:hodge}

The quadratic signature is the exact local criterion: scalar
contraction, complete contraction, and one positive eigenvalue are
equivalent. Its strict version controls the entire deletion chain.
Only the unconditioned two-particle law needs a quantitative gap;
the conditioned links retain the qualitative Lorentzian hypothesis.
All transition eigenvalues below are real by reversibility.

\begin{theorem}[Quadratic entropy criterion]\label{thm:quadratic-criterion}
Let $W$ be a symmetric entrywise nonnegative matrix with
$\one^{\mathsf T}W\one=1$, and remove its zero rows. Put
$p=W\one$. Regard $W$ as the law of an ordered pair $(U,V)$,
and let $F_{ij}=F_{ji}\succeq0$ be a matrix function of its unordered
occupation. The following conditions are equivalent.
\begin{enumerate}
\item $W$ has exactly one positive eigenvalue.
\item For scalar $F\geq0$,
  $\mathcal H_p(\E[F\mid U])\leq\tfrac12\mathcal H_W(F)$.
\item The inequality in (2) holds in every finite matrix dimension.
\end{enumerate}
Thus a scalar signature condition is sufficient for the complete
entropy inequality, with the same constant.
\end{theorem}
\begin{proof}
The proof of Theorem~\ref{r4l:main} gives (1)$\Rightarrow$(3),
and (3)$\Rightarrow$(2) is immediate. For necessity, put
$P=\operatorname{diag}(p)^{-1}W$. This stochastic matrix is
self-adjoint on $L^2(p)$ and has spectrum in $[-1,1]$. If $W$
has a second positive eigenvalue, the restriction of $P$ to the
orthogonal complement of the constants has an eigenvector $h$
with eigenvalue $\lambda>0$. For sufficiently small real $\varepsilon$,
set $F_{ij}=1+\varepsilon(h_i+h_j)$. Its mean is one and
$\E[F\mid U=i]=1+\varepsilon(1+\lambda)h_i$. Expansion at the
constant density gives
\[
 \mathcal H_W(F)
 =\varepsilon^2(1+\lambda)\|h\|_{L^2(p)}^2+O(\varepsilon^3),
 \qquad
 2\mathcal H_p(\E[F\mid U])
 =\varepsilon^2(1+\lambda)^2\|h\|_{L^2(p)}^2+O(\varepsilon^3).
\]
These contradict (2). Finally $\one^{\mathsf T}W\one=1$ ensures
that at least one eigenvalue is positive.
\end{proof}

For a normalized degree-$r$ Lorentzian polynomial with $r\geq2$, let $W$ be the
ordered two-particle law after uniform deletion, and let $p=W\one$.
Write $\kappa\geq0$ for any number satisfying
\begin{equation}\label{eq:hodge-gap}
 h^{\mathsf T}Wh\leq-\kappa\sum_i p_i h_i^2
 \quad\text{whenever }\sum_i p_i h_i=0.
\end{equation}
Equivalently, the nonconstant spectrum of
$\operatorname{diag}(p)^{-1}W$ is bounded above by $-\kappa$.
The value $\kappa=0$ is always available. For $k\geq1$ put
\[
 a_k=k+(k-1)\kappa.
\]

\begin{theorem}[The two-particle gap controls every deletion level]
\label{thm:all-level-gap}
Under \eqref{eq:hodge-gap}, in every finite matrix dimension,
\begin{equation}\label{eq:all-level-gap}
 \mathcal H_{\mu_\ell}(F_\ell)
 \leq\frac{a_\ell}{a_r}\mathcal H_{\mu_r}(F),
 \qquad 1\leq\ell\leq r.
\end{equation}
For deletion of one particle followed by conditional replacement,
the complete modified logarithmic Sobolev constant is at least
\[
 \frac{1+\kappa}{r+(r-1)\kappa}.
\]
The theorem uses only the gap of the unmarked two-particle law;
no quantitative bound on the conditioned links is required.
\end{theorem}
\begin{proof}
We first improve the rank-two estimate. Put
$G_i=\E[F_2\mid U=i]$, $M=\sum_i p_iG_i$, and $Z_i=G_i-M$.
Regularize by adding $\varepsilon I$ if necessary. The proof of
Theorem~\ref{r4l:main}, together with \eqref{eq:hodge-gap}, gives
\[
 H_2-2H_1\geq
 -\sum_{ij}W_{ij}\langle Z_i,Z_j\rangle_M
 \geq\kappa\sum_i p_i\langle Z_i,Z_i\rangle_M,
\]
where $H_k=\mathcal H_{\mu_k}(F_k)$ and
$\langle A,B\rangle_M=\Tr A\mathcal T_M(B)$.
Operator concavity of the logarithm implies
\[
 \log A-\log M\preceq\mathcal T_M(A-M).
\]
Taking the trace against $A$, and using
$\mathcal T_M(M)=I$, yields
\[
 \Tr[A(\log A-\log M)-A+M]
 \leq\Tr[(A-M)\mathcal T_M(A-M)].
\]
Averaging over $A=G_i$ proves $H_2\geq(2+\kappa)H_1$.

The entropy-profile chord inequality with lower endpoint one gives,
for $k\geq2$,
\begin{equation}\label{eq:refined-recursion}
 H_{k+1}\geq\frac{kH_k-H_1}{k-1}.
\end{equation}
Starting with $H_2\geq a_2H_1$, induction in
\eqref{eq:refined-recursion} gives $H_k\geq a_kH_1$ for every
$k$. Substitution of this last estimate back into the same recursion
gives
\[
 H_{k+1}\geq
 \left(\frac{k}{k-1}-\frac{1}{(k-1)a_k}\right)H_k
 =\frac{a_{k+1}}{a_k}H_k.
\]
Together with the rank-two estimate, these inequalities telescope to
\eqref{eq:all-level-gap}. Conditional replacement cannot increase
matrix entropy, so its one-step coefficient is at most
$a_{r-1}/a_r$. The supporting-hyperplane inequality for
$\Tr A\log A$ then gives the asserted entropy-production constant.
\end{proof}

\begin{corollary}[Exact structural criterion for the universal coefficient]
\label{cor:matroid-sharpness}
Let $M$ be a loopless matroid of rank $r\geq2$, with a fixed positive
product field, and let $1\leq\ell<r$. Its complete deletion
coefficient equals $\ell/r$ if and only if $M$ has a parallel class
containing at least two elements. If $M$ is simple, the coefficient
is strictly smaller than $\ell/r$, uniformly over the quantum
dimension.
\end{corollary}
\begin{proof}
For a simple matroid, the Hessian of its basis polynomial at every
positive point has signature $(+,-,\ldots,-)$ by
Murai--Nagaoka--Yazawa~\cite{MNY}. After absorbing the field,
the Hessian of the two-particle polynomial is a positive scalar
multiple of that Hessian. Indeed, for a homogeneous polynomial of
degree $r$,
\[
 \nabla^2\!\left(\frac{2}{r!}D_{\one}^{r-2}f\right)
 =\frac{2}{r(r-1)}\nabla^2 f(\one).
\]
Its restriction to $p^\perp$ is negative definite, so
\eqref{eq:hodge-gap} holds with $\kappa>0$. Since $\ell<r$,
$a_\ell/a_r<\ell/r$.

Conversely, let $C$ be a nontrivial parallel class and put
$q_i=\lambda_i/\sum_{j\in C}\lambda_j$. Choose nonconstant
positive matrices $A_i$ with mean $A=\sum_{i\in C}q_iA_i$, and set
\[
 F(B)=A_i\quad\text{if }B\cap C=\{i\},\qquad
 F(B)=A\quad\text{if }B\cap C=\varnothing.
\]
Conditional on the other elements of a basis meeting $C$, its
element in $C$ has law $q$. If a retained occupation misses $C$,
the conditional mean of $F$ is therefore $A$; if it contains $i$,
the conditional mean is $A_i$. The probability of retaining the
class element is $\ell/r$ times its original probability.
The output matrix entropy is consequently exactly $\ell/r$
times the input entropy. Scalar nonconstant $A_i$ already suffice.
\end{proof}

\section{Trace duality and the characterization of Lorentzian laws}
\label{sec:duality}

Entropy contraction and exponential trace comparison are dual forms
of the same finite-channel inequality. We keep an arbitrary
Hermitian offset in the exponential. This is the form that can be
iterated under conditioning in the concentration arguments.

Let
$\mu$ be a probability measure on $\mathcal X$, let
$D(y\mid x)$ be a stochastic kernel, and put $\nu=\mu D$.
For a matrix function $F$ on $\mathcal X$, define
\[
 TF(y)=\E[F(X)\mid Y=y],\qquad
 D^*H(x)=\sum_y D(y\mid x)H(y).
\]
Zero reference probabilities are omitted. The means of $F$ and $TF$
agree.

\begin{theorem}[Complete entropy--trace duality]
\label{thm:channel-duality}
Fix a matrix dimension and a real number $c\geq1$. The following
assertions are equivalent.
\begin{enumerate}
\item For every $F\succeq0$,
  $\mathcal H_\nu(TF)\leq c^{-1}\mathcal H_\mu(F)$.
\item For every positive definite matrix family $A_y$, with
  $M=\E_\nu A$, and every Hermitian $K$,
  \begin{equation}\label{eq:channel-trace}
  \E_\mu\Tr\exp\bigl(K+cD^*(\log A)\bigr)
  \leq\Tr\exp(K+c\log M).
  \end{equation}
\item For every such family,
  \begin{equation}\label{eq:channel-normalized}
  \E_\mu\Tr\exp\bigl(cD^*(\log A)-(c-1)\log M\bigr)
  \leq\Tr M.
  \end{equation}
\end{enumerate}
The statements hold at every matrix dimension simultaneously whenever
one of them does.
\end{theorem}
\begin{proof}
The same exponential family proves both directions.
Write $D_{\rm raw}(X\Vert Y)=\Tr X(\log X-\log Y)$.
This function is jointly convex, and Klein's inequality states
$D_{\rm raw}(X\Vert Y)\geq\Tr X-\Tr Y$.
Assume (1), and set
\[
 F_x=\exp(K+cD^*(\log A)(x)),\quad
 J_y=TF(y),\quad N=\E_\mu F=\E_\nu J.
\]
Adjointness gives
\[
 \E_\mu\Tr F\log F=\Tr NK+c\E_\nu\Tr J\log A.
\]
Consequently the entropy inequality implies
\[
 c\E_\nu D_{\rm raw}(J\Vert A)
 \leq\Tr NK+(c-1)\Tr N\log N.
\]
Joint convexity bounds the left side below by
$cD_{\rm raw}(N\Vert M)$. Rearrangement and Klein's inequality
against $\exp(K+c\log M)$ prove (2). Taking
$K=-(c-1)\log M$ proves (3).

Conversely, first take positive definite $F$, put $J=TF$ and
$N=\E F$, and use (3) with $A=J$. Define
\[
 E_x=\exp\bigl(cD^*(\log J)(x)-(c-1)\log N\bigr).
\]
The averaged Klein inequality gives
\begin{align*}
 \E_\mu\Tr F\log F
 &\geq\E_\mu\Tr F\log E+\Tr N-\E_\mu\Tr E\\
 &\geq c\E_\nu\Tr J\log J-(c-1)\Tr N\log N.
\end{align*}
This is (1). Replacing $F$ by $F+\varepsilon I$ and passing to
the limit handles positive semidefinite functions.
\end{proof}

This explicit duality is a finite-channel form of quantum
Brascamp--Lieb duality~\cite{R5BSW}. For the Lorentzian deletion
kernel it gives an inequality for matrices indexed by arbitrary
subconfigurations.

\begin{corollary}[Trace inequalities at every order]
\label{cor:all-order-trace}
Let $\mu$ be a degree-$r$ Lorentzian occupation law and let
$1\leq\ell\leq r$. For positive definite matrices $A_\beta$
indexed by retained $\ell$-occupations and any Hermitian $K$,
\begin{equation}\label{eq:all-order-trace}
 \E_\mu\Tr\exp\!\left(
 K+\frac r\ell\sum_\beta D_{r\to\ell}(\beta\mid\alpha)
                                      \log A_\beta\right)
 \leq\Tr\exp\!\left(K+\frac r\ell\log\E_{\mu_\ell}A_\beta\right).
\end{equation}
The coefficient $r/\ell$ is the largest universally valid coefficient
in this trace inequality over the class of Lorentzian laws.
Under the two-particle gap hypothesis \eqref{eq:hodge-gap} of Section~\ref{sec:hodge}, it may be increased to $a_r/a_\ell$, where $a_k=k+(k-1)\kappa$.
\end{corollary}
\begin{proof}
Apply Theorem~\ref{thm:channel-duality} to
Theorems~\ref{r4l:main} and~\ref{thm:all-level-gap}.
The product-of-pairs equality example in
Theorem~\ref{r4l:main}, together with the converse direction of
the duality, proves universal optimality.
\end{proof}

For $\ell=1$, writing $p_i=\E\alpha_i/r$, the inequality becomes
\[
 \E_\mu\Tr\exp\!\left(K+\sum_i\alpha_i\log A_i\right)
 \leq\Tr\exp\!\left(K+r\log\sum_i p_iA_i\right).
\]
For $r=1$ this is Lieb's concavity inequality. The general theorem
allows arbitrary noncommuting matrices on every retained level.

The scalar trace comparison has a geometric interpretation:
$f^{1/r}$ lies below each of its tangent affine functions. This
recovers concavity and, after taking derivatives, Lorentzianity.
For multiaffine polynomials, the scalar all-field equivalence below
is the $\alpha=1$ case of the entropic-independence characterization
in~\cite[Theorem 4]{AJKPV}. We include the tangent-majorant argument
to handle general homogeneous polynomials and to identify the
equivalent complete conditions.

\begin{theorem}[Entropic characterization of Lorentzian polynomials]
\label{thm:lorentzian-characterization}
Let $f$ be a nonzero homogeneous polynomial with nonnegative
coefficients. The following are equivalent.
\begin{enumerate}
\item $f$ is Lorentzian.
\item Every nonzero partial derivative $g=\partial^\gamma f$ of
  positive degree $k$, at every positive external field, has
  scalar selected-particle entropy contraction at most $1/k$.
\item Every law in (2) has complete selected-particle contraction
  at most $1/k$.
\item Every law in (2) has the complete deletion inequalities at
  every retained level, with coefficients $\ell/k$.
\end{enumerate}
If $f$ is multiaffine of positive degree $r$, it suffices in (2)
to test $f$ itself at every positive field, with coefficient $1/r$.
The same reduction applies to the complete conditions (3) and (4).
\end{theorem}
\begin{proof}
Lorentzian polynomials are closed under partial derivatives and
positive diagonal substitutions. Theorem~\ref{r4l:main} therefore
gives (1)$\Rightarrow$(4)$\Rightarrow$(3)$\Rightarrow$(2).

Assume (2). Fix a nonzero derivative $g$ of degree $k$ and a
positive point $a$. Its field law has selected-particle distribution
$p_i=a_i\partial_i g(a)/(kg(a))$. Apply the scalar version of
Theorem~\ref{thm:channel-duality} with $A_i=b_i/a_i$. For every
$b>0$,
\[
 \frac{g(b)}{g(a)}\leq
 \left(\frac{\sum_i b_i\partial_i g(a)}{kg(a)}\right)^k.
\]
Thus $\phi=g^{1/k}$ lies below its tangent affine function at every
positive point:
\[
 \phi(b)\leq\phi(a)+\nabla\phi(a)\cdot(b-a).
\]
The function $\phi$ is the infimum of these affine majorants, so
it is concave. In particular $g$ is log-concave on the positive
orthant. This holds for every nonzero partial derivative, proving
strong log-concavity. For homogeneous polynomials with nonnegative
coefficients, strong log-concavity is equivalent to Lorentzianity
by Br\"and\'en--Huh~\cite[Theorem 2.30]{BrandenHuh}.

For multiaffine $f$, the test on $f$ alone first proves its
log-concavity by the same tangent argument. For every variable,
\[
 \partial_i f=\lim_{u\to\infty}u^{-1}f|_{z_i=u}.
\]
Positive pointwise limits preserve log-concavity. Iterating this
identity over distinct variables proves log-concavity of every
nonzero partial derivative; repeated derivatives are zero.
Thus $f$ is strongly log-concave and hence Lorentzian.
\end{proof}

\begin{corollary}[Recognition at one positive field]
\label{cor:one-field-recognition}
Fix one positive field $x$. A nonzero homogeneous polynomial $f$
with nonnegative coefficients is Lorentzian if and only if every
nonzero coordinate derivative $g$ of degree $k\ge2$, at that field,
satisfies
\[
 \operatorname{Var}(T_{k\to1}u)\le k^{-1}\operatorname{Var}(u)
\]
for every linear count observable $u(\alpha)=\alpha\cdot h$.
Thus these infinitesimal scalar tests imply all the complete
inequalities at every positive field.
\end{corollary}
\begin{proof}
The forward direction follows by linearizing complete deletion.
For the converse absorb $x$ into $g$, normalize $g(\one)=1$, and put
$W=\nabla^2g(\one)/(k(k-1))$, $p=W\one$ and
$P=\diag(p)^{-1}W$ on the active labels.
Here $T_1^*h(\alpha)=\alpha\cdot h/k$ and
$T_1T_1^*=[I+(k-1)P]/k$: after conditioning on a tagged particle,
a uniformly selected particle is the tag with probability $1/k$,
and otherwise follows the pair transition $P$.
A positive nonconstant eigenvalue of $P$ would violate the test on
the linear function $T_1^*h$.
Hence each derivative Hessian has only one positive eigenvalue.
It is indecomposable: a decomposition in disjoint active variables
would give two positive Hessian blocks.
The quadratic derivative Hessians are constant, so their
log-concavity holds everywhere. The recognition theorem of
Anari--Liu--Oveis Gharan--Vinzant
\cite[Theorem 3.2]{ALOVIII} now gives complete log-concavity, hence
Lorentzianity. Degrees zero and one are automatic.
\end{proof}

\section{Translation directions and all scalar equality cases}
\label{sec:translation-equality}

Every scalar equality density comes from a translation that leaves
the generating polynomial unchanged.  For a Lorentzian polynomial,
these global translations are exactly the null directions of its
Hessian at any positive point.  This identifies all equality cases
and decides when the deletion coefficient is sharp.

Absorb the positive reference field into the variables and normalize
the Lorentzian generating polynomial as
\[
 f(z)=\sum_{|\alpha|=r}\mu(\alpha)z^\alpha,
 \qquad f(\one)=1,\qquad r\geq2.
\]
Remove variables with $\partial_i f=0$. Put
\[
 \mathcal N_f=\{h\in\mathbb R^n:D_hf=0\},
 \qquad D_h=\sum_i h_i\partial_i.
\]
Equivalently, $h\in\mathcal N_f$ means
$f(z+th)=f(z)$ as a polynomial identity in $z$ and $t$.
The space $\mathcal N_f$ is therefore intrinsic to $f$.

\begin{lemma}[The Hessian kernel is the translation space]
\label{lem:translation-kernel}
For every positive $x$,
\[
 \ker\nabla^2 f(x)=\mathcal N_f.
\]
Consequently the Hessian rank is constant on the positive orthant.
\end{lemma}
\begin{proof}
Differentiating $D_hf=0$ gives one inclusion.  For the converse,
the Lorentzian signature forces a Hessian null vector to lie in
the kernel of every derivative Hessian.  Iteration then turns a
null vector at one point into a polynomial identity.  This is the
kernel argument of Br\"and\'en--Huh~\cite[Lemma~2.15]{BrandenHuh},
also expressed by the degree-one Hodge--Riemann theorem
\cite[Theorem~3.8]{MNY}; the proof below includes identically zero
derivatives.

Induct on the degree, starting with the constant Hessian of a
quadratic.  For degree $d\geq3$, Euler's identities give
\[
 (d-2)\nabla^2f(x)=\sum_i x_i\nabla^2(\partial_i f)(x),
 \qquad
 x^{\mathsf T}\nabla^2(\partial_i f)(x)
     =(d-2)e_i^{\mathsf T}\nabla^2 f(x).
\]
Let $h\in\ker\nabla^2f(x)$.  The second identity makes $h$
orthogonal to $x$ with respect to each derivative Hessian $H_i$:
$x^{\mathsf T}H_i h=0$. For every nonzero derivative,
\[
 x^{\mathsf T}H_i x=(d-1)(d-2)\partial_i f(x)>0.
\]
Thus $H_i$ has exactly one positive eigenvalue. Its restriction to
this orthogonal complement is
therefore nonpositive, with zero square only on its kernel.
The first identity now expresses $0$ as a sum of these
nonpositive squares with strictly positive coefficients $x_i$.
Every square vanishes, so $h$ lies in every derivative Hessian
kernel.

Induction gives $D_h\partial_i f=0$ for each nonzero derivative;
the identity holds for the zero derivatives as well.  Hence every
partial derivative of $D_hf$ vanishes.  Since $D_hf$ is homogeneous
of positive degree, it is zero.  Thus $h$ generates a global
translation of $f$.
\end{proof}

\begin{theorem}[All scalar equality densities]
\label{thm:lorentzian-equality}
Let $f$ be as above and let $1\leq\ell<r$. For a probability
law $\nu$ supported on $\operatorname{supp}\mu$, equality in
\[
 D(\nu D_{r\to\ell}\Vert\mu D_{r\to\ell})
   \leq\frac\ell r D(\nu\Vert\mu)
\]
holds if and only if there is $a\in\mathbb R_{\geq0}^n$ such that
\begin{equation}\label{eq:translation-tilt}
 a-\one\in\mathcal N_f,
 \qquad \nu(\alpha)=\mu(\alpha)a^\alpha.
\end{equation}
These densities attain equality at every retained level. Writing
$p_i=\partial_i f(\one)/r$, their entropies are
\begin{equation}\label{eq:translation-entropy}
 \mathcal H_{\mu_k}(a^\beta)
       =k\sum_i p_i a_i\log a_i\qquad(0\leq k\leq r),
\end{equation}
with $0\log0=0$ and $0^0=1$.
\end{theorem}
\begin{proof}
Write $F=\nu/\mu$ and $H_k=\mathcal H_{\mu_k}(F_k)$.
Equality at one proper retained level makes the entire entropy
profile linear, by Theorem~\ref{thm:entropy-profile}:
$H_k=kH_1$ for every $k$.  The quadratic level will identify the
translation direction, and $H_r=rH_1$ will identify the density.

Set $a_i=F_1(i)$ and $h=a-\one$. Then $p^{\mathsf T}h=0$.
The quadratic proof of Theorem~\ref{r4l:main} gives
\[
 0=H_2-2H_1\geq-h^{\mathsf T}Wh\geq0,
 \qquad W=\frac{\nabla^2 f(\one)}{r(r-1)}.
\]
The form of $W$ is nonpositive on $p^\perp$.  Polarizing at the
zero square $h^{\mathsf T}Wh=0$ shows that $Wh$ annihilates $p^\perp$, so $Wh$
is proportional to $p$.  But
$\one^{\mathsf T}Wh=p^{\mathsf T}h=0$ and
$\one^{\mathsf T}p=1$, forcing $Wh=0$.
Lemma~\ref{lem:translation-kernel} therefore gives
$h\in\mathcal N_f$.

Translation invariance normalizes the candidate extremizer:
$f(a)=f(\one)=1$, so $\mu_a(\alpha)=\mu(\alpha)a^\alpha$ is a
probability law.  It contains the support of $\nu$: whenever
$F(\alpha)>0$ and $\alpha_i>0$, that occupation contributes
positively to $F_1(i)=a_i$.  Relative entropy now gives the
remaining rigidity:
\[
 D(\nu\Vert\mu_a)
 =H_r-\sum_i\mathbb E_\nu\alpha_i\log a_i
 =H_r-r\sum_i p_i a_i\log a_i
 =H_r-rH_1=0.
\]
Vanishing relative entropy forces $\nu=\mu_a$, proving necessity.

Conversely, suppose \eqref{eq:translation-tilt} holds.
Translation invariance passes to every partial derivative of $f$.
The conditional-mean formula for deletion therefore gives, for
each retained occupation $\beta$,
\[
 \mathbb E_\mu[a^\alpha\mid\beta]
 =a^\beta\frac{\partial^\beta f(a)}
                       {\partial^\beta f(\one)}
 =a^\beta.
\]
The identity extends to $a_i=0$ by polynomial continuity.
The normalized generating polynomial $f_k$ of $\mu_k$ is a
positive multiple of $D_{\one}^{r-k}f$, so
$f_k(a)=1$ and
$\partial_i f_k(a)=\partial_i f_k(\one)=kp_i$.
Its tilted mean occupation is consequently $kp_i a_i$.
Expanding $\log a^\beta$, with the stated boundary convention,
gives \eqref{eq:translation-entropy} and equality at every level.
\end{proof}

For an arbitrary density, its one-particle marginal still selects
a canonical field tilt.  The entropy deficit splits exactly into
relative entropy from that tilt and a normalization defect.

\begin{proposition}[The deficit controls an explicit field tilt]
\label{prop:lorentzian-deficit}
For any probability law $\nu$ on $\operatorname{supp}\mu$, set
$G_i=(\nu/\mu)_1(i)$, $Z=f(G)$ and
$\mu_G(\alpha)=Z^{-1}\mu(\alpha)G^\alpha$. Then $0<Z\leq1$ and
\begin{equation}\label{eq:scalar-defect}
 D(\nu\Vert\mu)-rD(\nu D_{r\to1}\Vert\mu D_{r\to1})
     =D(\nu\Vert\mu_G)-\log Z.
\end{equation}
Writing $H_k=D(\nu D_{r\to k}\Vert\mu D_{r\to k})$, every
$1\leq\ell<r$ also satisfies
\begin{equation}\label{eq:all-level-deficit}
 \frac\ell r H_r-H_\ell
 \geq\frac{r-\ell}{r(r-1)}
             \bigl(D(\nu\Vert\mu_G)-\log Z\bigr).
\end{equation}
Both terms in parentheses are nonnegative.
\end{proposition}
\begin{proof}
Any occupation of positive $\nu$ mass has $G_i>0$ wherever
its exponent is positive, so $Z>0$.  The scalar bound in
Corollary~\ref{cor:all-order-trace}, extended to the boundary by
continuity, gives $Z\leq(\sum_i p_iG_i)^r=1$.
Expansion of $D(\nu\Vert\mu_G)$ gives
\eqref{eq:scalar-defect}.
For $\ell>1$, the chord from level one to level $r$ gives
\[
 H_\ell\le
 \frac{\ell-1}{r-1}H_r+\frac{r-\ell}{r-1}H_1.
\]
Rearranging and substituting \eqref{eq:scalar-defect} proves
\eqref{eq:all-level-deficit}; for $\ell=1$ it is an identity.
\end{proof}

\begin{corollary}[Exact sharpness for every Lorentzian law]
\label{cor:lorentzian-sharpness}
For $1\leq\ell<r$, the complete deletion coefficient of a fixed
positive-field Lorentzian law equals $\ell/r$ if and only if its
active-variable Hessian is singular. Equivalently, its nonzero first
partial derivatives are linearly dependent. In the singular case,
nonconstant positive scalar densities attain the coefficient.
In the nonsingular case, the improvement in
Theorem~\ref{thm:all-level-gap} is strict in every matrix dimension.
This dichotomy is independent of the positive reference field.
\end{corollary}
\begin{proof}
A nonsingular Hessian is negative definite on $p^\perp$.
Thus \eqref{eq:hodge-gap} holds with $\kappa>0$, and the
complete contraction coefficient is strictly below $\ell/r$.
If $0\ne h\in\mathcal N_f$, take $a=\one+th>0$ for small
nonzero $t$. Euler's identity gives $p^{\mathsf T}h=0$.
Since every $p_i>0$, strict convexity of $u\log u$ gives
$\sum_i p_i a_i\log a_i>0$. The density $a^\alpha$ is therefore nonconstant and attains the
complete upper bound.  Finally, the translation space is
independent of the reference point, while positive fields act by
invertible diagonal coordinate changes.  The dichotomy holds
throughout the field family.
\end{proof}

\begin{corollary}[Sharpness over the whole field family]
\label{cor:field-uniform-sharpness}
Let $f$ be multiaffine Lorentzian of degree $r\ge2$ with at least
two supported bases. For $1\le\ell<r$,
\[
 \sup_{\lambda>0}\eta_{\mathrm{complete}}
                    (\mu_\lambda,D_{r\to\ell})=\ell/r .
\]
Thus strictness at every fixed positive field need not give an
improvement uniform over fields.
\end{corollary}
\begin{proof}
Expose an edge of the support matroid's base polytope and send
the positive field along an exposing functional.  The laws
converge to a nondegenerate law on two bases differing by one
exchange.  On this pair, deletion retains the distinguishing
particle with probability $\ell/r$ and otherwise erases the
choice.  Relative entropy therefore contracts by exactly $\ell/r$.

Fix a nonconstant scalar tilt of this limiting two-point law and
view it as a probability law $\nu$ on the original support.
Every point of positive $\nu$ mass has positive limiting reference
mass, and the same holds after deletion. Thus the scalar relative entropies converge to their two-point
values. Their output-to-input ratio tends to $\ell/r$, giving
the lower bound on the supremum. Complete deletion gives the upper bound.
\end{proof}

\subsection{Equality for product-field matroid bases}
\label{sec:equality}

For a matroid, the translation space has a concrete meaning:
one may redistribute mass inside a parallel class while leaving the
class itself unchanged. The equality theorem below shows that these
redistributions exhaust all scalar extremizers.

Let $M$ be a loopless rank-$r$ matroid, let $\lambda_i>0$, and
write $\mu$ for its product-field basis law. Partition the ground
set into parallel classes $C$. Put
\[
 w_C=\sum_{i\in C}\lambda_i,\qquad q_i=\lambda_i/w_C.
\]
Singleton classes are included. A basis contains at most one element
from each class. On the classes it induces the product-field basis
law of the simplification $\overline M$, with fields $w_C$.
Conditioned on this simplified basis, the chosen labels in its
classes are independent and have probabilities $q_i$.

\begin{theorem}[Complete description of scalar equality]
\label{thm:matroid-equality}
Suppose $r\geq2$ and $1\leq\ell<r$. For any probability law
$\nu$ on the bases of $M$, equality in
\[
 D(\nu D_{r\to\ell}\Vert\mu D_{r\to\ell})
 \leq\frac\ell r D(\nu\Vert\mu)
\]
holds if and only if there are numbers $a_i\geq0$ such that
\begin{equation}\label{eq:parallel-tilt}
 \sum_{i\in C}q_i a_i=1\quad\text{for every parallel class }C,
 \qquad
 \nu(B)=\mu(B)\prod_{i\in B}a_i.
\end{equation}
Every such density attains equality simultaneously at all retained
levels. Thus equality changes only the distribution of the label
inside each parallel class; it preserves the law of the simplified
basis. For a simple matroid, equality forces $\nu=\mu$.
\end{theorem}

\begin{proof}
Let $L$ be the class-sum map $(Lh)_C=\sum_{i\in C}\lambda_i h_i$.
The field-absorbed basis polynomial factors through $L$, and its
Hessian is a positive scalar multiple of
$L^{\mathsf T}\nabla^2 f_{\overline M}(w)L$.
The map $L$ is onto. The Hessian of the simple matroid
$\overline M$ is nonsingular by Murai--Nagaoka--Yazawa~\cite{MNY},
so the field-absorbed Hessian has kernel exactly $\ker L$.
Lemma~\ref{lem:translation-kernel} identifies this with the
translation space of the field-absorbed polynomial.
The condition $a-\one\in\ker L$ is exactly
$\sum_{i\in C}q_i a_i=1$ for each class.
Theorem~\ref{thm:lorentzian-equality} proves the assertion,
including zero coordinates and equality at every retained level.
\end{proof}

\subsection{The full extremal eigenspace}
\label{sec:extremal-dynamics}

Translation directions determine every eigenfunction that attains the
contraction bound $\ell/r$. A one-particle Gram identity detects these
directions, and the variance chord forces the same equality space at
every retained level. This also identifies the slow part of the
deletion--replacement heat flow.

Keep a normalized Lorentzian law of degree $r\ge2$ on its active
variables. Write $T_\ell F=\E[F\mid\beta]$ and
$K_\ell=T_\ell^*T_\ell$, the operator that deletes to level $\ell$
and then conditionally refills. All spaces use the reference inner
products. Put
\[
 W=\frac{\nabla^2 f(\one)}{r(r-1)},\qquad p=W\one,\qquad
 P=\operatorname{diag}(p)^{-1}W.
\]

\begin{theorem}[The entire extremal eigenspace]
\label{thm:extremal-eigenspace}
For $1\le\ell<r$,
\[
 \ker(K_\ell-\ell I/r)
   =\{\alpha\mapsto\alpha\cdot h:h\in\mathcal N_f\}.
\]
Its dimension is the active Hessian nullity.
If $0<\kappa\le1/(r-1)$ bounds the strictly negative nonconstant
spectrum of $P$ above by $-\kappa$, then, on the mean-zero
orthogonal complement of these linear modes,
\[
 K_\ell\preceq\frac{\ell-(r-\ell)\kappa}{r}\,I.
\]
If there is no strictly negative nonconstant eigenvalue, take
$\kappa=1/(r-1)$.
\end{theorem}
Here $\kappa$ bounds only the strictly negative spectrum; zero modes
are removed explicitly. In contrast, \eqref{eq:hodge-gap} bounds the
whole nonconstant spectrum and permits $\kappa>0$ only when that
kernel is zero.
\begin{proof}
Condition on one tagged particle, then select a particle uniformly
from the full configuration. With probability $1/r$ this selects
the tag; otherwise the label transition is $P$. Thus
\begin{equation}\label{eq:one-particle-gram}
 T_1T_1^*=\frac{I+(r-1)P}{r}.
\end{equation}
The nonconstant spectrum of $P$ is nonpositive, and
$\ker P=\mathcal N_f$ by Lemma~\ref{lem:translation-kernel}.
The Gram identity therefore identifies the $1/r$-eigenspace of
$K_1$ with the adjoint image of $\mathcal N_f$. Since
$T_1^*h(\alpha)=\alpha\cdot h/r$, these are the linear counts
$g_h(\alpha)=\alpha\cdot h$. Their means vanish. The
factorial-second-moment formula and $Wh=0$ give
\[
 \operatorname{Var}(g_h)=r\sum_i p_i h_i^2.
\]
All active $p_i$ are positive, so the map $h\mapsto g_h$ is injective.

To transfer this rigidity to level $\ell$, linearize
\eqref{eq:entropy-profile-chord} at the constant density. For
$V_j=\operatorname{Var}(T_jF)$, the result is
\[
 V_\ell\le\frac{\ell-1}{r-1}V_r+
                   \frac{r-\ell}{r-1}V_1.
\]
Combining this with $V_1\le V_r/r$ gives $V_\ell\le\ell V_r/r$.
Equality at level $\ell<r$ forces equality in the one-particle bound,
so every extremal mode belongs to the linear space already identified.
Conversely, translation invariance gives
$T_\ell g_h(\beta)=\beta\cdot h$, and uniform retention gives
$T_\ell^*(\beta\cdot h)=\ell g_h/r$. This proves the eigenspace
identity at every level.

On the mean-zero orthogonal complement of these modes,
\eqref{eq:one-particle-gram} improves the one-particle bound to
$V_1\le[1-(r-1)\kappa]V_r/r$; the remaining kernel of $K_1$
contributes only eigenvalue zero. Substitution in the same variance
chord proves the residual estimate.
\end{proof}

\begin{corollary}[Down-up equality; Alaniz Pintos's degree-one conjecture]
\label{cor:pintos-down-up}
Let $M$ be a loopless rank-$r$ matroid with at least two bases,
$r\ge1$, and let $\mu_\lambda$ be its positive product-field basis
law. Let $P_{\mathrm{du}}$ delete a uniformly chosen element and
then refill conditionally according to $\mu_\lambda$.
Writing $C$ for the parallel classes, one has
\begin{equation}\label{eq:matroid-down-up-full-eigenspace}
 \ker\bigl(P_{\mathrm{du}}-(1-1/r)I\bigr)
 =\left\{B\mapsto\sum_{i\in B}h_i:
        \sum_{i\in C}\lambda_i h_i=0\text{ for every }C\right\}.
\end{equation}
Its dimension is $\sum_C(|C|-1)$. The down-up gap equals $1/r$
if and only if $M$ has a nontrivial parallel class; otherwise it
is strictly larger. In particular every extremal eigenfunction is
linear, resolving \cite[Conjecture 8.1]{Pintos2026} for the uniform
law and extending the conclusion to all positive fields.
\end{corollary}
\begin{proof}
For $r\ge2$, one has $P_{\mathrm{du}}=K_{r-1}$. The proof of
Theorem~\ref{thm:matroid-equality} identifies the translation space
with the kernel of the weighted class-sum map. Substitution in
Theorem~\ref{thm:extremal-eigenspace} proves the displayed identity
and its dimension. The same theorem gives the gap bound and makes
equality equivalent to a nonzero translation direction. If $r=1$,
all active elements form one parallel class and each step resamples
the entire basis. Thus $P_{\mathrm{du}}$ is projection onto constants;
its zero eigenspace is precisely the weighted mean-zero functions
on the singleton bases, which are the displayed linear functions.
\end{proof}

The exact parallel-class criterion gives the classification sought
in \cite[Conjecture 8.2]{Pintos2026}, on the original matroid. Having
a parallel pair in some minor is insufficient: $U_{2,3}$ contracts
to $U_{1,2}$, but its down-up gap is $3/4>1/2$. Indeed, its
three-state transition matrix has diagonal entries $1/2$ and
off-diagonal entries $1/4$.

The gap characterization also follows by combining
the product bound of Alev--Lau~\cite[Theorem 1.5]{AlevLau}
with the strict Hessian theorem of Murai--Nagaoka--Yazawa~\cite{MNY}.
Theorem~\ref{thm:extremal-eigenspace} identifies every extremal mode
and extends the conclusion to general Lorentzian occupation laws.

The residual estimate also separates the slow part of the continuous
flow. Let $\Pi_{\mathcal N}$ be the orthogonal projection onto the
linear modes in Theorem~\ref{thm:extremal-eigenspace}. For $t\ge0$
and every mean-zero $g$,
\[
 \left\|e^{t(K_{r-1}-I)}g-e^{-t/r}\Pi_{\mathcal N}g\right\|_2
 \le e^{-(1+\kappa)t/r}\|g-\Pi_{\mathcal N}g\|_2.
\]
Indeed, the positive generator $I-K_{r-1}$ equals $1/r$ on the projected
space and is
at least $(1+\kappa)/r$ on its orthogonal complement; apply its
spectral decomposition.

\section{Retained-level reconstruction and sharp stability}
\label{sec:quantum-reconstruction}

Average the retained conditional logarithms and correct for their
common mean. The resulting Gibbs candidate expresses the deletion
deficit as relative entropy plus a normalization defect, giving sharp
stability and self-consistency even for singular inputs.

Let $f$ be Lorentzian of degree $r\ge2$, absorb the positive field,
normalize $f(\one)=1$, and remove inactive variables. Write
\[
 \rho=\sum_\alpha\mu(\alpha)|\alpha\rangle\langle\alpha|
                                      \otimes F_\alpha,
 \qquad \sigma=\E_\mu F,\qquad\Tr\sigma=1,
\]
and put $G_i=F_1(i)$ and $p_i=\mu_1(i)$. The exponential variational
identities below belong to the framework of quantum entropy remainders
and Brascamp--Lieb duality \cite{CarlenLieb2014,R5BSW}. Here all their
matrix parameters are the observed conditional means.

For positive definite conditional matrices, reconstruction at level
$\ell$ with reciprocal contraction coefficient $c$ has the form
\[
 E_\alpha=\exp\!\left(
 c\sum_\beta D_{r\to\ell}(\beta\mid\alpha)\log F_\ell(\beta)
 -(c-1)\log\sigma\right).
\]
For semidefinite inputs, the same formula is evaluated on the common
support of the conditional matrices seen from $\alpha$.
Compress the register to $\operatorname{supp}\sigma$, so that
$\sigma>0$. For $A\succeq0$, write $\log_s A$ for its logarithm on
$\operatorname{supp}A$, extended by zero on the kernel. For a
projection $P$, $\exp_P$ denotes the exponential on $\operatorname{ran}P$;
its extension to the whole register is zero on $\ker P$.

\begin{theorem}[Reconstruction from a retained level]
\label{thm:retained-level-reconstruction}
Fix $1\le\ell<r$, and let $c\ge1$ be a valid reciprocal complete
contraction coefficient, meaning $H_r(F)\ge cH_\ell(F_\ell)$
for every matrix input.
For each occupation $\alpha$, let $P_\alpha$ project onto
\[
 \bigcap_{\beta:\,D_{r\to\ell}(\beta\mid\alpha)>0}
                         \operatorname{supp}F_\ell(\beta),
\]
and define
\begin{align*}
 K^{(\ell,c)}_\alpha
 &=c\sum_\beta D_{r\to\ell}(\beta\mid\alpha)\log_s F_\ell(\beta)
                                      -(c-1)\log\sigma,\\
 E^{(\ell,c)}_\alpha
 &=P_\alpha\exp_{P_\alpha}
       (P_\alpha K^{(\ell,c)}_\alpha P_\alpha)P_\alpha,\\
 Z_{\ell,c}&=\E_\mu\Tr E^{(\ell,c)}_\alpha,\qquad
 \tau_{\ell,c}=Z_{\ell,c}^{-1}
 \sum_\alpha\mu(\alpha)|\alpha\rangle\langle\alpha|
                                      \otimes E^{(\ell,c)}_\alpha .
\end{align*}
Then $0<Z_{\ell,c}\le1$ and
\begin{equation}\label{eq:retained-level-gibbs-deficit}
 H_r-cH_\ell
 =D(\rho\Vert\tau_{\ell,c})-\log Z_{\ell,c}.
\end{equation}
Consequently
\[
 \|\rho-\tau_{\ell,c}\|_1
 \le\sqrt{2(H_r-cH_\ell)},\qquad
 Z_{\ell,c}\ge e^{-(H_r-cH_\ell)}.
\]
Equality $H_r=cH_\ell$ holds exactly when
$Z_{\ell,c}=1$ and $\rho=\tau_{\ell,c}$.
For positive definite inputs, every $P_\alpha=I$ and the formula
is the ordinary exponential of $K^{(\ell,c)}_\alpha$.

For the universal choice $c=r/\ell$, writing
$\delta_\ell=\ell H_r/r-H_\ell$, this gives
\begin{equation}\label{eq:retained-level-gibbs-pinsker}
 \|\rho-\tau_{\ell,r/\ell}\|_1
 \le\sqrt{\frac{2r}{\ell}\,\delta_\ell}.
\end{equation}
Under \eqref{eq:hodge-gap}, one may instead use
$c=a_r/a_\ell$ and the improved deficit $H_r-(a_r/a_\ell)H_\ell$.

Conversely, start with any positive semidefinite matrices $G_\beta$ of mean
$\sigma=\E_{\mu_\ell}G$, with $\Tr\sigma=1$, and form the same
exponentials with $F_\ell(\beta)$ replaced by $G_\beta$.
If $Z_{\ell,c}=1$, then
\[
 \E_\mu E^{(\ell,c)}=\sigma,\qquad
 T_\ell E^{(\ell,c)}(\beta)=G_\beta.
\]
Thus trace saturation is sufficient for self-consistency.
\end{theorem}
\begin{proof}
The conditional mean $F_\ell(\beta)$ is a positive mixture containing
$F_\alpha$ whenever $D_{r\to\ell}(\beta\mid\alpha)>0$. Therefore
$\operatorname{supp}F_\alpha\subseteq\operatorname{ran}P_\alpha$.
At least one $F_\alpha$ is nonzero, so $Z_{\ell,c}>0$.

For a Hermitian $H$ and a projection $P$,
\[
 \Tr_P\exp_P(PHP)\le\Tr\exp H.
\]
Indeed, in an eigenbasis of $PHP$, scalar Jensen gives
$e^{\langle u,Hu\rangle}\le\langle u,e^Hu\rangle$; summing proves
the inequality. Apply it to the exponent obtained by replacing
every $F_\ell(\beta)$ with $F_\ell(\beta)+\varepsilon I$ and
$\sigma$ with $\sigma+\varepsilon I$. The normalized trace inequality
of Theorem~\ref{thm:channel-duality} bounds the sum of the
uncompressed exponential traces by $1+\varepsilon\dim Q$.
On $\operatorname{ran}P_\alpha$, the compressed logarithms converge
to those defining $K^{(\ell,c)}_\alpha$. Passing to the limit
therefore proves $Z_{\ell,c}\le1$.

The support inclusion also gives
$\Tr(F_\alpha\log E^{(\ell,c)}_\alpha)
=\Tr(F_\alpha K^{(\ell,c)}_\alpha)$.
Adjointness of deletion now gives
\[
 \E_\mu\Tr(F_\alpha\log E^{(\ell,c)}_\alpha)
 =c\,\E_{\mu_\ell}\Tr(F_\ell\log F_\ell)
             -(c-1)\Tr(\sigma\log\sigma).
\]
Expanding the relative entropy proves
\eqref{eq:retained-level-gibbs-deficit}. Its two terms on the right
are nonnegative. Quantum Pinsker and their separate bounds prove
the stability, normalization, and equality assertions. The choices
of $c$ follow from Theorems~\ref{r4l:main} and~\ref{thm:all-level-gap}.

For the converse abbreviate $E=E^{(\ell,c)}$, put
$N=\E_\mu E$ and $J=T_\ell E$, and use
$D_g(A\Vert B)=\Tr[A(\log A-\log B)-A+B]$.
By construction $\operatorname{supp}J_\beta
\subseteq\operatorname{supp}G_\beta$, so these relative entropies
are finite.
Since $\Tr N=Z_{\ell,c}=1$, the same expansion gives
\[
 H_r(E)-cH_\ell(J)
 =-c\,\E_{\mu_\ell}D_g(J_\beta\Vert G_\beta)
                         +(c-1)D_g(N\Vert\sigma)
 \le-D_g(N\Vert\sigma).
\]
The last inequality is joint convexity. Complete contraction makes
the left side nonnegative, so $N=\sigma$. The identity then forces
$D_g(J_\beta\Vert G_\beta)=0$ for every positive-mass output,
which gives $J_\beta=G_\beta$.
\end{proof}

\begin{corollary}[Sharpness of retained-level stability]
\label{cor:retained-level-stability-sharpness}
For every $r\ge2$ and $1\le\ell<r$, the coefficient $2r/\ell$
in \eqref{eq:retained-level-gibbs-pinsker} is optimal uniformly over
Lorentzian laws, even when the classical marginal is fixed to the
reference and the encoding register is a classical bit.
\end{corollary}
\begin{proof}
Take the uniform basis law of a partition matroid with $r$ binary
classes, identified with independent unbiased bits
$X=(X_1,\ldots,X_r)$. Put $\chi(X)=(-1)^{X_1+\cdots+X_r}$ and,
for $0<t<1$, encode $X=x$ by
\[
 F_x=\frac12
 \begin{pmatrix}1+t\chi(x)&0\\0&1-t\chi(x)\end{pmatrix}.
\]
The mean is $\sigma=I/2$. Every proper retention omits at least one
independent bit, so $F_\ell=\sigma$ and $H_\ell=0$ for $\ell<r$.
Its reconstruction is exactly
$\tau_{\ell,r/\ell}=\mu\otimes\sigma$. Therefore
\[
 \|\rho-\tau_{\ell,r/\ell}\|_1=t,\qquad
 H_r=\frac{(1+t)\log(1+t)+(1-t)\log(1-t)}2
     =\frac{t^2}{2}+O(t^4).
\]
Since $\delta_\ell=\ell H_r/r$, the ratio of the squared trace
distance to $\delta_\ell$ tends to $2r/\ell$ as $t\downarrow0$.
\end{proof}

Reconstruction at level $\ell$ uses the whole conditional family
$F_\ell$. The next theorem uses only the one-particle family
$(G_i)_i$, even when the given deficit is measured at a higher level.

\begin{theorem}[Gibbs reconstruction and stability]
\label{thm:quantum-gibbs-reconstruction}
For a positive definite input define
\[
 E_\alpha=\exp\!\left(\sum_i\alpha_i\log G_i-(r-1)\log\sigma\right),
 \quad Z=\E_\mu\Tr E_\alpha,\quad
 \tau=Z^{-1}\sum_\alpha\mu(\alpha)|\alpha\rangle\langle\alpha|
                                                    \otimes E_\alpha.
\]
For a general input use the support-compressed version of this
formula from Theorem~\ref{thm:retained-level-reconstruction};
its common support at $\alpha$ is
$\bigcap_{i:\,\alpha_i>0}\operatorname{supp}G_i$.
This is the $\ell=1$, $c=r$ case of
Theorem~\ref{thm:retained-level-reconstruction}, so
\begin{equation}\label{eq:quantum-gibbs-deficit}
 H_r-rH_1=D(\rho\Vert\tau)-\log Z,
 \qquad 0<Z\le1.
\end{equation}
Its saturation and self-consistency conclusions apply.
For $1\le\ell<r$, let $\delta_\ell=\ell H_r/r-H_\ell$. Then
\begin{align}
 D(\rho\Vert\tau)-\log Z
 &\le\frac{r(r-1)}{r-\ell}\delta_\ell,
 \label{eq:quantum-gibbs-all-level}\\
 \|\rho-\tau\|_1
 &\le\sqrt{\frac{2r(r-1)}{r-\ell}\delta_\ell}.
 \label{eq:quantum-gibbs-pinsker}
\end{align}
Equality at any proper deletion level is equivalent to $Z=1$ and
$\rho=\tau$, and every equality state satisfies
\begin{equation}\label{eq:quantum-null-directions}
 (G_i-\sigma)_i\in\mathcal N_f\otimes\operatorname{Herm}(Q).
\end{equation}
\end{theorem}
\begin{proof}
Take $\ell=1$ and $c=r$ in
Theorem~\ref{thm:retained-level-reconstruction}. Since
$rD_{r\to1}(i\mid\alpha)=\alpha_i$, its reconstruction is exactly
$E_\alpha$, and it gives \eqref{eq:quantum-gibbs-deficit} and the
self-consistency assertion. The profile chord bound yields
\[
 \delta_\ell\ge
 \frac{r-\ell}{r(r-1)}(H_r-rH_1).
\]
Both terms on the right of \eqref{eq:quantum-gibbs-deficit} are
nonnegative. Quantum Pinsker proves the norm estimate, and
Theorem~\ref{thm:entropy-profile} gives the equality equivalences.

For the null-direction assertion, we record all three quadratic
remainders. Let $F_{2,ij}=F_2(e_i+e_j)$ and
$W=\nabla^2f(\one)/(r(r-1))$. Let $C_{ij}$ be the support-compressed
exponential of $\log_sG_i+\log_sG_j-\log\sigma$ on
$\operatorname{supp}G_i\cap\operatorname{supp}G_j$.
Trace compression followed by Lieb's inequality, with the same
regularization as above, gives
$\Tr C_{ij}\le\Tr[G_i\mathcal T_\sigma(G_j)]$.
Putting $Z_i=G_i-\sigma$, direct expansion gives
\begin{align}
 H_2-2H_1
 &=\sum_{ij}W_{ij}D_g(F_{2,ij}\Vert C_{ij})\notag\\
 &\quad+\sum_{ij}W_{ij}
       \bigl(\Tr[G_i\mathcal T_\sigma(G_j)]-\Tr C_{ij}\bigr)\notag\\
 &\quad-\sum_{ij}W_{ij}
                  \Tr[Z_i\mathcal T_\sigma(Z_j)].
 \label{eq:quadratic-entropy-remainders}
\end{align}
Only pairs with $W_{ij}>0$ are included. The first two sums are
nonnegative, and the last term is nonnegative by the Lorentzian
signature and $\sum_i p_iZ_i=0$. At equality it follows that
\[
 0\ge-\sum_{ij}W_{ij}\Tr[Z_i\mathcal T_\sigma(Z_j)]\ge0.
\]
Expand the $Z_i$ in a Hermitian basis orthonormal for
$\langle A,B\rangle_\sigma=\Tr[A\mathcal T_\sigma(B)]$.
Every coefficient vector is $p$-centered and has zero quadratic form
under $W$. The Lorentzian signature puts it in $\ker W=\mathcal N_f$,
proving \eqref{eq:quantum-null-directions}.
\end{proof}

The common-support convention makes both reconstruction identities
exact for singular inputs as well. It also keeps every term in
\eqref{eq:quadratic-entropy-remainders} finite: a positive-mass
$F_{2,ij}$ is supported in both $G_i$ and $G_j$.

\section{Quadratic models of quantum entropy equality}
\label{sec:universal-equality}

The reconstruction theorem makes an equality density unique once its
one-particle matrices are known. We now decide exactly which matrices
can occur. Only their two-particle compatibility matters, in every
degree and every finite matrix dimension.

Normalize $f(\one)=1$, remove inactive labels, and suppose its degree
is $r\ge2$. Put $W=\nabla^2f(\one)/(r(r-1))$ and $p=W\one$.
Write $\mathcal E_f$ for the positive semidefinite densities with
$\Tr\E_\mu F=1$ attaining equality at a proper deletion level.
By Theorem~\ref{thm:entropy-profile}, these are exactly the densities
with $H_k=kH_1$ at every level.

\Needspace{24\baselineskip}
\begin{samepage}
\begin{theorem}[Quadratic extension and the complete equality model]
\label{thm:quadratic-extension}
\label{thm:two-particle-universality}
Taking the two-particle conditional density gives a bijection
\[
 T_2:\mathcal E_f\xrightarrow{\ \sim\ }\mathcal E_{f_2},
 \qquad f_2(z)=z^{\mathsf T}Wz.
\]
More explicitly, let $G_i\succeq0$, put $\sigma=\sum_i p_iG_i$, and
assume $\Tr\sigma=1$. On $\supp\sigma$ define
\[
 B_i=\sigma^{-1/2}G_i\sigma^{-1/2},\qquad
 L_i=\sigma^{-1/2}G_i^{1/2}.
\]
These are the one-particle matrices of an equality density if and
only if
\begin{equation}\label{eq:local-equality-test}
 \sum_jW_{ij}G_j=p_i\sigma\quad\text{for every }i,
 \qquad [B_i,L_j]=0\quad\text{whenever }W_{ij}>0.
\end{equation}
Supported self-pairs are included. The unique density is
\begin{equation}\label{eq:quadratic-extension-product}
 F_\alpha=\sigma^{1/2}\left(\prod_iB_i^{\alpha_i}\right)\sigma^{1/2}.
\end{equation}
The factors in every supported occupation commute, so the product
is unambiguous and positive. For a quadratic input $E_{ij}$, its
extension uses $G_i=p_i^{-1}\sum_jW_{ij}E_{ij}$ and satisfies
$H_r(F)/r=H_2(E)/2=H_1(G)$.

Consequently any two homogeneous Lorentzian laws with the same
normalized pair matrix $W$ have canonically identified equality
families, preserving $\sigma$, all $G_i$, and entropy per particle.
Their degrees may differ. Here a quadratic density is written
$E_{ij}=E(e_i+e_j)$, so $E_{ij}=E_{ji}$.
\end{theorem}
\end{samepage}

The distinction from reconstruction is existence: every quadratic
extremizer extends, and \eqref{eq:local-equality-test} is a finite
test for arbitrary proposed marginals. We prove it by identifying
where noncommuting matrices can occur.

Define an equivalence relation by declaring distinct labels $i,j$
equivalent when $W_{ij}=0$. Denote its classes by $\mathcal C$, and set
\[
 p_C=\sum_{i\in C}p_i,\qquad q_C(i)=p_i/p_C.
\]
The following lemma proves, in particular, that the declared relation is
indeed transitive.

\begin{lemma}[The missing-pair geometry of a Lorentzian law]
\label{lem:lorentzian-missing-pair-geometry}
If $W_{ij}=0$ for active distinct labels $i,j$, then
\[
 W_{ii}=W_{jj}=0,
 \qquad W\left(\frac{e_i}{p_i}-\frac{e_j}{p_j}\right)=0.
\]
Consequently every class containing more than one label is mutually
exclusive: a supported occupation contains at most one particle from
that entire class. There is a normalized homogeneous Lorentzian
polynomial $g$ such that
\begin{equation}\label{eq:universal-exclusive-splitting}
 f(z)=g((y_C)_{C\in\mathcal C}),
 \qquad y_C=\sum_{i\in C}q_C(i)z_i.
\end{equation}
Its Hessian has strictly positive entries between every two distinct
classes. Each class is either exclusive, meaning that $g$ has degree at
most one in $y_C$, or a singleton whose label can occur repeatedly.
\end{lemma}

\begin{proof}
Put $h=e_i/p_i-e_j/p_j$. It lies in $p^\perp$, where the form
$W$ is negative semidefinite. Since $W_{ij}=0$,
\[
 0\ge h^{\mathsf T}Wh
       =\frac{W_{ii}}{p_i^2}+\frac{W_{jj}}{p_j^2}\ge0.
\]
Thus both diagonal entries vanish and $h$ is null for the restricted
form. Polarization makes it $W$-orthogonal to $p^\perp$. It is also
$W$-orthogonal to $\mathbf1$, since $W\mathbf1=p$, so $Wh=0$.
The rows $W_{i\bullet}/p_i$ and $W_{j\bullet}/p_j$ are therefore equal.
This immediately proves transitivity of the missing-pair relation.

Nonnegative coefficients show that a vanishing diagonal prohibits a
repetition, while $W_{ij}=0$ prohibits simultaneous occurrence.
Lemma~\ref{lem:translation-kernel} upgrades each displayed kernel
vector to
\[
 \frac1{p_i}\partial_i f=\frac1{p_j}\partial_j f.
\]
Thus $f$ depends on variables inside each class only through the stated
weighted sum. One obtains $g$ by setting all variables in a class equal;
closure under nonnegative linear substitution makes $g$ Lorentzian.
The remaining claims follow directly from the definition of the classes.
\end{proof}

Write $\mathcal E\subseteq\mathcal C$ for the exclusive classes.
For an occupation $\alpha$, let $\gamma_C=\sum_{i\in C}\alpha_i$.
For $C\in\mathcal E$ with $\gamma_C=1$, write $x_C$ for the selected
representative. The reference law is the quotient occupation law
$\nu(\gamma)=[y^\gamma]g(y)$ followed by independent representative
choices $q_C$ in occupied exclusive classes.

We first derive the row and commutator conditions
\eqref{eq:local-equality-test} from equality. For the converse we start
with arbitrary matrices satisfying those conditions. Their class means
need not agree: the relative class densities first determine central
blocks, preserved by the faithful mean state. Scalar normalization
inside each block then gives ensembles with a common mean, to which
the Markov gluing lemma applies. Translation invariance verifies the
candidate's conditional means and entropy profile; Gibbs reconstruction
gives uniqueness.

\subsection{Recovery on supported pairs}

Vanishing of the quadratic remainders gives equality in relative-entropy
data processing. A faithful invariant state then turns recovery into
commutation.

\begin{lemma}[Petz compatibility on every supported pair]
\label{lem:lorentzian-petz-pairs}
Let $f$ be any normalized Lorentzian occupation law of degree $r\ge2$.
For a matrix density attaining equality in a proper deletion inequality,
compress to the support of $\sigma=\mathbb EF$, write $G_i=T_1F(i)$,
and set $B_i=\sigma^{-1/2}G_i\sigma^{-1/2}$ and
$L_i=\sigma^{-1/2}G_i^{1/2}$.
Let $W$ be its ordered two-particle law and $p=W\mathbf1$.
Then, whenever $W_{ij}>0$,
\begin{align}
 F_2(e_i+e_j)
   &=G_j^{1/2}\sigma^{-1/2}G_i\sigma^{-1/2}G_j^{1/2},
          \label{eq:lorentzian-petz-pair-density}\\
 [B_i,L_j]&=0.\label{eq:lorentzian-petz-pair-commutator}
\end{align}
In particular $B_i$ and $B_j$ commute on every supported pair.
The statement includes $i=j$ whenever repeated labels are supported.
\end{lemma}

\begin{proof}
Equality at a proper retained level makes the entropy profile linear
by Theorem~\ref{thm:entropy-profile},
so $H_2-2H_1=0$. The three nonnegative terms in
\eqref{eq:quadratic-entropy-remainders} give
\[
 (G_i-\sigma)_i\in\mathcal N_f\otimes\operatorname{Herm}(Q)
\]
and identify $F_{ij}:=F_2(e_i+e_j)$ with the support-compressed
exponential of $\log G_i+\log G_j-\log\sigma$.
Thus
\begin{equation}\label{eq:pair-log-identity-trace}
 \operatorname{Tr}F_{ij}\log F_{ij}
 =\operatorname{Tr}F_{ij}(\log_sG_i+\log_sG_j-\log\sigma).
\end{equation}
The support of $F_{ij}$ is contained in the supports of both $G_i$
and $G_j$.

Put $P_{ij}=W_{ij}/p_i$. Translation nullity gives
$\sum_jP_{ij}G_j=\sigma$, while conditional expectation gives
$\sum_jP_{ij}F_{ij}=G_i$. For
\[
 D_g(A\Vert B)=\operatorname{Tr}[A(\log A-\log B)-A+B],
\]
expanding \eqref{eq:pair-log-identity-trace} therefore yields the exact
joint-convexity equality
\begin{equation}\label{eq:rowwise-petz-equality}
 \sum_jP_{ij}D_g(F_{ij}\Vert G_j)=D_g(G_i\Vert\sigma).
\end{equation}
If $G_i=0$, all $F_{ij}$ on this row vanish and both conclusions are
immediate. Otherwise normalize by $\operatorname{Tr}G_i$ and apply
the equality case of relative-entropy data processing~\cite[Theorem~3]{HJPW} to the channel
that forgets the classical index $j$. The Petz recovery map gives
\eqref{eq:lorentzian-petz-pair-density}. This remains valid for singular
$G_j$ because the recovered reference mean $\sigma$ is faithful.

Consequently $B_i$ is fixed by the unital completely positive map
\[
 \Phi_i(X)=\sum_jP_{ij}L_jXL_j^*.
\]
Indeed, the recovery identity averaged over $j$ gives $\Phi_i(B_i)=B_i$;
$\sum_jP_{ij}G_j=\sigma$ gives unitality; and
\[
 \operatorname{Tr}[\sigma\Phi_i(X)]=\operatorname{Tr}(\sigma X)
\]
shows that $\sigma$ is a faithful invariant state.
Kadison--Schwarz and invariance imply
$\Phi_i(B_i^2)-B_i^2=0$. Equivalently,
\[
 \sum_jP_{ij}(B_iL_j-L_jB_i)(B_iL_j-L_jB_i)^*=0.
\]
Every term is positive semidefinite, proving
\eqref{eq:lorentzian-petz-pair-commutator}.
Since $B_j=L_jL_j^*$, the two relative densities commute as well.
\end{proof}

\subsection{Pairwise compatibility gives a global quantum Markov state}

\begin{lemma}[A common Markov extension of the class ensembles]
\label{lem:global-markov-gluing}
Suppose a family of ensembles $\{q_C(i),G_i\}_{i\in C}$ has common
faithful mean $\sigma$ and satisfies
$[B_i,L_j]=0$ whenever $i,j$ belong to distinct classes. For each set $J$ of classes,
put
\[
 F_J(x_J)=\sigma^{1/2}\prod_{C\in J}B_{x_C}\sigma^{1/2},
 \quad
 \Omega_{X_JQ}=\sum_{x_J}\left(\prod_{C\in J}q_C(x_C)\right)
       |x_J\rangle\langle x_J|\otimes F_J(x_J).
\]
These are consistent states with marginal $\sigma$. For every
$C\notin J$,
\[
 I(X_C:X_J\mid Q)_{\Omega}=0.
\]
Thus the full state satisfies
\[
 \sum_C S(X_C\mid Q)_\Omega-S(X_{\mathcal C}\mid Q)_\Omega=0,
\]
and admits a decomposition
\[
 Q=\bigoplus_z\bigotimes_C Q_C^z,\qquad
 \Omega_{X_{\mathcal C}Q}=\bigoplus_z w_z
                        \bigotimes_C\Omega_{X_CQ_C^z}^{(z)}.
\]
\end{lemma}

\begin{proof}
The commutator condition implies pairwise commutation of relative
densities from distinct classes, so every displayed product is positive.
The relations $\sum_iq_C(i)B_i=I$ show consistency and normalization.
For each class define the Petz channel
\[
 \mathcal R_C(X)=\sum_{i\in C}q_C(i)|i\rangle\langle i|\otimes
       G_i^{1/2}\sigma^{-1/2}X\sigma^{-1/2}G_i^{1/2}.
\]
It is trace preserving because $\sum_iq_C(i)G_i=\sigma$.
If $C\notin J$, the product $\prod_{D\in J}B_{x_D}$ commutes with
$L_i$ and $L_i^*$ for every $i\in C$. It follows directly that
\[
 (\operatorname{id}_{X_J}\otimes\mathcal R_C)(\Omega_{X_JQ})
      =\Omega_{X_{J\cup\{C\}}Q},
 \qquad \mathcal R_C(\sigma)=\Omega_{X_CQ}.
\]
The same channel therefore recovers both the joint state and the
reference state in relative-entropy data processing under
$\operatorname{Tr}_{X_C}$. Equality gives
$I(X_C:X_J\mid Q)=0$. The chain rule proves vanishing of the full
conditional total correlation. In particular,
\[
 \operatorname{TC}(X_1,\ldots,X_m\mid Q)
 =I(X_1:X_2\cdots X_m\mid Q)
      +\operatorname{TC}(X_2,\ldots,X_m\mid Q),
\]
where $\operatorname{TC}(X_1,\ldots,X_m\mid Q)
=\sum_jS(X_j\mid Q)-S(X_1\cdots X_m\mid Q)$.

Apply the Hayden--Jozsa--Petz--Winter equality theorem~\cite[Theorem~6]{HJPW} to
$X_{C_1}$, $Q$, and all remaining classical registers together.
It splits $Q$ into a direct sum of tensor products, separating
$X_{C_1}$ from the remaining registers inside every block.
The residual conditional total correlation is an average of the
conditional total correlations in those blocks, and each is nonnegative.
It therefore vanishes in every positive-weight block. Repeat this
argument for the remaining classes. This yields
\[
 \Omega_{X_{\mathcal C}Q}
   =\bigoplus_zp_z\bigotimes_C\Omega_{X_CQ_C^z}^{(z)},
\]
with all factors normalized. Dephasing the classical registers preserves
the state, so each factor is cq. Write
$\Omega_{X_CQ_C^z}^{(z)}=
 \sum_iq_C(i)|i\rangle\langle i|\otimes A_i^{C,z}$;
this defines positive semidefinite $A_i^{C,z}$ and gives
$\sum_iq_C(i)A_i^{C,z}=\sigma_C^z$.
\end{proof}

\subsection{An elementary invariant-state lemma}

\begin{lemma}[A faithful invariant state respects the Kraus center]
\label{lem:faithful-state-kraus-center}
Let
\[
 \Phi(X)=\sum_i t_i L_iXL_i^*,\qquad t_i>0,
\]
be a unital completely positive map with a faithful invariant state
$\sigma$, so $\operatorname{Tr}(\sigma\Phi(X))=\operatorname{Tr}(\sigma X)$.
Let $\mathcal A=C^*(L_i:i)$. Then $\sigma$ commutes with every central
element of $\mathcal A$.
\end{lemma}

\begin{proof}
Kadison--Schwarz, faithfulness, and the commutator sum used in
Lemma~\ref{lem:lorentzian-petz-pairs} show that the fixed-point space of
$\Phi$ is precisely the commutant $\mathcal A'$. Let
\[
 \mathsf E(X)=\lim_{N\to\infty}\frac1N\sum_{k=0}^{N-1}\Phi^k(X)
\]
be the finite-dimensional mean ergodic projection onto that space.
It preserves the state $\sigma$ and is an $\mathcal A'$-bimodule map,
because each Kraus operator commutes with $\mathcal A'$.
If $P$ is a central projection of $\mathcal A$, then
$P\in\mathcal A'$ and $P$ commutes with every value of $\mathsf E$.
Therefore
\[
 \mathsf E(PX(I-P))=P\mathsf E(X)(I-P)=0.
\]
State invariance gives $\operatorname{Tr}(\sigma PX(I-P))=0$ for every
$X$, proving that $\sigma$ has no off-diagonal block between $P$ and
$I-P$. Spectral decomposition proves the assertion for the whole center.
\end{proof}

\subsection{The operator normal form}

\begin{theorem}[Universal operator structure of Lorentzian entropy equality]
\label{thm:universal-lorentzian-matrix-equality}
Let $F_\alpha\succeq0$, let
$\sigma=\mathbb E_\mu F$ have trace one, and compress to its support.
Equality
\[
 \mathcal H_{\mu_\ell}(T_\ell F)
       =\frac\ell r\mathcal H_\mu(F)
 \qquad\text{for some }0<\ell<r
\]
holds if and only if the following structure exists.

There are an orthogonal decomposition
\[
 Q=\bigoplus_z\left(Q_0^z\otimes
                  \bigotimes_{C\in\mathcal E}Q_C^z\right),
\]
positive weights $w_z$ summing to one, faithful density matrices
$\tau_0^z,\tau_C^z$, scalar vectors $a^z\in[0,\infty)^{\mathcal C}$
with
\begin{equation}\label{eq:universal-skeleton-translation-tilt}
 a^z-\mathbf1\in\mathcal N_g,
\end{equation}
and positive semidefinite class matrices $A_i^{C,z}$ satisfying
\[
 \sum_{i\in C}q_C(i)A_i^{C,z}=\tau_C^z,
\]
such that
\begin{equation}\label{eq:universal-quantum-equality-form}
 F_\alpha=
 \bigoplus_z w_z(a^z)^{\gamma(\alpha)}\,
 \tau_0^z\otimes\bigotimes_{C\in\mathcal E}
 \begin{cases}
 A_{x_C}^{C,z},&\gamma_C=1,\\
 \tau_C^z,&\gamma_C=0.
 \end{cases}
\end{equation}
The conventions are $0^0=1$ and that an empty tensor product is scalar.
The mean is
\[
 \sigma=\bigoplus_z w_z\tau_0^z\otimes
                         \bigotimes_{C\in\mathcal E}\tau_C^z.
\]
Every such density has equality at every retained level.
\end{theorem}

\begin{proof}[Proof of Theorem~\ref{thm:universal-lorentzian-matrix-equality}]
We prove at the same time the sufficiency of
\eqref{eq:local-equality-test}. Start either with an equality density
$F$, or directly with matrices $G_i$ satisfying that test. In the
first case Lemma~\ref{lem:lorentzian-petz-pairs} and
\eqref{eq:quantum-null-directions} give the test. In both cases put
$B_i=\sigma^{-1/2}G_i\sigma^{-1/2}$, and
$L_i=\sigma^{-1/2}G_i^{1/2}$. Set
\[
 \Gamma_C=\sum_{i\in C}q_C(i)G_i,
 \qquad \overline B_C=\sigma^{-1/2}\Gamma_C\sigma^{-1/2}
                      =\sum_{i\in C}q_C(i)L_iL_i^*.
\]
For distinct classes, every pair is supported, so
$[B_i,L_j]=0$ whenever $i,j$ belong to different classes.
Hence $\overline B_C$ commutes with every $L_j$ outside $C$.

Put $P_{ij}=W_{ij}/p_i$. Let $\overline W$ and
$\overline P_{CD}=\overline W_{CD}/p_C$ be the
ordered pair law and conditional kernel of $g$. The splitting identity
gives $P_{ij}=\overline P_{CD}q_D(j)$ for $i\in C,j\in D$.
Translation nullity and the common mean imply
\[
 \sum_D\overline P_{CD}\Gamma_D=\sigma,
 \qquad \sum_Dp_D\Gamma_D=\sigma.
\]
For an exclusive class $C$, one has $\overline P_{CC}=0$, and therefore
\[
 p_C\Gamma_C
   =\sum_{D\ne C}(\overline P_{CD}-p_D)\Gamma_D.
\]
After conjugating by $\sigma^{-1/2}$, the right side commutes with
every $L_i$ inside $C$. Thus $\overline B_C$ commutes with these operators
too. For a repeated singleton class, the supported self-pair in
\eqref{eq:local-equality-test} gives the same commutation directly.
Consequently every $\overline B_C$ lies in the center of
$\mathcal A=C^*(L_i:i)$.

The map $\Phi(X)=\sum_i p_iL_iXL_i^*$ is unital and has faithful
invariant state $\sigma$. Lemma~\ref{lem:faithful-state-kraus-center}
therefore makes all the $\overline B_C$ commute with $\sigma$.
Take their joint spectral decomposition. On a block $z$, write
\[
 \overline B_C=a_C^z I,
 \qquad \sigma_z=w_z\tau_z,
 \qquad \operatorname{Tr}\tau_z=1.
\]
All $G_i$ preserve these blocks, and
$\Gamma_C|_z=a_C^z\sigma_z$. The averaged translation-null identity
descends through \eqref{eq:universal-exclusive-splitting} to
\[
 (\Gamma_C-\sigma)_C\in\mathcal N_g\otimes\operatorname{Herm}(Q).
\]
Restriction to $z$ and faithfulness of $\sigma_z$ give
$a^z-\mathbf1\in\mathcal N_g$.

If $a_C^z>0$, define the normalized ensemble on block $z$ by
\[
 H_i^z=\frac{G_i|_z}{w_za_C^z},\qquad i\in C.
\]
Its $q_C$-mean is $\tau_z$, and the cross-class Petz commutators are
unchanged after this scalar normalization. If $a_C^z=0$, positivity
forces all $G_i|_z=0$ in that class; use the harmless constant ensemble
$H_i^z=\tau_z$. Repeated singleton classes also have constant ensembles.
Lemma~\ref{lem:global-markov-gluing}, applied to the exclusive classes,
gives a common quantum Markov extension with marginal $\tau_z$.
Its iterated strong-subadditivity decomposition, refining the index $z$
if necessary, has the tensor factors and ensembles in
\eqref{eq:universal-quantum-equality-form}. A spectator factor $Q_0^z$
allows the case of no exclusive classes.

It remains to prove that the proposed form is an equality density with
the original one-particle matrices. This will also prove sufficiency.
For clarity, fix one block, suppress $z$, and let $R_J(x_J)$ denote its
common quantum Markov family on a set $J$ of exclusive classes, with
mean $\tau$. The candidate block density is
\[
 \widetilde F_\alpha=a^{\gamma(\alpha)}R_{J(\alpha)}(x_{J(\alpha)}),
 \qquad J(\alpha)=\{C\in\mathcal E:\gamma_C=1\}.
\]
Translation invariance gives, for every multiindex $\delta$,
\begin{equation}\label{eq:all-derivatives-skeleton-tilt}
 \partial^\delta g(a)=\partial^\delta g(\mathbf1).
\end{equation}
In particular $g(a)=1$, and the tilted quotient law
$\nu^a(\gamma)=\nu(\gamma)a^\gamma$ is a probability law.
Fix a supported retained occupation $\beta$, with quotient counts
$\delta$. Its posterior quotient law weights $\gamma\ge\delta$ by the
falling factorial $(\gamma)_\delta=\prod_C\gamma_C!/(\gamma_C-\delta_C)!$.
Consequently
\[
 \mathbb E_\nu[a^\gamma\mid\delta]
 =a^\delta\frac{\partial^\delta g(a)}{\partial^\delta g(\mathbf1)}
 =a^\delta.
\]
Every unretained representative ensemble averages to its class mean,
whether that class is absent or present but unretained. The retained
class factors remain fixed. Combining these averages gives
\begin{equation}\label{eq:universal-candidate-conditional-mean}
 T_k\widetilde F(\beta)
      =a^{\delta(\beta)}R_{J(\beta)}(x_{J(\beta)}).
\end{equation}
This identity includes zero coordinates of $a$ by polynomial continuity.
At level zero it gives mean $\tau$, and at level one it gives
$a_C H_i$, exactly the original one-particle data after restoring $w_z$.

Let
\[
 e_C=\sum_{i\in C}q_C(i)
       \operatorname{Tr}[H_i(\log_s H_i-\log\tau)]
       \qquad(C\in\mathcal E).
\]
For fixed $\delta$, the representative average of $R_{J(\beta)}$ is
$\tau$, with trace one. The scalar factor therefore contributes
$\operatorname{Ent}_{\nu_k}(a^\delta)$. Additivity of the common Markov
family gives $\sum_{C\in J(\beta)}e_C$ for its matrix part. This part is
weighted by $\nu_k(\delta)a^\delta$, the retained law of $\nu^a$.
Thus
\begin{equation}\label{eq:universal-entropy-additivity}
 H_k(\widetilde F)
   =\operatorname{Ent}_{\nu_k}(a^\delta)
       +\sum_{C\in\mathcal E}
          \Pr_{(\nu^a)_k}(C\text{ retained})e_C.
\end{equation}
Here the tilted retained law equals the retention of $\nu^a$ by
\eqref{eq:all-derivatives-skeleton-tilt}. The scalar equality theorem
scales the first term by $k/r$. Exclusive classes contain at most one
particle, so their retention probabilities in the second term also
scale by $k/r$. Thus $H_k=(k/r)H_r$ in every block and hence in their
weighted orthogonal sum.

The candidate is therefore an equality density with the prescribed mean $\sigma$ and one-particle matrices $G_i$.
This proves sufficiency in \eqref{eq:local-equality-test}. If an
equality density $F$ was given, Theorem~\ref{thm:quantum-gibbs-reconstruction}
makes it unique, so $F=\widetilde F$. The same calculation proves
sufficiency of every displayed normal-form parameter.
\end{proof}

\begin{proof}[Completion of the proof of Theorem~\ref{thm:quadratic-extension}]
Necessity of \eqref{eq:local-equality-test} follows from
\eqref{eq:quantum-null-directions} and
Lemma~\ref{lem:lorentzian-petz-pairs}; the normal-form construction
just proved its sufficiency. In each block the relative one-particle
matrix is the scalar $a_C$ times its class likelihood factor. Their
product is exactly \eqref{eq:quadratic-extension-product}.

The pair matrix $W$ determines $p$, the exclusive classes, their laws
$q_C$, and the quotient pair matrix
$\overline W_{CD}=\sum_{i\in C,j\in D}W_{ij}$. For every law with
this pair matrix, Lemma~\ref{lem:translation-kernel} gives
$\mathcal N_g=\ker\overline W$. Thus the admissible normal-form
parameters are the same, and their product formula has the same
one-particle means. Uniqueness from Gibbs reconstruction makes this
correspondence independent of the chosen block decomposition.
Taking the second law to be $f_2$ proves that the correspondence is
$T_2$ and is bijective. Equality gives $H_r/r=H_1$ in every law.
\end{proof}

\begin{corollary}[Matroid equality and parallel choices]
\label{cor:matroid-quantum-parallel}
Let $\mu$ be a positive product-field basis law of a loopless matroid
of rank $r\ge2$. For every normalized matrix entropy equality density,
the normal form in Theorem~\ref{thm:universal-lorentzian-matrix-equality}
has $a^z=\mathbf1$ in every block.
If $\mathsf I(B)$ is the basis of the simplification selected by $B$, then
\[
 \E_\mu[F_B\mid\mathsf I]=\sigma,
 \qquad \rho_{\mathsf I Q}=\mu_{\mathsf I}\otimes\sigma.
\]
Thus equality retains no information about the simplified basis:
all nonconstant factors describe choices within parallel classes.
This holds even when the input classical marginal differs from $\mu$.
If the matroid is simple, every equality density is constant.
\end{corollary}
\begin{proof}
The simplified positive-field basis polynomial has nonsingular Hessian,
as used in Theorem~\ref{thm:matroid-equality}. Hence
$\mathcal N_g=0$, so every quotient tilt is $\mathbf1$.
Conditioned on $\mathsf I$, the representatives are independent with
laws $q_C$. Averaging each class factor in
\eqref{eq:universal-quantum-equality-form} gives $\tau_C^z$,
and therefore the same matrix $\sigma$ for every simplified basis.
For a simple matroid, each class has one representative and its
factor is already constant.
\end{proof}

\begin{corollary}[Quantum equality is classical when no alternatives are exclusive]
\label{cor:simple-support-classical-equality}
If every two distinct active labels can occur together, every matrix
entropy equality family is simultaneously diagonalizable. More
explicitly, there are a unitary $U$, positive weights $s_a$ summing to
one, and vectors $v^a\ge0$ with $v^a-\mathbf1\in\mathcal N_f$ such that
\[
 F_\alpha=U\operatorname{diag}
       \bigl(s_a(v^a)^\alpha\bigr)_aU^*.
\]
This applies in particular to every multiaffine Lorentzian law whose
support matroid is simple, and to every multinomial law.
\end{corollary}

\begin{proof}
Every class is a singleton. Its ensemble is constant, so the tensor
part of the universal formula is independent of $\alpha$. Diagonalize
the fixed faithful matrix in each orthogonal block and absorb its
eigenvalues into the scalar weights.
\end{proof}

\subsection{Fixed classical marginals remove every quotient tilt}

\begin{corollary}[Fixed-marginal rigidity and exact attainment]
\label{cor:fixed-marginal-no-quotient-tilt}
Let $F$ be an equality density in
Theorem~\ref{thm:universal-lorentzian-matrix-equality}, and suppose
$\operatorname{Tr}F_\alpha=1$ for every supported occupation.
Then every scalar vector in its universal decomposition satisfies
$a^z=\mathbf1$.
Thus the quotient occupation is independent of the quantum register,
and within each classical block all nonconstant tensor factors
belong to exclusive classes.

An informative trace-one equality encoding exists in a finite
register if and only if two distinct active labels never co-occur,
equivalently $W_{ij}=0$ for some $i\ne j$. A qubit suffices.
If every exclusive class is a singleton, the only equality encoding
is constant: $F_\alpha=\sigma$ for every supported occupation.
If, in addition, the active Hessian is singular, the fixed-marginal
Holevo contraction coefficient is nevertheless $\ell/r$, and is
approached by nonconstant encodings but never attained by one.
\end{corollary}

\begin{proof}
All conditional densities, including the one- and two-particle
densities, have trace one. Average the one-particle identity over
representatives inside a class. In the universal decomposition this gives
\[
 \sum_z w_z a_C^z=1.
\]
For every two distinct classes $C,D$, their joint pair probability is
positive. Average the corresponding two-particle identity over the
independent representative laws. The common Markov extension has mean
one in trace, so
\[
 \sum_z w_z a_C^za_D^z=1\qquad(C\ne D).
\]
Translation nullity also gives the pointwise identity
$\sum_C p_Ca_C^z=1$ in every block. Regard $z$ as a random index with
law $w_z$. The preceding pair identities say that all distinct
$a_C,a_D$ have zero covariance. Therefore
\[
 0=\operatorname{Var}\!\left(\sum_Cp_Ca_C\right)
   =\sum_Cp_C^2\operatorname{Var}(a_C).
\]
Since every $p_C>0$, all $a_C$ are constant and hence equal to one.
The same conclusion is immediate if there is only one class.

When every exclusive class is a singleton, its ensemble has only one
member and is constant. The universal formula with $a=\mathbf1$ then
gives $F_\alpha=\sigma$.

Conversely, let $C$ contain at least two labels. Choose positive definite
trace-one qubit states $A_i$, not all equal, with mean
$\tau=\sum_{i\in C}q_C(i)A_i$.
Use $A_{x_C}$ when $C$ is occupied and $\tau$ otherwise. This is the
normal form with one block, $a=\mathbf1$, and only the $C$ factor
nonconstant, so it attains equality at every level. The class has
positive mass and every $q_C(i)>0$, hence its information is positive.

For sharpness under a singular Hessian, choose a nonzero
$h\in\mathcal N_f$ and set $u(\alpha)=\alpha\cdot h$. On a qubit put
\[
 F_\varepsilon(\alpha)=\frac12
 \begin{pmatrix}
  1+\varepsilon u(\alpha)&0\\
  0&1-\varepsilon u(\alpha)
 \end{pmatrix}.
\]
For small nonzero $\varepsilon$ these are positive, have trace one,
and have common mean $I/2$. The second-order entropy expansion and
the exact variance mode $u$ make their contraction ratio tend to
$\ell/r$. Their information is nonzero, whereas the preceding
classification excludes a nonconstant extremizer.
\end{proof}

\begin{remark}[The pointwise trace constraint]
Theorem~\ref{thm:quadratic-extension} normalizes $\Tr\E F=1$.
The extra condition $\Tr F_\alpha=1$ need not transfer between
different supports with the same $W$. For example, double each
element of the rank-three Fano matroid and of $U_{3,7}$ into a parallel
pair, choosing representatives uniformly. Their uniform quotient
laws have the same ordered pair matrix, $\overline W_{CD}=1/42$
for $C\ne D$. Use a classical register of seven bits, uniform
subject to even parity on one Fano line, and representative
likelihoods $2\mathbf1_{\{\text{register bit}=\text{chosen bit}\}}$.
Every Fano basis omits a bit of that line, so its product density has
trace one. On the uniform-matroid basis equal to the line, the trace
is zero or two according to parity. Both are complete-entropy
extremizers under \eqref{eq:quadratic-extension-product}; only the
Fano family is pointwise trace one on its whole support.
\end{remark}

\section{Power entropies and orthogonal equality structure}
\label{sec:power-entropies}

Let $f(z)=\sum_{|\alpha|=r}\mu(\alpha)z^\alpha$ be normalized,
homogeneous, and Lorentzian, with $r\ge2$. Let $T_\ell$ be conditional
expectation after uniform retention of $\ell$ particles. For a positive
semidefinite matrix-valued function $F$, put
\[
 \mathcal H_p(F)=\operatorname{Tr}\mathbb E F^p
                   -\operatorname{Tr}(\mathbb EF)^p,
 \qquad 1<p\le2.
\]
For any trace function $\varphi$ used below, write more generally
$\mathcal H_\varphi(F)=\Tr\E\varphi(F)-\Tr\varphi(\E F)$.
In particular $\mathcal H_{\log}$ denotes $\varphi(x)=x\log x$.
No trace normalization is imposed. The logarithmic equality results
apply after positive scalar rescaling; a zero mean forces $F=0$.

\Needspace{11\baselineskip}
\begin{samepage}
\begin{theorem}[Complete power-entropy contraction]\label{thm:complete-power-contraction}
For every matrix dimension, every $1<p\le2$, and $0\le\ell\le r$,
\[
 \mathcal H_p(T_\ell F)\le\frac\ell r\mathcal H_p(F).
\]
The complete coefficient of a fixed positive-field law equals
$\ell/r$, for $0<\ell<r$, if and only if its active Hessian is singular.
\end{theorem}
\end{samepage}

\begin{proof}
For $s>0$ set
\[
 h_s(x)=(x+s)\log(x+s)-s\log s-(1+\log s)x.
\]
For $1<p<2$ one has the positive integral representation
\[
 x^p=\frac{p(p-1)\sin(\pi(p-1))}{\pi}
             \int_0^\infty s^{p-2}h_s(x)\,ds.
\]
Indeed, both sides and their first derivatives vanish at zero, while
their second derivatives agree by
$\int_0^\infty s^{p-2}/(x+s)\,ds
 =\pi x^{p-2}/\sin(\pi(p-1))$.
The integral converges at zero and infinity, also for $x=0$.
Functional calculus and finite expectation give
\[
 \mathcal H_p(F)=c_p\int_0^\infty
 s^{p-2}\mathcal H_{\log}(F+sI)\,ds,
\]
where $c_p$ is the positive coefficient above and
$\mathcal H_{\log}(A)=\operatorname{Tr}\mathbb E A\log A
 -\operatorname{Tr}(\mathbb EA)\log(\mathbb EA)$.
Since $T_\ell(F+sI)=T_\ell F+sI$, integration of complete entropy
contraction proves the assertion. The case $p=2$ is its variance
linearization, or follows by continuity as $p\uparrow2$.

If the Hessian is nonsingular, Theorem~\ref{thm:all-level-gap}
improves entropy contraction uniformly in matrix dimension. Apply
it to every shift and integrate. If the Hessian is singular, choose nonzero
$h\in\mathcal N_f$ and $g(\alpha)=\alpha\cdot h$.
For $F_\varepsilon=1+\varepsilon g>0$,
\[
 \mathcal H_p(F_\varepsilon)
 =\frac{p(p-1)}2\varepsilon^2\mathbb E g^2+O(\varepsilon^3).
\]
The variance identity in Theorem~\ref{thm:extremal-eigenspace} shows that the contraction ratio
tends to $\ell/r$. Scalar functions already prove complete sharpness.
\end{proof}

\subsection{Tensor and orthogonal equality}

Take the partition matroid with two classes of three elements and
uniform representative laws. Choose positive definite matrices
$A_1,A_2,A_3$ on $\mathbb C^{d_A}$ and $B_1,B_2,B_3$ on
$\mathbb C^{d_B}$ with averages $I_{d_A}$ and $I_{d_B}$.
For instance, in dimension two one may take
\[
 A_1=I+tX,\quad A_2=I+tZ,\quad A_3=I-t(X+Z),
 \qquad 0<t<1/\sqrt2,
\]
where $X,Z$ are the Pauli matrices, and use the same choice for $B$.
The matrices in each class do not commute.

For the tensor-product density $F_{ij}=A_i\otimes B_j$, logarithmic
entropy satisfies $\mathcal H_{\log}(T_1F)=\mathcal H_{\log}(F)/2$.
For powers put
\[
 a_p=\operatorname{Tr}\mathbb EA_i^p-d_A,
 \qquad b_p=\operatorname{Tr}\mathbb EB_j^p-d_B.
\]
Then
\[
 \mathcal H_p(F)=d_Ba_p+d_Ab_p+a_pb_p,
 \qquad
 \mathcal H_p(T_1F)=\tfrac12(d_Ba_p+d_Ab_p),
\]
so the power deficit is exactly $a_pb_p/2>0$ for nonconstant ensembles.

For the orthogonal direct-sum density $\widetilde F_{ij}=A_i\oplus B_j$,
one has
\[
 \mathcal H_p(\widetilde F)=a_p+b_p,
 \qquad
 \mathcal H_p(T_1\widetilde F)=\tfrac12(a_p+b_p)
\]
for every $1<p\le2$, and the logarithmic equality holds as well.
The tensor product preserves logarithmic equality; the orthogonal sum
preserves every trace profile. The next theorem identifies all power
extremizers and explains this distinction.

\subsection{Orthogonal structure of all power extremizers}

Put $\sigma=\E_\mu F$, $G_i=(T_1F)(i)$ and $F_2=T_2F$,
and restrict all operators to $\supp\sigma$.

\begin{theorem}[Universal power equality]
\label{thm:universal-matrix-power-equality}
Let $f$ be any normalized homogeneous Lorentzian law of degree $r\ge2$,
and let $\mathcal E$ be its exclusive classes from
Lemma~\ref{lem:lorentzian-missing-pair-geometry}.
For $1<p<2$ and $0<\ell<r$, a matrix density $F$ attains equality in
\[
 \mathcal H_p(T_\ell F)\le\frac\ell r\mathcal H_p(F)
\]
if and only if there is an orthogonal decomposition into
$\sigma$-reducing subspaces
\[
 Q=Q_0\oplus\bigoplus_{C\in\mathcal E}Q_C,
 \qquad \sigma=\sigma_0\oplus\bigoplus_{C\in\mathcal E}\sigma_C,
\]
and centered Hermitian ensembles $D_i^C$ satisfying
\[
 \sum_{i\in C}q_C(i)D_i^C=0,\qquad
 \sigma_C+D_i^C\succeq0,
\]
such that
\begin{equation}\label{eq:universal-power-class-blocks}
 F_\alpha=\sigma_0\oplus\bigoplus_{C\in\mathcal E}
       \left(\sigma_C+
          \mathbf1_{\{\gamma_C(\alpha)=1\}}D_{x_C}^C\right).
\end{equation}
Zero-dimensional factors are permitted. These densities have exact
linear profiles for every trace function $\varphi$ for which the
quantities are finite.

In particular, if every exclusive class is a singleton, there are no
nonconstant matrix power-entropy extremizers. This conclusion holds in
every finite matrix dimension, even when the sharp complete coefficient
equals $\ell/r$ because the Hessian is singular.
\end{theorem}

\begin{proof}
Power equality makes $F+sI$ a logarithmic equality density for every
$s>0$. Write $D_i=G_i-\sigma$, $S_s=\sigma+sI$, and
$B_i(s)=I+S_s^{-1/2}D_iS_s^{-1/2}$.
For every supported ordered pair, the Petz identity and its commutator
in Lemma~\ref{lem:lorentzian-petz-pairs} give the exact formula
\begin{align*}
 F_{2,ij}+sI
 &=S_s^{1/2}B_i(s)B_j(s)S_s^{1/2}\\
 &=S_s+D_i+D_j+D_iS_s^{-1}D_j.
\end{align*}
Letting $s\to\infty$ yields $F_{2,ij}=\sigma+D_i+D_j$, so
\begin{equation}\label{eq:all-supported-pair-resolvent-annihilation}
 D_i(\sigma+sI)^{-1}D_j=0
 \quad\text{for every $s>0$ whenever $W_{ij}>0$}.
\end{equation}
If a label can repeat, its self-pair is supported and this identity
forces $D_i=0$.

For each exclusive class form the cyclic subspace
\[
 Q_C=\operatorname{span}
       \{\sigma^aD_iv:i\in C,\ a\ge0,\ v\in Q\}.
\]
Distinct classes have supported cross pairs. Expanding
\eqref{eq:all-supported-pair-resolvent-annihilation} at infinity and
using
\[
 \langle \sigma^aD_i v,\sigma^bD_jw\rangle
   =\langle v,D_i\sigma^{a+b}D_jw\rangle=0
\]
shows that the $Q_C$
are pairwise orthogonal and reduce $\sigma$, and that every $D_i$ in
class $C$ is supported on $Q_C$. The identity
\[
 0=\sum_i p_iD_i
   =\sum_{C\in\mathcal E}p_C\sum_{i\in C}q_C(i)D_i
\]
is an orthogonal block identity, hence every class mean vanishes.
Positivity of $G_i=\sigma+D_i$ gives the required block positivity.

Construct the candidate $\widetilde F$ by
\eqref{eq:universal-power-class-blocks}. Its class-centered form and the
reference splitting \eqref{eq:universal-exclusive-splitting} give mean
$\sigma$ and the same one-particle matrices $G_i$ as the original
density. Conditioning on a retained configuration keeps precisely the
perturbations in its occupied exclusive classes; unrevealed
representatives average to zero. Orthogonal-sum functional calculus
therefore gives
\[
 \mathcal H_\varphi(T_k\widetilde F)
 =\sum_{C\in\mathcal E}\Pr_{\mu_k}(C\text{ retained})
   \left\{\sum_{i\in C}q_C(i)
       \operatorname{Tr}\varphi(\sigma_C+D_i^C)
       -\operatorname{Tr}\varphi(\sigma_C)\right\}.
\]
Exclusive classes occur at most once, so each probability scales by
$k/r$. This proves exact linearity for every trace function.

The original $F$ is a logarithmic equality density by letting
$s\downarrow0$ in the shifted equalities. The candidate is also such
a density and has the same mean and one-particle data. Exact Gibbs
reconstruction therefore gives $F=\widetilde F$. Conversely, the
displayed profile formula proves power equality for every candidate.
If an exclusive class is a singleton, its centered perturbation is
zero, proving the final assertion.
\end{proof}

\subsection{An exact fourth-order criterion}
\label{sec:quartic}

Let $f$ be a normalized homogeneous polynomial with nonnegative
coefficients, of degree
$r\ge2$, with occupation law $\mu$, and put
\[
 W=\frac{\nabla^2f(\mathbf1)}{r(r-1)},\qquad
 p=W\mathbf1=\frac1r\nabla f(\mathbf1).
\]
Let $D_i$ be Hermitian $d\times d$ matrices satisfying the
matrix-valued translation identity
\begin{equation}\label{eq:matrix-direction-translation-identity}
 \sum_i D_i\partial_i f=0.
\end{equation}
Thus every matrix-entry vector of $(D_i)_i$ is a polynomial
translation direction. Define
\[
 H_\alpha=\sum_i\alpha_iD_i,\qquad F_\varepsilon=I+\varepsilon H.
\]
Then $\mathbb EH=0$, and $F_\varepsilon$ is positive definite for all
sufficiently small real $\varepsilon$.

\begin{theorem}[The fourth-order obstruction is sufficient as well as necessary]
\label{thm:complete-quartic-matrix-obstruction}
Set
\begin{align}
 \mathcal Q_f(D)
 &=\sum_{ij}W_{ij}
       \left(2\operatorname{Tr}D_i^2D_j^2
                    +\operatorname{Tr}D_iD_jD_iD_j\right)
       \label{eq:matrix-quartic-obstruction}\\
 &=\sum_{ij}W_{ij}
       \left(\|D_iD_j\|_{\mathrm{HS}}^2
              +\frac12\|D_iD_j+D_jD_i\|_{\mathrm{HS}}^2\right)
       \ge0.\notag
\end{align}
For $1<p<2$ and $0<\ell<r$, the power-entropy deficit has expansion
\begin{align}
 \frac\ell r\mathcal H_p(F_\varepsilon)
     -\mathcal H_p(T_\ell F_\varepsilon)
 &=\frac{p(p-1)(2-p)(3-p)}{24}\,
    \ell(r-\ell)\mathcal Q_f(D)\varepsilon^4
       +O(\varepsilon^5).
       \label{eq:matrix-quartic-power-deficit}
\end{align}
For logarithmic entropy the corresponding expansion is
\begin{equation}\label{eq:matrix-quartic-log-deficit}
 \frac\ell r\mathcal H_{\log}(F_\varepsilon)
     -\mathcal H_{\log}(T_\ell F_\varepsilon)
 =\frac{\ell(r-\ell)}{12}
        \mathcal Q_f(D)\varepsilon^4+O(\varepsilon^5).
\end{equation}
Here the mean is $I$; trace-one normalization is not imposed.

The following are equivalent:
\begin{enumerate}
\item $\mathcal Q_f(D)=0$.
\item $D_iD_j=0$ whenever $W_{ij}>0$, including supported self-pairs.
\item For every admissible $\varepsilon$ and every trace function
$\varphi$ for which the quantities are finite,
\[
 \mathcal H_\varphi(T_kF_\varepsilon)
       =\frac{k}{r}\mathcal H_\varphi(F_\varepsilon)
                  \qquad(0\le k\le r).
\]
\item The affine family $F_\varepsilon$ attains equality for one fixed
power $1<p<2$ and one fixed proper retained level for all sufficiently
small $\varepsilon$.
\end{enumerate}
In (3), admissibility means $F_\varepsilon(\alpha)\succeq0$ for every
supported $\alpha$. Thus no higher-order obstruction remains after the
quartic obstruction vanishes.
\end{theorem}

\begin{proof}
Realize an occupation $\alpha$ by an exchangeable ordered list of
particles $(X_1,\ldots,X_r)$, uniformly ordered conditionally on
$\alpha$. The law of $X_1$ is $p$ and that of $(X_1,X_2)$ is $W$.
For any list of previously specified distinct particle positions, let
$\beta$ be their occupation count. Differentiating
\eqref{eq:matrix-direction-translation-identity} by $\partial^\beta$
shows that the conditional mean of the matrix at a new particle
position vanishes:
\begin{equation}\label{eq:singleton-matrix-block-vanishes}
 \mathbb E[D_{X_a}\mid\text{the specified other particle labels}]=0.
\end{equation}
Indeed, its numerator is
$\sum_iD_i\partial^{\beta+e_i}f(\mathbf1)=0$.
The normalizing denominator is positive whenever the conditioning
event has positive probability. This also gives
\[
 T_kH(\beta)=\sum_i\beta_iD_i.
\]

Expand $\operatorname{Tr}(\sum_{a=1}^rD_{X_a})^j$ by partitions of
the $j$ positions according to which particle they use. Every pattern
with a singleton block vanishes by
\eqref{eq:singleton-matrix-block-vanishes}. This argument is valid
without commutativity: after conditioning, the singleton matrix occurs
linearly between fixed matrix factors. Consequently
\begin{align}
 \mathbb E\operatorname{Tr}H^2
    &=r\sum_i p_i\operatorname{Tr}D_i^2,\notag\\
 \mathbb E\operatorname{Tr}H^3
    &=r\sum_i p_i\operatorname{Tr}D_i^3,\notag\\
 \mathbb E\operatorname{Tr}H^4
    &=r\sum_i p_i\operatorname{Tr}D_i^4
                    +r(r-1)\mathcal Q_f(D).
       \label{eq:matrix-translation-moment-formulas}
\end{align}
At order four the only remaining patterns are a single block of size
four and the three pair partitions. Two pair partitions have trace
$\operatorname{Tr}D_i^2D_j^2$, while the alternating partition has
$\operatorname{Tr}D_iD_jD_iD_j$, yielding the stated formula.

Uniform retention keeps the same one- and two-particle laws, replacing
$r$ by $k$ in \eqref{eq:matrix-translation-moment-formulas}. Thus the
second- and third-order deficit moments vanish, while
\[
 \frac\ell r\mathbb E\operatorname{Tr}H^4
     -\mathbb E\operatorname{Tr}(T_\ell H)^4
       =\ell(r-\ell)\mathcal Q_f(D).
\]
Taylor expansion around the scalar matrix $I$ now proves
\eqref{eq:matrix-quartic-power-deficit} and
\eqref{eq:matrix-quartic-log-deficit}: the fourth coefficients of
$(1+x)^p$ and $(1+x)\log(1+x)$ are respectively
$p(p-1)(p-2)(p-3)/24$ and $1/12$.

For Hermitian $A,B$, cyclicity gives
\[
 \|AB\|_{\mathrm{HS}}^2=\operatorname{Tr}A^2B^2,
 \qquad
 \frac12\|AB+BA\|_{\mathrm{HS}}^2
       =\operatorname{Tr}A^2B^2+\operatorname{Tr}ABAB.
\]
This proves the positive sum-of-squares formula and the equivalence
of (1) and (2).

Assume (2). If a supported occupation repeats label $i$, its self-pair
has positive mass, and $D_i^2=0$ forces $D_i=0$. Distinct nonzero
matrices appearing in the same supported occupation have orthogonal
ranges, since $D_iD_j=0$ and the matrices are Hermitian. Therefore
functional calculus splits on those ranges and their orthogonal
complement. Writing
\[
 e_i^\varphi(\varepsilon)
   =\operatorname{Tr}\varphi(I+\varepsilon D_i)
                         -\operatorname{Tr}\varphi(I),
\]
one has the exact identity on every supported occupation
\[
 \operatorname{Tr}\varphi(I+\varepsilon H_\alpha)
                         -\operatorname{Tr}\varphi(I)
       =\sum_i\alpha_i e_i^\varphi(\varepsilon).
\]
The same identity holds for every retained occupation. Since
$\mathbb E\beta_i=(k/r)\mathbb E\alpha_i$ and the retained mean is
again $I$, averaging proves (3). It immediately implies (4), and the
positive fourth coefficient in
\eqref{eq:matrix-quartic-power-deficit} shows that (4) implies (1).
\end{proof}

\begin{remark}[Generality of the moment argument]
No Lorentzian hypothesis is needed for the displayed moment identities
or exact-profile criterion. For Lorentzian laws, these translations
are precisely the Hessian null directions, and the surrounding
contraction theorem holds for every positive semidefinite input.
\end{remark}

\begin{corollary}[Scalar and trace-one qubit forms]
\label{cor:quartic-scalar-qubit}
For a scalar $h\in\mathcal N_f$, put $u(\alpha)=\alpha\cdot h$ and
\[
 Q_f(h)=(h\odot h)^TW(h\odot h).
\]
Then $\mathcal Q_f(h)=3Q_f(h)$, and
\begin{align*}
 \frac\ell r\mathcal H_p(1+\varepsilon u)
       -\mathcal H_p(1+\varepsilon T_\ell u)
 &=\frac{p(p-1)(2-p)(3-p)}8\,
       \ell(r-\ell)Q_f(h)\varepsilon^4+O(\varepsilon^5).
\end{align*}
For the trace-one qubit encoding
\[
 R_\varepsilon(\alpha)=\frac12
     \operatorname{diag}(1+\varepsilon u(\alpha),
                         1-\varepsilon u(\alpha)),
\]
one has
\begin{align*}
 \frac\ell r\mathcal H_{\log}(R_\varepsilon)
       -\mathcal H_{\log}(T_\ell R_\varepsilon)
 &=\frac14\ell(r-\ell)Q_f(h)\varepsilon^4+O(\varepsilon^6),\\
 \frac\ell r\mathcal H_p(R_\varepsilon)
       -\mathcal H_p(T_\ell R_\varepsilon)
 &=2^{1-p}\frac{p(p-1)(2-p)(3-p)}8\,
       \ell(r-\ell)Q_f(h)\varepsilon^4+O(\varepsilon^6).
\end{align*}
The remainders are even because replacing $\varepsilon$ by
$-\varepsilon$ exchanges the two qubit eigenvalues.
\end{corollary}

\begin{proof}
The scalar identity is immediate from the definition of
$\mathcal Q_f$. For the qubit, use $D_i=h_i\operatorname{diag}(1,-1)$,
so $\mathcal Q_f(D)=6Q_f(h)$, and then divide the density by two.
Logarithmic entropy is positively homogeneous, whereas power entropy
is homogeneous of degree $p$.
\end{proof}

\section{Polynomial modes and deletion under label splitting}
\label{sec:deletion-dynamics}

The linear extremal modes belong to a whole family of exact polynomial
eigenfunctions. For independently split labels, this family sits in
an explicit decomposition of every deletion--refill operator.

Let $f(z)=\sum_{|\alpha|=r}\mu(\alpha)z^\alpha$ be a normalized occupation law of degree $r\ge1$. Write $T_\ell F=\E[F\mid\beta]$ for conditional expectation after uniform retention of $\ell$ particles, and $K_\ell=T_\ell^*T_\ell$. The $K_\ell$ are self-adjoint operators on the reference $L^2$ space.

The scalar equality densities already contain exact eigenfunctions of
all orders. The mechanism uses only homogeneity and polynomial
translation invariance; Lorentzianity identifies these translations
from the Hessian and identifies the generating densities as entropy
extremizers.

\begin{theorem}[Polynomial eigenfunctions from translations]
\label{thm:translation-modes}
Let $f(z)=\sum_{|\alpha|=r}\mu(\alpha)z^\alpha$, $r\ge1$, have nonnegative
coefficients and $f(\mathbf1)=1$, and set
$\mathcal N_f=\{h:D_hf=0\}$. For $h\in\mathcal N_f$, define
\[
 \Phi_h(t;\alpha)=\prod_i(1+t h_i)^{\alpha_i}
       =\sum_{j=0}^r t^j P_{j,h}(\alpha),
 \qquad
 \mathcal E_j=\operatorname{span}\{P_{j,h}:h\in\mathcal N_f\}.
\]
For every $0\le\ell\le r$ and $0\le j\le r$,
\[
 K_\ell\big|_{\mathcal E_j}
       =\frac{\binom\ell j}{\binom rj}I,
 \qquad \binom\ell j=0\quad(j>\ell).
\]
The spaces $\mathcal E_j$ are pairwise orthogonal. More precisely, with
$u\odot v=(u_i v_i)_i$,
\[
 \langle P_{j,u},P_{k,v}\rangle_{L^2(\mu)}
   =\mathbf1_{j=k}[z^j]f(\mathbf1+z(u\odot v)).
\]
Consequently the one-particle replacement semigroup acts on
$\mathcal E_j$ by $e^{-jt/r}$. When $f$ is Lorentzian,
$\bigoplus_j\mathcal E_j$ is exactly the linear span of all strictly
positive scalar equality densities.
\end{theorem}

\begin{proof}
Write $a=\mathbf1+t h$. Translation invariance passes to every
partial derivative, so the conditional-mean identity gives
\[
 T_\ell\Phi_h(t;\cdot)(\beta)
 =a^\beta\frac{\partial^\beta f(a)}{\partial^\beta f(\mathbf1)}
 =\Phi_h(t;\beta).
\]
It suffices to prove this for small real $t$, when $a>0$; the
identity then holds polynomially. Each degree-$j$ term in
$\Phi_h(t;\alpha)$ is a product of the labels of $j$ distinct particle
positions. A uniformly retained $\ell$-subset contains those positions
with probability $\binom\ell j/\binom rj$. Thus
\[
 T_\ell^*P_{j,h}(\beta)
   =\frac{\binom\ell j}{\binom rj}P_{j,h}(\alpha),
\]
which proves the eigenvalue formula. The Gram identity follows from
\[
 \mathbb E_\mu[\Phi_u(s;\alpha)\Phi_v(t;\alpha)]
 =f(\mathbf1+su+tv+st(u\odot v))
 =f(\mathbf1+st(u\odot v)).
\]
Compare coefficients of $s^jt^k$. For $\ell=r-1$, the eigenvalue is
$1-j/r$, proving the semigroup assertion. Finally, Theorem~\ref{thm:lorentzian-equality}
identifies the strictly positive extremizers with
$\Phi_h(1;\cdot)$ whenever $\mathbf1+h>0$. Evaluating at $r+1$
sufficiently small distinct values of $t$ and inverting the Vandermonde
matrix recovers each $P_{j,h}$. Conversely each positive extremizer is
the sum of its homogeneous coefficients. In rank one, the same equality
classification follows directly from $\sum_i\mu(e_i)h_i=0$.
\end{proof}

\begin{remark}[Classical specialization]
For $f(z)=(\sum_i p_i z_i)^r$, the translation space is $p^\perp$ and
the Gram kernel is
$(1+st\sum_i p_i u_i v_i)^r$. The spaces above give the classical
multivariate Krawtchouk decomposition of a multinomial law \cite{DiaconisGriffiths}. In general
some $\mathcal E_j$ can vanish, and these spaces need not exhaust the
full eigenspaces of a single $K_\ell$. Their construction is canonical
and simultaneous over all retained levels.
\end{remark}

\subsection{The full decomposition under multinomial splitting}
\label{sec:multinomial-splitting}

The exact spectral decomposition has a probabilistic form that includes
repeated labels. It requires only independent splitting of each coarse
label. The quotient occupation law may be arbitrary.

Let
\[
 g(y)=\sum_{|s|=r}\nu(s)y^s,\qquad g(\one)=1,
\]
be homogeneous with nonnegative coefficients. For each coarse label
$C$, choose a finite set $\mathcal C_C$, with these sets pairwise disjoint,
and a positive probability
law $q_C$ on that set. Define
\[
 f(z)=g\left(\left(\sum_{i\in\mathcal C_C}q_C(i)z_i\right)_C\right).
\]
The law of $f$ first samples the coarse occupation $s$ from $\nu$,
then independently assigns colors with law $q_C$ to the $s_C$
particles of class $C$. Write $K_\ell^f$ for uniform retention of
$\ell$ particles followed by conditional refill.

For a multiindex $a=(a_C)_C$ with $\partial^a g\ne0$, put
\[
 d=|a|,\qquad
 g_a(y)=\frac{\partial^a g(y)}{\partial^a g(\one)},\qquad
 V_C=L^2_0(\mathcal C_C,q_C),\qquad
 V_a=\bigotimes_C\Sym^{a_C}V_C.
\]
The polynomial $g_a$ is a normalized occupation law of degree $r-d$.
Symmetric powers here are the symmetric tensor subspaces with their
product Hilbert inner products.
The index $a_C$ counts centered color slots in class $C$. A conditional
statistic in this sector can survive only if all those slots are retained.

\begin{theorem}[Exact derivative-law decomposition]
\label{thm:multinomial-blocks}
There is an orthogonal decomposition
\[
 L^2(\mu_f)\cong
 \bigoplus_{\substack{a\ge0\\\partial^ag\ne0}}
         L^2(\mu_{g_a})\otimes V_a
\]
under which, simultaneously for every $0\le\ell\le r$,
\[
 K_\ell^f\big|_a\cong
 \begin{cases}
 \displaystyle
 \frac{\binom{\ell}{d}}{\binom rd}
       K_{\ell-d}^{g_a}\otimes I_{V_a},&d\le\ell,\\[6pt]
 0,&d>\ell.
 \end{cases}
\]
For a constant polynomial $g_a$, its rank-zero operator is the
identity on the one-point probability space. Zero-dimensional
$V_a$ may be omitted. The multiplicity of the derivative-law block is
\[
 \dim V_a=\prod_C\dim\Sym^{a_C}V_C.
\]
If every color class has at least two elements, this is
$\prod_C\binom{|\mathcal C_C|+a_C-2}{a_C}$. A singleton class
contributes zero when $a_C>0$, and one when $a_C=0$.
\end{theorem}

\begin{proof}
Fix a coarse occupation $s$. The space of functions of the color
counts in class $C$ identifies with the symmetric functions in
$L^2(q_C)^{\otimes s_C}$. The decomposition
$L^2(q_C)=\R\one\oplus V_C$ gives the orthogonal sum
\[
 \Sym^{s_C}L^2(q_C)
      \cong\bigoplus_{a_C=0}^{s_C}\Sym^{a_C}V_C.
\]
We make the normalization explicit.
For a symmetric tensor $v\in\Sym^{a_C}V_C$, define the symmetric
statistic
\[
 U_{a_C,v}(x_1,\ldots,x_{s_C})
    =\sum_{\substack{J\subseteq[s_C]\\|J|=a_C}}v(x_J).
\]
For two distinct subsets $J,J'$, a coordinate in their symmetric
difference is centered in exactly one factor; integrating it gives
zero. Hence
\[
 \|U_{a_C,v}\|^2_{L^2(q_C^{\otimes s_C})}
       =\binom{s_C}{a_C}\|v\|^2.
\]
Take tensor products across classes, and extend linearly to
$v\in V_a$; write the resulting statistic as $U_{a,v}$.

For a coefficient function $\psi(s-a)$, the squared norm of
$\psi(s-a)U_{a,v}$ is
\[
 \sum_{s\ge a}\nu(s)\prod_C\binom{s_C}{a_C}
          |\psi(s-a)|^2\|v\|^2.
\]
The normalization factor is
\[
 c_a=\sum_s\nu(s)\prod_C\binom{s_C}{a_C}
       =\frac{\partial^a g(\one)}{\prod_C a_C!}.
\]
After division by $c_a$, the coefficient law is precisely $\mu_{g_a}$,
because the coefficient of $y^{s-a}$ in $\partial^ag$ is
$\nu(s)\prod_C(s_C)_{a_C}$, where $(u)_v=u(u-1)\cdots(u-v+1)$.
Thus
\[
 J_a(\psi\otimes v)
       =c_a^{-1/2}\psi(s-a)U_{a,v}
\]
is isometric. The fixed-$s$ tensor decompositions show that the
images of all $J_a$ are orthogonal and exhaustive.

Now condition on a retained detailed color occupation, with coarse
counts $b_C$ and total size $\ell$. If $b_C<a_C$ for some $C$, the
conditional statistic is zero. Otherwise $b-a$ is a retained occupation
for the derivative law. Every term of $U_{a,v}$ that
uses an unretained centered color has conditional mean zero.
The surviving terms are exactly $U_{a,v}$ evaluated on the
retained colors. The conditional distribution of the coefficient
occupation $s-a$ is the posterior for $g_a$ after retention of
$b-a$. Indeed, the two posterior weights agree up to a
factor independent of $s$, by
\[
 (s)_a\binom{s-a}{b-a}
      =(b)_a\binom{s}{b},
\]
with products taken coordinatewise.

The corresponding normalization at level $\ell$ is
\[
 c_{a,\ell}
   =\E\prod_C\binom{b_C}{a_C}
   =\frac{\binom{\ell}{d}}{\binom rd}\,c_a.
\]
This follows by counting a fixed set of $d$ marked particle
positions: uniform retention keeps it with probability
$\binom{\ell}{d}/\binom rd$. Therefore, under the isometries at
the input and output levels, conditional expectation on the
$a$-sector is
\[
 \sqrt{\frac{\binom{\ell}{d}}{\binom rd}}\,
          T_{\ell-d}^{g_a}\otimes I.
\]
Taking its adjoint followed by itself gives the claimed $K_\ell$
block. If $\ell<d$, no such marked set survives and the block is
zero.
\end{proof}

\begin{corollary}[Isospectrality under color redistribution]
\label{cor:general-splitting-isospectral}
Fix $g$ and the cardinality of each color class. Varying the positive
laws $q_C$ gives one unitary equivalence intertwining all the
operators $K_\ell^f$ simultaneously.
\end{corollary}
\begin{proof}
The derivative laws $g_a$ are unchanged. Choose an isometry between
the old and new centered spaces $V_C$, and take its symmetric
and tensor powers in the decomposition.
\end{proof}

For one coarse class, $g(y)=y^r$, every derivative law is a point
mass. The decomposition becomes the usual multivariate Krawtchouk
decomposition of a multinomial law \cite{DiaconisGriffiths}. For a multiaffine
$g$ with matroid support, each nonzero multiindex $a$ is the indicator of an independent
set in the support matroid, and differentiation becomes contraction.
This gives the next section.

\subsection{Contraction minors of a parallel extension}
Let $M$ be a loopless rank-$r$ matroid with positive fields
$\lambda_i$, and let $S$ be its simplification. Its elements are
the parallel classes $C$. Put
\[
 w_C=\sum_{i\in C}\lambda_i,\qquad q_C(i)=\lambda_i/w_C,
 \qquad V_C=L^2_0(C,q_C).
\]
The simplified basis $I$ has law proportional to
$\prod_{C\in I}w_C$; given $I$, its representatives are independent
with laws $q_C$. For $A\in\mathcal I(S)$, let $\pi^A$ be the
basis law of $S/A$ with fields $w_C$, equivalently the simplified
law conditioned to contain $A$.

\begin{corollary}[All retained-level walks through contraction minors]
\label{thm:parallel-blocks}
One orthogonal decomposition, independent of $\ell$, identifies
\[
 L^2(\mu_M)\cong\bigoplus_{A\in\mathcal I(S)}
       L^2(\pi^A)\otimes\bigotimes_{C\in A}V_C.
\]
On the summand indexed by $A$, with $j=|A|\le\ell$,
\[
 K_\ell^M\cong
       \frac{\binom\ell j}{\binom rj}
             K_{\ell-j}^{S/A}\otimes I.
\]
The summand is killed when $j>\ell$. A rank-zero contracted walk
is the identity on a one-point space. Redistributing positive
fields within parallel classes while preserving every $w_C$ gives
one unitary equivalence intertwining all retained-level walks.
In particular, every strictly positive scalar entropy equality tilt
preserves their entire spectra simultaneously.
\end{corollary}
\begin{proof}
In Theorem~\ref{thm:multinomial-blocks}, multiaffinity makes $a$
the indicator of an independent set $A$, and the derivative law
is $\pi^A$. This gives the decomposition. The simplification and all its
conditional laws depend only on $w_C$; the representative spaces
retain dimension $|C|-1$ under redistribution. Apply
Corollary~\ref{cor:general-splitting-isospectral}. Theorem~\ref{thm:matroid-equality}
identifies positive scalar equality tilts with exactly these
redistributions.
\end{proof}

The constants on the contracted spaces give the polynomial modes
of degree $j$, with dimension
\[
 m_j=\sum_{\substack{A\in\mathcal I(S)\\|A|=j}}
                         \prod_{C\in A}(|C|-1).
\]
Deleting one representative from each parallel class gives a matroid
on $N'=\sum_C(|C|-1)$ elements with independent-set numbers $m_j$.
Thus $(m_j)_{j=0}^{N'}$ is ultra log-concave by Mason's theorem
\cite{ALOVIII,BrandenHuh}: $(m_j/\binom{N'}j)_{j=0}^{N'}$ is log-concave. These are the dimensions of the canonical
mode spaces, not necessarily the full multiplicities of a given
eigenvalue: other contraction blocks can contribute that eigenvalue.

The operator decomposition also applies to every multiaffine Lorentzian
coefficient law. Lemma~\ref{lem:lorentzian-missing-pair-geometry}
gives independent representative choices after merging parallel
classes; the contraction laws are then normalized derivatives of
the quotient polynomial, without a product-field assumption.

For such coefficient laws, quotient translations can give additional
polynomial modes beyond the within-class spaces counted by $m_j$.

\section{Exchange edges, isoperimetry, and rigidity}
\label{sec:exchange-edges}

Conditioning on a retained face leaves a clique of possible completions.
Factor entropy inside each such clique, then count how often each pair
of configurations occurs. This converts deletion contraction into a
complete inequality on individual exchange edges. Equality remembers
the clique sizes and identifies exactly the binary factors of the
basis graph.

Let $\mu$ be positive on a family $\mathcal B$ of $r$-element sets.
For an unordered pair $e=\{B,B'\}$, let $K_e$ replace the two matrix
values by their $\mu$-weighted mean and leave every other value fixed.
Define its entropy loss by
\[
 J_e(F)=\mathcal H_\mu(F)-\mathcal H_\mu(K_eF).
\]

\begin{lemma}[Complete-graph factorization]
\label{lem:matrix-clique-factorization}
For every probability law $\nu$ on a finite set and every
positive semidefinite matrix function $F$,
\[
 \mathcal H_\nu(F)\le
 \sum_{x<y}H(\nu(x),\nu(y);F(x),F(y)),
\]
where
\[
 H(a,b;A,B)=\Tr\left[aA\log A+bB\log B
 -(aA+bB)\log\frac{aA+bB}{a+b}\right].
\]
\end{lemma}
\begin{proof}
The two-point loss is jointly convex in $(A,B)$, homogeneous in
$(a,b)$, and nondecreasing in each mass. To see the last assertion,
write $M=(aA+bB)/(a+b)$ and
$D_g(A\Vert M)=\Tr[A(\log A-\log M)-A+M]$. Then
$H=aD_g(A\Vert M)+bD_g(B\Vert M)$ and
$\partial_aH=D_g(A\Vert M)\ge0$. Joint convexity follows from
relative entropy; regularization covers singular matrices.

Induct on the support size. Split off $x$ of mass $a$, and put
$b=1-a$. The entropy chain rule separates the remaining conditional
entropy from the loss between $F(x)$ and the remaining mean. Convexity
bounds this latter loss by
\[
 \sum_{y\ne x}\frac{\nu(y)}b H(a,b;F(x),F(y))
 =\sum_{y\ne x}H(a\nu(y)/b,\nu(y);F(x),F(y)).
\]
Since $\nu(y)\le b$, mass monotonicity bounds these terms by the
required cross pairs. Homogeneity and induction treat the remaining
conditional entropy.
\end{proof}

\begin{theorem}[Intersection-weighted edge factorization]
\label{thm:deletion-edge-factorization}
Suppose uniform deletion to $\ell$ elements, $0\le\ell<r$, satisfies
$\mathcal H_{\mu_\ell}(F_\ell)\le\eta\mathcal H_\mu(F)$ completely,
with $\eta<1$. Then
\begin{equation}\label{eq:intersection-edge-factorization}
 \mathcal H_\mu(F)\le
 \frac1{\binom r\ell(1-\eta)}
 \sum_{B<B'}\binom{|B\cap B'|}\ell J_{\{B,B'\}}(F).
\end{equation}
For every homogeneous multiaffine Lorentzian law, in particular,
\begin{equation}\label{eq:lorentzian-edge-factorization}
 \mathcal H_\mu(F)\le\sum_{e\in E}J_e(F),
 \qquad E=\{\{B,B'\}:|B\cap B'|=r-1\}.
\end{equation}
Under the gap hypothesis \eqref{eq:hodge-gap}, the right side can be
multiplied by
\begin{equation}\label{eq:strict-edge-factor}
 C_\kappa=\frac{r+(r-1)\kappa}{r(1+\kappa)}.
\end{equation}
\end{theorem}
\begin{proof}
The exact chain rule gives
\[
 \mathcal H_\mu(F)-\mathcal H_{\mu_\ell}(F_\ell)
 =\sum_I\mu_\ell(I)\mathcal H_{\mu(\cdot\mid I)}(F).
\]
Apply Lemma~\ref{lem:matrix-clique-factorization} to each completion
law. For $I\subseteq B$, Bayes' formula gives
$\mu_\ell(I)\mu(B\mid I)=\mu(B)/\binom r\ell$.
Homogeneity therefore makes every shared $I$ contribute
$J_{\{B,B'\}}(F)/\binom r\ell$. There are
$\binom{|B\cap B'|}\ell$ such faces, proving
\eqref{eq:intersection-edge-factorization}.
For Lorentzian laws take $\ell=r-1$ and $\eta=(r-1)/r$;
for $r=1$ the output is the mean and $\eta=0$.
The intersection weight is exactly one on exchange edges.
Using $\eta=a_{r-1}/a_r$ from Theorem~\ref{thm:all-level-gap} instead
produces $C_\kappa$.
\end{proof}

\begin{corollary}[Edge heat flow and matrix concentration]
\label{cor:edge-heat-flow}
Let $L_E=\sum_{e\in E}(I-K_e)$, so each edge updates at rate one.
For every Lorentzian basis law,
\[
 \mathcal H_\mu(e^{-tL_E}F)\le e^{-t}\mathcal H_\mu(F).
\]
The complete modified logarithmic Sobolev constant is at least one,
and at least $C_\kappa^{-1}$ under \eqref{eq:hodge-gap}.
Theorem~\ref{thm:nonlinear-laplace} applies with this constant and
the actual $\Gamma_{L_E}$. The exponent one is optimal uniformly
over positive product fields.
\end{corollary}
\begin{proof}
Convexity gives
$J_e(F)\le\E_\mu\Tr[((I-K_e)F)\log F]$.
Sum the edge factorization and integrate entropy production along the
semigroup. For sharpness use the rank-one two-base law
$\mu=(\varepsilon,1-\varepsilon)$. Starting from $(a,1-a)$,
$0<a<1$, the first coordinate at time $t$ is
$\varepsilon+(a-\varepsilon)e^{-t}$. As $\varepsilon\downarrow0$,
the output-to-input relative-entropy ratio tends to $e^{-t}$.
\end{proof}

For the uniform law, the unit-edge graph Laplacian is $L_G=2L_E$.
When $E\ne\varnothing$, the discrete kernel
$|E|^{-1}\sum_eK_e$ has complete contraction at most $1-1/|E|$:
convexity bounds its output entropy by the average output entropy
of the edge updates. These clocks differ from one
deletion--replacement update; we keep their normalizations explicit.

\subsection{Complete entropy decay for the unit-edge clock}
\label{subsec:discovery-complete-four}

For the uniform law, the sharp matrix midpoint estimate below
converts edge factorization into complete entropy production, in the
modified logarithmic Sobolev framework of \cite{GJL2020}.
No normalization of $\E\Tr F$ is required.

\begin{lemma}[Complete midpoint entropy bound]
\label{lem:discovery-midpoint-production}
For positive definite matrices $A,B$, put $M=(A+B)/2$. Then
\begin{equation}
 J(A,B):=
 \frac12\Tr(A\log A+B\log B)-\Tr(M\log M)
 \le\frac18\Tr[(A-B)(\log A-\log B)].
 \label{eq:discovery-midpoint-one-eighth}
\end{equation}
The constant $1/8$ is sharp. Semidefinite inputs follow by regularization,
with the usual extended-value interpretation of entropy production.
\end{lemma}

\begin{proof}
Set $D=(A-B)/2$ and
\[
 \varphi(t)=\frac12\Tr\bigl[(M+tD)\log(M+tD)
                           +(M-tD)\log(M-tD)\bigr].
\]
For $s\ge0$ let
$C_s=(M+sI)^{-1/2}D(M+sI)^{-1/2}$. Since $M\pm D>0$,
all eigenvalues of $C_s$ lie in $(-1,1)$. The resolvent formula for the
Fr\'echet derivative of the logarithm and cyclicity of trace give
\[
 \varphi''(t)=\int_0^\infty\frac12\Tr
 \left[C_s^2\bigl((I+tC_s)^{-2}+(I-tC_s)^{-2}\bigr)\right]ds.
\]
Diagonalize $C_s$ inside this integrand. Its even power-series coefficients
are nonnegative. Since $\varphi'(0)=0$, integrating twice and using
monotone convergence yields the exact expansion
\begin{equation}
 \varphi(t)-\varphi(0)
 =\sum_{m\ge1}\frac{t^{2m}}{2m}
       \int_0^\infty\Tr(C_s^{2m})\,ds,
 \qquad 0\le t\le1.
 \label{eq:discovery-positive-midpoint-expansion}
\end{equation}
The coefficients are finite: $C_s=O(s^{-1})$ at infinity, while positive
definiteness controls $s$ near zero. Nonnegativity of the coefficients gives
$\varphi(1)-\varphi(0)\le\varphi'(1)/2$. Direct differentiation gives
\[
 \varphi'(1)=\frac14\Tr[(A-B)(\log A-\log B)],
\]
which proves \eqref{eq:discovery-midpoint-one-eighth}. Taking scalar
$A=1+\varepsilon$, $B=1-\varepsilon$ and sending $\varepsilon$ to zero
shows sharpness. Adding $\varepsilon I$ and passing to the limit handles
the boundary.
\end{proof}

\begin{theorem}[Complete entropy decay and the sharp universal rate]
\label{thm:discovery-complete-rate-four}
Let $G$ be a finite graph with uniform law $\mu$. For each edge $e$ let
$K_e$ average its two matrix values and leave all other values fixed.
Assume, for some $C>0$, the complete edge entropy factorization
\begin{equation}
 \mathcal H_\mu(F)
 \le C\sum_{e\in E(G)}
       \bigl[\mathcal H_\mu(F)-\mathcal H_\mu(K_eF)\bigr].
 \label{eq:discovery-edge-factorization-hypothesis}
\end{equation}
Then the unit-edge Laplacian satisfies
\begin{align}
 \mathbb E_\mu\Tr[(L_GF)\log F]&\ge\frac4C\mathcal H_\mu(F),
 \label{eq:discovery-complete-production-four}\\
 \mathcal H_\mu(e^{-tL_G}F)&\le e^{-4t/C}\mathcal H_\mu(F)
 \qquad(t\ge0).
 \label{eq:discovery-complete-decay-four}
\end{align}
Both statements hold in every finite matrix dimension. In particular,
Theorem~\ref{thm:deletion-edge-factorization} gives $C=1$ for every
uniform matroid basis graph. The resulting rate four is optimal
uniformly over this class, already on a binary edge.
For rank $r\ge2$ under \eqref{eq:hodge-gap}, one may take
$C=C_\kappa$, giving the stronger rate
\begin{equation}\label{eq:strict-unit-edge-rate}
 \frac4{C_\kappa}
   =\frac{4r(1+\kappa)}{r+(r-1)\kappa}.
\end{equation}
Every simple matroid of rank at least two admits $\kappa>0$ by
\cite{MNY}, so this guaranteed rate is strictly greater than four.
\end{theorem}

\begin{proof}
If $N$ is the number of vertices, the loss on $e=\{x,y\}$ is
$J_e(F)=(2/N)J(F(x),F(y))$. Lemma~\ref{lem:discovery-midpoint-production}
therefore gives
\begin{align*}
 \sum_eJ_e(F)
 &\le\frac1{4N}\sum_{\{x,y\}\in E(G)}
   \Tr[(F(x)-F(y))(\log F(x)-\log F(y))]\\
 &=\frac14\mathbb E_\mu\Tr[(L_GF)\log F].
\end{align*}
Combine this with \eqref{eq:discovery-edge-factorization-hypothesis}.
The derivative of entropy along $e^{-tL_G}$ is the negative of this
production, so integration proves
\eqref{eq:discovery-complete-decay-four}. Regularization gives the
semidefinite case.

On a binary edge take $F=1+\varepsilon h$, with
$\mathbb E h=0$ and $L_Gh=2h$. Entropy is quadratic to leading order,
so production divided by entropy tends to four; for each fixed $t>0$,
the output-to-input entropy ratio tends to $e^{-4t}$.
\end{proof}

\subsection{Small sets and scalar hypercontractivity}

The scalar consequences below follow from scalar deletion
\cite{R4CGM} and the clique lemma. Montenegro--Son established the
small-set profile for balanced matroids \cite{MontenegroSon}, and
Anari--Liu--Oveis Gharan--Vinzant proved expansion at least one for
all matroid basis graphs \cite{ALOV}. The complete factorization
retains matrix-valued inputs throughout the same argument.

\begin{corollary}[Small sets and scalar hypercontractivity]
\label{cor:matroid-edge-isoperimetry}
Let $G$ be the basis graph of a matroid with $N$ bases and uniform
law $\mu$. Every vertex set $U$, with
$a=\min(|U|,N-|U|)>0$, satisfies
\begin{equation}\label{eq:matroid-small-set}
 |\partial U|\ge a\log_2(N/a).
\end{equation}
For $(L_Gg)(B)=\sum_{B'\sim B}(g(B)-g(B'))$,
\begin{equation}\label{eq:matroid-gross}
 \operatorname{Ent}_\mu(g^2)
 \le\langle g,L_Gg\rangle_\mu
 =\frac1N\sum_{\{B,B'\}\in E}(g(B)-g(B'))^2.
\end{equation}
Thus $e^{-tL_G}$ contracts $L^p(\mu)$ to $L^q(\mu)$ whenever
$1<p\le q$ and $q-1\le(p-1)e^{4t}$.
Both constants are sharp uniformly over matroids.
Under \eqref{eq:hodge-gap}, the right sides of
\eqref{eq:matroid-small-set} and \eqref{eq:matroid-gross} improve by
the factors $C_\kappa^{-1}$ and $C_\kappa$, respectively, and
hypercontractivity holds when $q-1\le(p-1)e^{4t/C_\kappa}$.
\end{corollary}
\begin{proof}
Apply \eqref{eq:lorentzian-edge-factorization} to the indicator of the
smaller side. Its entropy is $(a/N)\log(N/a)$; each crossing edge
contributes $(\log2)/N$. For $g^2$, use the two-point inequality
\[
 u^2\log u^2+v^2\log v^2
 -(u^2+v^2)\log\frac{u^2+v^2}{2}\le(u-v)^2.
\]
This inequality is strict unless $u=v$. For a proof, absolute
values reduce to $u,v\ge0$, and homogeneity reduces a nonzero pair
to $u^2+v^2=2$. Write $u^2=1+s$, $v^2=1-s$, with $0\le s\le1$.
The difference between the right and left sides is
\[
 R(s)=2-2\sqrt{1-s^2}
       -(1+s)\log(1+s)-(1-s)\log(1-s).
\]
It satisfies $R(0)=R'(0)=0$ and
$R''(s)=2(1-s^2)^{-3/2}-2(1-s^2)^{-1}>0$ for $0<s<1$.
Continuity treats $s=1$, proving the inequality and its strictness.

For hypercontractivity, differentiate
$\log\|e^{-tL_G}g\|_{q(t)}$ for positive $g$. The Gross inequality
and the edge bound
\[
 (u^{q-1}-v^{q-1})(u-v)
 \ge\frac{4(q-1)}{q^2}(u^{q/2}-v^{q/2})^2
\]
make this derivative nonpositive when $q'\le4(q-1)$.
Approximation and positivity of the semigroup treat arbitrary inputs.
Retaining $C_\kappa$ in the same edge factorization gives all three
improvements, with $q'\le4(q-1)/C_\kappa$.
Direct sums of $U_{1,2}$ have cube basis graphs: subcubes attain
\eqref{eq:matroid-small-set}, and perturbations of a constant along
one binary factor approach equality in \eqref{eq:matroid-gross}.
\end{proof}

\subsection{A forbidden spectral interval for matroid basis graphs}
\label{subsec:discovery-spectral-interval}

We continue with the uniform basis law and unit-edge Laplacian $L_G$.
Write $\mathcal E_G(g,g)=\langle g,L_Gg\rangle_\mu$; by
\eqref{eq:matroid-gross},
\begin{equation}\label{eq:discovery-unit-edge-energy}
 \mathcal E_G(g,g)=\frac1N\sum_{\{B,B'\}\in E(G)}(g(B)-g(B'))^2.
\end{equation}
Loops do not affect the basis graph, and coloops contribute no varying
coordinate.

\begin{theorem}[Universal spectral separation]
\label{thm:discovery-spectral-separation}\label{thm:basis-graph-spectral-rigidity}
Let $M$ be a finite matroid with exactly $k$ direct summands isomorphic to
$U_{1,2}$. Then
\begin{equation}
 \Spec(L_G)\subseteq\{0,2\}\cup[3,\infty).
 \label{eq:discovery-spectral-exclusion}
\end{equation}
The eigenvalue-two space has dimension $k$ and is
\begin{equation}\label{eq:basis-graph-eigenvalue-two}
 \ker(L_G-2I)=\operatorname{span}\left\{
 \mathbf1_{\{a\in B\}}-\mathbf1_{\{b\in B\}}:
 \{a,b\}\text{ is a }U_{1,2}\text{ direct summand}\right\}.
\end{equation} In particular, if
$M$ has no $U_{1,2}$ direct summand and has at least two bases, then
$\lambda_{\mathrm{gap}}(L_G)\ge3$.
The threshold three is sharp: the basis graphs of $U_{1,3}$ and $U_{2,3}$
are triangles.
\end{theorem}

The proof uses the following weighted statement. Its rates must be kept
explicit: they are the rates produced by collapsing parallel classes in
a uniform basis graph.

\begin{lemma}[Weighted upper-star factorization]
\label{lem:discovery-weighted-simple}
Let $M$ be a loopless matroid with at least two bases, and let $m_i>0$.
On its bases put
\[
 \pi(B)=\frac1Z\prod_{i\in B}m_i,
 \qquad
 (L_mg)(B)=
 \sum_{\substack{i\in B,\ j\notin B\cr B-i+j\in\mathcal B(M)}}
 m_j\bigl(g(B)-g(B-i+j)\bigr).
\]
If $g(M)$ denotes the girth of $M$, then
\begin{equation}
 \mathcal E_m(g,g):=\langle g,L_mg\rangle_\pi
 \ge g(M)\min_i m_i\,\Var_\pi(g).
 \label{eq:discovery-weighted-girth}
\end{equation}
In particular, a simple matroid with all $m_i\ge1$ has gap at least three
for this weighted generator.
\end{lemma}

\begin{proof}
Write $r=\rank(M)$ and $s=|E(M)|-r>0$. Under complementation, the law
$\pi$ is the positive product-field basis law of the dual matroid,
with fields $1/m_i$. Delete one particle from a dual basis, and let
$I^*$ be the retained dual face of size $s-1$. Its complement $U$ has
size $r+1$. The possible primal bases inside $U$ are
$B_i=U\setminus\{i\}$, where $i$ ranges over the unique circuit
$C\subseteq U$. Define
\[
 W_U=\sum_{B\subset U}\pi(B),\qquad
 S_C=\sum_{i\in C}\frac1{m_i},\qquad
 \pi_U(i)=\frac{1/m_i}{S_C}.
\]
The retained-face probability is $W_U/s$. The scalar Lorentzian deletion
inequality, or its variance linearization, gives
\[
 \Var\bigl(\mathbb E[g\mid I^*]\bigr)
 \le\frac{s-1}{s}\Var_\pi(g).
\]
The variance chain rule therefore yields
\begin{equation}
 \Var_\pi(g)
 \le\sum_U W_U\Var_{\pi_U}\bigl(g(B_i)\bigr).
 \label{eq:discovery-upper-star-variance}
\end{equation}
Thus the dual rank cancels: the weights in
\eqref{eq:discovery-upper-star-variance} are $W_U$, rather than the
retained-face probabilities $W_U/s$.

On the clique indexed by $C$, take transition rates $q(i,j)=m_i$ for
$i\ne j$. They are reversible for $\pi_U$, with conductance $1/S_C$
on every unordered pair. With $c=|C|$ and $D_m=\operatorname{diag}(m_i)_{i\in C}$,
the positive generator is
\[
 L_C=D_m(cI-\mathbf1\mathbf1^{\mathsf T}).
\]
Its nonzero eigenvalues agree with those of the symmetric matrix
\[
 (cI-\mathbf1\mathbf1^{\mathsf T})^{1/2}D_m
             (cI-\mathbf1\mathbf1^{\mathsf T})^{1/2}.
\]
Since $D_m\succeq(\min_{i\in C}m_i)I$, each such eigenvalue is at least
$c\min_{i\in C}m_i\ge g(M)\min_i m_i$. Consequently
\[
 \Var_{\pi_U}(g)
 \le\frac{\mathcal E_C(g,g)}{g(M)\min_i m_i}.
\]
The local energies sum exactly to the global energy. Indeed,
\begin{equation}
 \frac{W_U}{S_C}
 =\frac{\prod_{i\in U}m_i}{Z}
 =\pi(U\setminus\{i\})m_i.
 \label{eq:discovery-upper-star-conductance}
\end{equation}
This is the conductance of the global edge
$U\setminus\{i\}\leftrightarrow U\setminus\{j\}$, and every exchange
edge has a unique union $U$. Hence
$\sum_UW_U\mathcal E_C(g,g)=\mathcal E_m(g,g)$.
Substitution in \eqref{eq:discovery-upper-star-variance} proves the claim.
\end{proof}

\begin{corollary}[Weighted girth and cogirth bounds]
\label{cor:discovery-girth-cogirth}
Under the hypotheses of Lemma~\ref{lem:discovery-weighted-simple},
\begin{equation}
 \lambda_{\mathrm{gap}}(L_m)
 \ge\max\left\{
       g(M)\min_i m_i,
       \min_{C^*\text{ a cocircuit}}\sum_{i\in C^*}m_i
                    \right\}.
 \label{eq:discovery-weighted-cogirth}
\end{equation}
In particular, for the unit-edge basis graph,
$\lambda_{\mathrm{gap}}(L_G)\ge\max\{g(M),g(M^*)\}$.
\end{corollary}

\begin{proof}
The first bound is Lemma~\ref{lem:discovery-weighted-simple}.
For the second, use primal deletion to a face $I$ of size $r-1$.
Its completion set is the cocircuit $C^*=E(M)\setminus\operatorname{cl}(I)$.
Writing $T=\sum_{i\in C^*}m_i$, the conditional law on completions is
$m_i/T$. The local generator with rates $q(i,j)=m_j$ is
$T(I-\mathbb E_{m/T})$, so its gap is exactly $T$ when the clique has
at least two vertices. A singleton completion has zero conditional variance.
The same variance factorization and the identity
\[
 W_I\frac{m_i}{T}m_j
 =\frac{(\prod_{a\in I}m_a)m_im_j}{Z}
 =\pi(I\cup\{i\})m_j
\]
show that the clique energies again sum exactly to $\mathcal E_m$.
This proves the second bound. Set every $m_i=1$ for the last assertion.
\end{proof}

\begin{lemma}[Centered sectors on parallel classes]
\label{lem:discovery-parallel-sectors}
If $M$ has no $U_{1,2}$ direct summand, its unit-edge basis graph has
gap at least three whenever it has more than one vertex.
\end{lemma}

\begin{proof}
Apply Corollary~\ref{thm:parallel-blocks} at uniform original fields.
The simplification $S$ has fields $m_C=|C|$, and
$V_C=L^2_0(C,\mathrm{unif})$. Write $I(B)$ for the simplified basis
and $x_C(B)$ for its representative in $C$. With
$c_A=\Pr_\mu(A\subseteq I(B))^{-1/2}$, an elementary tensor is realized by
\[
 u(B)=
 \begin{cases}
 c_Ah(I(B))\displaystyle\prod_{C\in A}v_C(x_C(B)),&A\subseteq I(B),\\
 0,&A\not\subseteq I(B),
 \end{cases}
 \qquad \sum_{x\in C}v_C(x)=0.
\]
The coefficient law on $I\supseteq A$ is proportional to
$\prod_{C\in I}m_C$. A nonzero contrast requires $m_C\ge2$.

We compute the generator on this function. For an internal swap
in a required class $C\in A$,
\[
 \sum_{\substack{y\in C\\y\ne x_C}}
       \bigl(v_C(x_C)-v_C(y)\bigr)=m_Cv_C(x_C).
\]
Thus internal swaps contribute $\sum_{C\in A}m_C$; swaps in other
classes leave $u$ unchanged. An exchange $I-C+D$ retaining $A$
preserves every contrast factor. Summing over the $m_D$ possible
inserted representatives gives
\[
 (L^{(A)}h)(I)=
 \sum_{\substack{C\in I\setminus A,\ D\notin I\\
                  I-C+D\in\mathcal B(S)}}
 m_D\bigl(h(I)-h(I-C+D)\bigr).
\]
If $C\in A$ is removed, the new value is zero. These exchanges
therefore give the multiplication operator
\[
 V_A(I)=\sum_{C\in A}
       \sum_{\substack{D\notin I\\ I-C+D\in\mathcal B(S)}}m_D.
\]
Finally take $A\not\subseteq I(B)$. A neighbor can contain $A$
only by inserting the unique missing class $C$. All its representatives
occur with unit rate, so their contribution sums to zero by
$\sum_{x\in C}v_C(x)=0$. This proves invariance, including at the
boundary, and gives
\begin{equation}\label{eq:discovery-sector-generator}
 L_G\big|_{\mathcal H_A}\cong
 \left[\left(\sum_{C\in A}m_C\right)I+L^{(A)}+V_A\right]\otimes I.
\end{equation}

Both $L^{(A)}$ and $V_A$ are nonnegative. For $A=\varnothing$,
Lemma~\ref{lem:discovery-weighted-simple} gives gap at least three
on the simplification; if it has one basis, this sector consists of constants.
Every other sector has lower bound at least three unless
$A=\{C\}$ and $m_C=2$. This is the only possible sector below three. Here $C$ is not a coloop of $S$, since
$M$ has no $U_{1,2}$ direct summand. Exchange with a basis omitting
$C$ gives an exit from $C$ at every simplified basis containing it.
Thus $V_A\ge1$, and the lower bound is again $2+1=3$.
\end{proof}

\begin{proof}[Proof of Theorem~\ref{thm:discovery-spectral-separation}]
Split off the binary direct factors:
\[
 M=U_{1,2}^{\oplus k}\oplus M_0.
\]
The basis graph is the Cartesian product of a $k$-cube and the basis
graph of $M_0$, and the unit-edge Laplacians add. The cube has eigenvalues
$2j$, $0\le j\le k$, with the usual degree-$j$ binary characters.
By Lemma~\ref{lem:discovery-parallel-sectors}, the remaining factor is
either a single vertex or has all its nonconstant eigenvalues at least
three. The product spectrum consists of their sums. Its only positive
eigenvalue below three is therefore two, attained precisely by one
binary contrast times constants on the other factors. Triangles show
the sharpness of three.
\end{proof}

\begin{corollary}[Feder--Mihail equality; Alaniz Pintos's Problem~8.4]
\label{cor:pintos-feder-mihail}
Let $M$ have at least two bases and let $G$ be its basis graph.
For $0<p\le1/\Delta(G)$, put $P_{\mathrm{FM}}=I-pL_G$.
Then
\[
 \operatorname{gap}(P_{\mathrm{FM}})=2p
 \quad\Longleftrightarrow\quad
 M\text{ has a }U_{1,2}\text{ direct summand}.
\]
The entire eigenspace at $1-2p$ is spanned by the contrasts on
these binary summands. If there is no such summand, then
$\operatorname{gap}(P_{\mathrm{FM}})\ge3p$.
Thus the extremal matroids asked for in
\cite[Problem 8.4]{Pintos2026} are exactly those with a binary
direct component.
\end{corollary}
\begin{proof}
Since $I-P_{\mathrm{FM}}=pL_G$, the positive gap is
$p\lambda_{\mathrm{gap}}(L_G)$ and the eigenspaces are unchanged.
Apply Theorem~\ref{thm:discovery-spectral-separation}. This includes
rank one: its basis graph is $K_N$, with nonzero Laplacian eigenvalue
$N$, so equality requires $N=2$.
\end{proof}

\begin{corollary}[Quantitative rigidity, heat flow, and expansion]
\label{cor:discovery-spectral-consequences}
Let $\Pi_{\mathrm{bin}}$ be the orthogonal projection onto the
binary-coordinate contrasts. For every real function $g$,
\begin{equation}
 \mathcal E_G(g,g)-2\Var_\mu(g)
 \ge\|g-\mathbb E_\mu g-\Pi_{\mathrm{bin}}g\|_2^2.
 \label{eq:discovery-poincare-stability}
\end{equation}
The coefficient one is sharp. For mean-zero $g$ and $t\ge0$,
\begin{equation}
 \|e^{-tL_G}g-e^{-2t}\Pi_{\mathrm{bin}}g\|_2
 \le e^{-3t}\|g-\Pi_{\mathrm{bin}}g\|_2.
 \label{eq:discovery-graph-remainder}
\end{equation}
If there is no binary direct factor, then every vertex set $U$ with
$0<a=|U|\le N/2$ satisfies
\begin{equation}
 |\partial U|\ge3a(1-a/N)\ge\tfrac32 a.
 \label{eq:discovery-three-halves-expansion}
\end{equation}
Together with \eqref{eq:matroid-small-set}, this gives
$|\partial U|\ge a\max\{\log_2(N/a),3(1-a/N)\}$.
\end{corollary}

\begin{proof}
Apply the spectral theorem to
\eqref{eq:discovery-spectral-exclusion}. On the orthogonal complement
of constants and binary contrasts, $L_G-2I\succeq I$, which proves
\eqref{eq:discovery-poincare-stability}; an eigenfunction on a triangle
shows sharpness. The same decomposition proves the heat-flow estimate.
For the expansion bound apply the gap-three inequality to $\mathbf1_U$:
$\Var_\mu(\mathbf1_U)=(a/N)(1-a/N)$ and
$\mathcal E_G(\mathbf1_U,\mathbf1_U)=|\partial U|/N$.
\end{proof}

\begin{remark}
The spectral threshold three is optimal. The expansion conclusion applies,
in particular, to every connected matroid with more than two elements.
\end{remark}

\begin{corollary}[Spanning-tree exchange gaps]
\label{cor:spanning-tree-exchange-gaps}
Let $H$ be a finite connected loopless graph with more than one spanning tree.
For $m_e>0$, give its trees law proportional to $\prod_{e\in T}m_e$
and use the generator
\[
 (L_{\mathrm{ex},m}u)(T)=
 \sum_{\substack{e\in T,\ f\notin T\\T-e+f\ \mathrm{a\ tree}}}
 m_f\bigl(u(T)-u(T-e+f)\bigr).
\]
Then
\[
 \lambda_{\mathrm{gap}}(L_{\mathrm{ex},m})
 \ge\max\left\{\operatorname{girth}(H)\min_e m_e,\,
                  \min_{C\ \mathrm{a\ bond}}\sum_{e\in C}m_e\right\}.
\]
For unit rates this is at least the maximum of the girth and edge
connectivity of $H$.
\end{corollary}
\begin{proof}
Apply Corollary~\ref{cor:discovery-girth-cogirth} to the graphic matroid:
its circuits are cycles and its cocircuits are bonds.
\end{proof}

\begin{theorem}[All edge-isoperimetric equality sets]
\label{thm:basis-graph-isoperimetric-rigidity}
Separate all binary components as
$M=U_{1,2}^{\oplus b}\oplus M_0$.
For $0<|U|\le N/2$, equality in \eqref{eq:matroid-small-set} holds
if and only if
\[
 U=C\times\mathcal B(M_0),
\]
where $C$ is a subcube of $\{0,1\}^b$. Equivalently, $U$ fixes one
choice in some binary components and leaves every other component
unrestricted. In particular $|U|/N=2^{-k}$ for an integer $k\ge1$.
Every nontrivial cut in a connected matroid with more than two active
elements satisfies strict inequality.
\end{theorem}
\begin{proof}
For rank one, $G=K_N$ and the claim follows directly from
$|\partial U|=|U|(N-|U|)$; equality requires $N=2$, $|U|=1$.
Assume $r\ge2$. Equality in the edge bound requires equality in
deletion and in each conditional clique bound, since every discarded
term is nonnegative. Apply Theorem~\ref{thm:matroid-equality} to
$(N/|U|)\mathbf1_U$. Its product-tilt form shows that $U$ consists
of the bases avoiding selected forbidden labels; every parallel class
contains an allowed label.

If a completion clique has $c$ vertices, with $0<u<c$ inside $U$,
equality in its indicator factorization requires
\[
 \log_2(c/u)=c-u.
\]
Writing $k=c-u\ge1$ gives
$\log_2(1+k/u)\le\log_2(1+k)\le k$, with equality only when
$u=k=1$. Thus every mixed completion clique has exactly two vertices.

A parallel class containing both allowed and forbidden labels must
therefore have size two: extend that class to a simplified basis and
choose allowed labels in its other classes; the resulting facet
has a mixed completion clique containing the whole class.
If its simplified element is not a coloop,
choose an exchange removing that class from a simplified basis, and
choose allowed representatives in all the other classes. The resulting
facet has a mixed completion clique containing the two labels of the
cut class and an allowed label from another class. This contradicts
its required size. Every cut class is consequently a binary direct
summand, giving the asserted cylinder. Conversely, fixing $k$ binary
choices gives $N2^{-k}$ bases, each with $k$ outgoing edges, and attains
the bound.
\end{proof}

This recovers the classical subcube equality cases on the discrete cube
\cite{EllisKellerLifshitz} and identifies their exact occurrence inside
every matroid basis graph.

\begin{corollary}[Sharpness of the Gross constant]
\label{cor:basis-graph-gross-rigidity}
For a matroid with at least two bases, the optimal constant $C_G$ in
$\operatorname{Ent}_\mu(g^2)\le C_G\langle g,L_Gg\rangle_\mu$
equals one if and only if $M$ has a $U_{1,2}$ direct summand.
\end{corollary}
\begin{proof}
A perturbation of a constant by a binary-factor contrast gives ratios
tending to one. Conversely, the two-point inequality in the proof of
\eqref{eq:matroid-gross} is strict on every edge with different endpoint
values. The basis graph is connected, so the ratio is strictly below
one for every nonconstant $g$. Normalize $\|g\|_2=1$.
Near either constant function write $g=a\one+h$, where
$\E_\mu h=0$. In this finite-dimensional space, Taylor expansion
gives, uniformly in the direction of $h$,
\[
 \operatorname{Ent}_\mu(g^2)=2\|h\|_2^2+o(\|h\|_2^2),\qquad
 \langle g,L_Gg\rangle_\mu=\langle h,L_Gh\rangle_\mu.
\]
Thus the limiting upper bound for the ratio is
$2/\lambda_{\mathrm{gap}}(L_G)$. If there is no binary summand,
Theorem~\ref{thm:basis-graph-spectral-rigidity} makes this strictly
less than one. The rest of the normalized sphere is compact and
has a continuous ratio strictly below one, completing the proof.
\end{proof}

\section{Matrix concentration for arbitrary nonlinear observables}
\label{sec:nonlinear-concentration}

We now turn entropy decay into concentration for an arbitrary
Hermitian observable. The proof has two stages. First, follow a
matrix density along the reversible heat flow and bound the change
in its expectation of the observable. Then optimize over densities
in the Gibbs variational formula. The matrix mean survives both
stages, keeping the variance matrix inside the final exponential
and allowing an arbitrary Hermitian offset.

Let $\mu$ be a positive probability measure on a finite set $\Omega$,
and let
\[
 (LA)(x)=\sum_y q(x,y)(A(x)-A(y)),\qquad
 \mu(x)q(x,y)=\mu(y)q(y,x),\quad q(x,y)\geq0\quad(x\ne y),
\]
be an irreducible reversible Markov generator, with $P_t=e^{-tL}$.
For a Hermitian matrix function $A$, its matrix carr\'e du champ is
\[
 \Gamma_L(A)(x)=\frac12\sum_y q(x,y)(A(x)-A(y))^2.
\]
We use the entropy-production convention
\begin{equation}\label{eq:nonlinear-cmlsi}
 \alpha\mathcal H_\mu(G)
 \leq \mathcal I_L(G):=\E_\mu\Tr[(LG)\log G],\qquad G>0.
\end{equation}
When \eqref{eq:nonlinear-cmlsi} holds in every matrix dimension with
the same $\alpha>0$, it is a complete modified logarithmic Sobolev
inequality.

\subsection{Transport and the matrix Laplace comparison}

For a density with matrix mean $M$, variance is measured by
$\Tr(MV)$. This state-dependent quantity is the link between
transport and the matrix Laplace transform.

\begin{theorem}[Matrix transport and Laplace comparison]
\label{thm:nonlinear-laplace}
Suppose \eqref{eq:nonlinear-cmlsi} holds in matrix dimension $d$.
Let $A$ take values in the Hermitian $d\times d$ matrices, put
$\overline A=\E_\mu A$, and
suppose a fixed positive semidefinite matrix $V$ satisfies
\begin{equation}\label{eq:nonlinear-variance-bound}
                    \Gamma_L(A)(x)\preceq V\qquad(x\in\Omega).
\end{equation}
For every $G\succeq0$ with $\Tr\E_\mu G=1$, writing $M=\E_\mu G$,
\begin{equation}\label{eq:nonlinear-transport}
 \left|\E_\mu\Tr[(A-\overline A)G]\right|
 \leq 2\sqrt{\frac{\Tr(MV)}{\alpha}\mathcal H_\mu(G)}.
\end{equation}
For every Hermitian $K$ and every real $\theta$,
\begin{equation}\label{eq:nonlinear-laplace}
 \E_\mu\Tr e^{K+\theta A}
 \leq\Tr\exp\left(K+\theta\overline A+
                                  \frac{\theta^2}{\alpha}V\right).
\end{equation}
Both conclusions hold uniformly in $d$ under the complete inequality.
\end{theorem}

\begin{proof}
Subtract $\overline A$ first.

\emph{The edge estimate.}
Entropy production pairs a difference of densities with a difference
of logarithms. Their logarithmic mean converts this pairing into
a Hilbert-space Cauchy--Schwarz inequality.
For positive definite matrices $B,C$,
define the positive self-adjoint operator on the Hilbert--Schmidt space
\[
 \mathcal J_{B,C}(X)=\int_0^1 B^sXC^{1-s}\,ds.
\]
The operators of left multiplication by $B$ and right multiplication
by $C$ commute and are positive. Differentiating $B^sC^{1-s}$ gives
\begin{equation}\label{eq:nonlinear-logmean}
 B-C=\mathcal J_{B,C}(\log B-\log C).
\end{equation}
For Hermitian $X$, the weighted arithmetic--geometric mean inequality,
applied in eigenbases of $B$ and $C$, gives
\begin{equation}\label{eq:nonlinear-arithmetic}
 \Tr[X\mathcal J_{B,C}(X)]
 \leq\frac12\Tr[X^2(B+C)].
\end{equation}
Indeed, before integration the coefficient of the squared absolute
value of the $(i,j)$ matrix entry of $X$ is $b_i^sc_j^{1-s}$,
which is at most $sb_i+(1-s)c_j$.

Assume first that every $G(x)$ is positive definite. Reversibility,
\eqref{eq:nonlinear-logmean}, and Cauchy--Schwarz in the direct sum
of the Hilbert--Schmidt spaces with weights
$\mu(x)q(x,y)/2$ yield
\begin{align}
 \left|\E_\mu\Tr[A\,LG]\right|^2
 &\leq\mathcal I_L(G)\,
 \frac12\sum_{x,y}\mu(x)q(x,y)
 \Tr\bigl[(A(x)-A(y))
    \mathcal J_{G(x),G(y)}(A(x)-A(y))\bigr]
 \nonumber\\
 &\leq\mathcal I_L(G)\,
                  \E_\mu\Tr[G\Gamma_L(A)]
 \leq\mathcal I_L(G)\Tr(MV).
 \label{eq:nonlinear-fisher-cs}
\end{align}
The second inequality is \eqref{eq:nonlinear-arithmetic}; the two
endpoint terms agree after exchanging $x$ and $y$.

\emph{Integration along the heat flow.}
Set $G_t=P_tG$, $h(t)=\mathcal H_\mu(G_t)$, and
$i(t)=\mathcal I_L(G_t)$. Its matrix mean remains $M$,
$h'(t)=-i(t)$, and $i(t)\geq\alpha h(t)$. Moreover
$G_t\to M$ and $h(t)\to0$. Integrating
\eqref{eq:nonlinear-fisher-cs} along the semigroup gives
\[
 \left|\E_\mu\Tr[AG]\right|
 \leq\sqrt{\Tr(MV)}\int_0^\infty\sqrt{i(t)}\,dt
 \leq\sqrt{\frac{\Tr(MV)}{\alpha}}
       \int_0^\infty\frac{i(t)}{\sqrt{h(t)}}\,dt
 =2\sqrt{\frac{\Tr(MV)}{\alpha}h(0)}.
\]
The integrals stop if $h$ reaches zero. If $h(0)=0$, the density
is constant and the assertion is immediate. Replacing $G$ by
$(G+\varepsilon I)/(1+\varepsilon d)$ and passing to the limit
proves \eqref{eq:nonlinear-transport} for all $G\succeq0$.

\emph{Optimization over the matrix mean.}
Young's inequality applied to \eqref{eq:nonlinear-transport} gives
\[
 \theta\E_\mu\Tr[AG]-\mathcal H_\mu(G)
 \leq\frac{\theta^2}{\alpha}\Tr(MV).
\]
The Gibbs variational formula on the finite direct sum of matrix
algebras, with the weights $\mu(x)$, reads
\begin{align*}
 \log\E_\mu\Tr e^{K+\theta A}
 &=\sup_{\substack{G\succeq0\\\Tr\E_\mu G=1}}
   \{\Tr(MK)+\theta\E_\mu\Tr[AG]
                       -\mathcal H_\mu(G)-\Tr(M\log M)\}\\
 &\leq\sup_{\substack{M\succeq0\\\Tr M=1}}
   \{\Tr[M(K+\theta^2V/\alpha)]-\Tr(M\log M)\}\\
 &=\log\Tr e^{K+\theta^2V/\alpha}.
\end{align*}
This proves the centered version of \eqref{eq:nonlinear-laplace}.
Replacing $K$ by $K+\theta\overline A$ proves the stated formula.
\end{proof}

Huang--Tropp~\cite{HuangTropp} obtain nonlinear matrix subgaussian
concentration from the Bakry--\'Emery curvature criterion.
Gao--Junge--LaRacuente prove noncommutative transport from modified
logarithmic Sobolev inequalities and its equivalence to conditional
$L_p$ concentration~\cite[Theorems 6.7 and 6.10]{GJL2020}.
See also Gao--Junge--LaRacuente--Li~\cite[Section 5.4]{GJLL}
and De Palma--Rouz\'e~\cite{DePalmaRouze}.
The proof above follows the same entropy-decay argument on a finite
reversible chain. It keeps $\Tr(MV)$ in the transport estimate,
where the Lipschitz seminorm in the cited results uses $\|V\|$.
Optimizing the matrix mean $M$ then retains both $V$ and the
arbitrary offset $K$ in the trace comparison.

\begin{corollary}[Subgaussian matrix concentration]
\label{cor:nonlinear-gaussian}
Under Theorem~\ref{thm:nonlinear-laplace}, put $v=\|V\|$.
If $v>0$, then for every $t>0$,
\begin{align}
 \Pr\{\lambda_{\max}(A-\overline A)\geq t\}
 &\leq\inf_{\theta>0}e^{-\theta t}
                         \Tr e^{\theta^2V/\alpha}
 \leq d\exp\left(-\frac{\alpha t^2}{4v}\right),
 \label{eq:nonlinear-tail}\\
 \Pr\{\|A-\overline A\|\geq t\}
 &\leq2d\exp\left(-\frac{\alpha t^2}{4v}\right).
 \label{eq:nonlinear-norm-tail}
\end{align}
The one-sided bound has the intrinsic-dimension refinement
\begin{equation}\label{eq:nonlinear-intrinsic}
 \Pr\{\lambda_{\max}(A-\overline A)\geq t\}
 \leq e\,\frac{\Tr V}{\|V\|}
                      \exp\left(-\frac{\alpha t^2}{4v}\right).
\end{equation}
If $V=0$, the observable is constant on the reference support.
\end{corollary}
\begin{proof}
Take $K=0$ in \eqref{eq:nonlinear-laplace} for $A-\overline A$,
apply the trace Laplace bound~\cite{R5Tropp}, and optimize at
$\theta=\alpha t/(2v)$. Apply the same argument to $-A$ for the
opposite tail. If $V=0$, every edge difference vanishes.

For \eqref{eq:nonlinear-intrinsic}, let
$\psi(u)=e^u-u-1$ and $X=A-\overline A$. Since $\E X=0$,
\[
 \E\Tr\psi(\theta X)
 \leq\Tr(e^{\theta^2V/\alpha}-I)
 \leq\frac{\Tr V}{v}(e^{\theta^2v/\alpha}-1).
\]
On the upper-tail event the left-hand random variable is at least
$\psi(\theta t)$. Put $E=\alpha t^2/(4v)$ and use the same
$\theta$, so that $\theta t=2E$ and $\theta^2v/\alpha=E$.
For $E\geq1$,
$(e^E-1)/(e^{2E}-2E-1)\leq e^{1-E}$.
For $E<1$ the claimed bound exceeds one because $\Tr V/v\geq1$.
\end{proof}

\subsection{Affine variance and Bernstein bounds}

An affine variance bound allows the local fluctuations to grow
with the observable itself. Along the heat flow, its square root
has a bounded rate of change. Integrating that bound and completing
one square gives the Bernstein denominator, while retaining the
full negative-parameter range.

\begin{theorem}[An affine carr\'e-du-champ bound]
\label{thm:nonlinear-affine}
Suppose \eqref{eq:nonlinear-cmlsi} holds in dimension $d$.
Let $A$ be a centered Hermitian observable and suppose
\begin{equation}\label{eq:nonlinear-affine-assumption}
 \Gamma_L(A)(x)\preceq aA(x)+B\qquad(x\in\Omega),
 \qquad a\geq0,\quad B\succeq0,
\end{equation}
where $a$ and $B$ are fixed. For every Hermitian $K$,
\begin{equation}\label{eq:nonlinear-affine-laplace}
 \E_\mu\Tr e^{K+\theta A}
 \leq\Tr\exp\left(K+\frac{\theta^2}{\alpha-a\theta}B\right)
 \qquad(\theta<\alpha/a).
\end{equation}
When $a=0$, the parameter range is all of $\mathbb R$.
When $a>0$, every normalized matrix density $G$ satisfies
\begin{equation}\label{eq:nonlinear-affine-transport}
 \left|\sqrt{a\E_\mu\Tr(AG)+\Tr(MB)}-\sqrt{\Tr(MB)}\right|
 \leq a\sqrt{\mathcal H_\mu(G)/\alpha},\qquad M=\E_\mu G.
\end{equation}
Thus the comparison includes all negative $\theta$ as well.
\end{theorem}
\begin{proof}
The case $a=0$ is Theorem~\ref{thm:nonlinear-laplace}.
Suppose $a>0$ and first take $G>0$. Set $G_t=P_tG$,
$f(t)=\E_\mu\Tr(AG_t)$, $h(t)=\mathcal H_\mu(G_t)$,
$i(t)=\mathcal I_L(G_t)$, and $b=\Tr(MB)$.
The Cauchy--Schwarz step in \eqref{eq:nonlinear-fisher-cs}
and \eqref{eq:nonlinear-affine-assumption} give
\[
 |f'(t)|^2\leq i(t)\E_\mu\Tr[G_t\Gamma_L(A)]
                     \leq i(t)[af(t)+b].
\]
The quantity $af(t)+b$ is nonnegative, since
$aA+B\succeq\Gamma_L(A)\succeq0$.
For $\varepsilon>0$, let $g_\varepsilon(t)=
\sqrt{af(t)+b+\varepsilon}$. Then
$|g_\varepsilon'(t)|\leq(a/2)\sqrt{i(t)}$. Since $f(t)\to0$,
integration as in Theorem~\ref{thm:nonlinear-laplace} yields
\[
 |\sqrt{af(0)+b+\varepsilon}-\sqrt{b+\varepsilon}|
 \leq\frac a2\int_0^\infty\sqrt{i(t)}\,dt
 \leq a\sqrt{h(0)/\alpha}.
\]
Letting $\varepsilon\downarrow0$ and then regularizing $G$
proves \eqref{eq:nonlinear-affine-transport}, including $b=0$.

For brevity put $f=\E_\mu\Tr(AG)$ and $h=\mathcal H_\mu(G)$.
Writing $u=\sqrt{af+b}$, the transport bound gives
$h\geq(\alpha/a^2)(u-\sqrt b)^2$. For every $\theta<\alpha/a$,
completing the square gives
\begin{align*}
 \theta f-h
 &\leq\frac\theta a(u^2-b)-\frac\alpha{a^2}(u-\sqrt b)^2\\
 &=\frac{\theta^2b}{\alpha-a\theta}
   -\frac{\alpha-a\theta}{a^2}
        \left(u-\frac{\alpha\sqrt b}{\alpha-a\theta}\right)^2
 \leq\frac{\theta^2}{\alpha-a\theta}\Tr(MB).
\end{align*}
The Gibbs variational formula, followed by optimization over $M$,
now proves \eqref{eq:nonlinear-affine-laplace}.
\end{proof}

Affine conditional variance bounds also appear in the matrix
Stein-pair inequalities of
Mackey--Jordan--Chen--Farrell--Tropp~\cite[Theorem 4.1]{MackeyEtAl}.
Theorem~\ref{thm:nonlinear-affine} uses complete entropy production
for an arbitrary nonlinear observable and retains a non-scalar
matrix $B$ and every Hermitian offset $K$.

\begin{corollary}[Nonlinear matrix Bernstein bounds]
\label{cor:nonlinear-affine-tails}
Under Theorem~\ref{thm:nonlinear-affine}, put $v=\|B\|$.
If $v>0$, then for $t>0$,
\begin{align}
 \Pr\{\lambda_{\max}(A)\geq t\}
 &\leq\inf_{0<\theta<\alpha/a}e^{-\theta t}
               \Tr e^{\theta^2B/(\alpha-a\theta)}
 \leq d\exp\left(-\frac{\alpha t^2}{4v+2at}\right),
 \label{eq:nonlinear-bernstein}\\
 \Pr\{\lambda_{\min}(A)\leq-t\}
 &\leq d\exp\left(-\frac{\alpha t^2}{4v}\right).
 \label{eq:nonlinear-affine-lower}
\end{align}
For $a>0$, the optimized scalar exponent in the upper tail is
\begin{equation}\label{eq:nonlinear-affine-exact-tail}
 \Pr\{\lambda_{\max}(A)\geq t\}
 \leq d\exp\left[-\frac\alpha{a^2}
                         (\sqrt{v+at}-\sqrt v)^2\right].
\end{equation}
For $a>0$ and $0<t<v/a$, the lower-tail exponent improves
to $-\alpha(\sqrt v-\sqrt{v-at})^2/a^2$; for $t\ge v/a$ the lower
event is impossible. The elementary upper bound
\eqref{eq:nonlinear-bernstein} also holds with prefactor
$e\Tr(B)/\|B\|$. If $B=0$, then $A=0$ on the reference support.
\end{corollary}
\begin{proof}
Apply the trace Laplace bound to
\eqref{eq:nonlinear-affine-laplace}. The positive choice
$\theta=\alpha t/(2v+at)$ gives
\eqref{eq:nonlinear-bernstein}; direct optimization for $a>0$
gives $\theta=(\alpha/a)(1-\sqrt{v/(v+at)})$ and
\eqref{eq:nonlinear-affine-exact-tail}. For negative $\theta$,
$\theta^2/(\alpha-a\theta)\leq\theta^2/\alpha$, which gives
\eqref{eq:nonlinear-affine-lower}. Direct optimization gives
the stronger stated lower exponent. If $a>0$, then
$aA+B\succeq0$ gives $A\succeq-B/a\succeq-(v/a)I$.
Equality on a unit vector would force every adjacent difference of $A$
to annihilate that vector, since $\Gamma_L(A)\preceq aA+B$.
The same value would propagate through the irreducible support,
contradicting $\E A=0$. Hence $\lambda_{\min}(A)>-v/a$.

For the intrinsic-dimension bound, apply
$\psi(u)=e^u-u-1$ and $\E A=0$ to obtain
\[
 \E\Tr\psi(\theta A)
 \leq\frac{\Tr B}{v}
                  (e^{\theta^2v/(\alpha-a\theta)}-1).
\]
At the positive parameter above, put $E=\alpha t^2/(4v+2at)$.
Then $\theta t=2E$ and
$\theta^2v/(\alpha-a\theta)=E$. The same elementary estimate
as in Corollary~\ref{cor:nonlinear-gaussian} gives the claimed
prefactor. Finally, if $B=0$, taking expectations of
\eqref{eq:nonlinear-affine-assumption} gives
$\E\Gamma_L(A)\preceq a\E A=0$, so every edge difference
vanishes and the centered observable is zero.
\end{proof}

\subsection{Natural matroid basis laws and Lorentzian occupations}

\begin{corollary}[All nonlinear Lorentzian observables]
\label{cor:nonlinear-lorentzian}
Let $\mu$ be a degree-$r$ Lorentzian occupation law, and let $K_0$
be deletion of one particle followed by conditional replacement.
Set $L=I-K_0$. For every Hermitian observable $A$ satisfying
$\Gamma_L(A)\preceq V$, every Hermitian offset $K$, and every
real $\theta$,
\begin{equation}\label{eq:nonlinear-lorentzian-laplace}
 \E_\mu\Tr e^{K+\theta A}
 \leq\Tr\exp(K+\theta\E_\mu A+r\theta^2V).
\end{equation}
If every allowed one-particle replacement changes $A$ in operator
norm by at most $a>0$, then
\begin{equation}\label{eq:nonlinear-exchange-lipschitz}
 \Pr\{\lambda_{\max}(A-\E A)\geq t\}
 \leq d\exp\left(-\frac{t^2}{2ra^2}\right).
\end{equation}
Both statements hold, in particular, for arbitrary Hermitian functions
on the bases of every matroid at every positive product field.
Under the affine assumptions of Theorem~\ref{thm:nonlinear-affine},
a centered observable with $\Gamma_L(A)\preceq aA+B$ satisfies
\begin{equation}\label{eq:nonlinear-lorentzian-affine}
 \E_\mu\Tr e^{K+\theta A}
 \leq\Tr\exp\left(K+\frac{r\theta^2}{1-ra\theta}B\right)
 \qquad(\theta<1/(ra)).
\end{equation}
For $a=0$, every real $\theta$ is allowed.
Under \eqref{eq:hodge-gap}, replace $r$ in these conclusions by
$[r+(r-1)\kappa]/(1+\kappa)$.
\end{corollary}
\begin{proof}
Corollary~\ref{r4l:replacement} gives
\eqref{eq:nonlinear-cmlsi} with $\alpha=1/r$. Theorem~\ref{thm:all-level-gap}
gives the improved constant. The replacement Lipschitz assumption
implies $\Gamma_L(A)\preceq a^2I/2$, since $K_0$ is stochastic.
Apply Theorem~\ref{thm:nonlinear-laplace} and
Corollary~\ref{cor:nonlinear-gaussian}.
The affine comparison follows from
Theorem~\ref{thm:nonlinear-affine}.
\end{proof}

This includes the nonlinear observable class discussed in
Kathuria~\cite[Sections 4--5]{Kathuria}: every Hermitian function
with bounded exchange differences is allowed. The bound is subgaussian
throughout its range and applies to arbitrary matroids.
Exponential matrix concentration for Lorentzian measures was already
obtained from the scalar-to-matrix Poincar\'e principle
of Garg--Kathuria--Srivastava~\cite{GargKathuriaSrivastava}.

\begin{corollary}[Rank and corank in each component]
\label{cor:component-gaussian}
Let $\mu$ be a positive product-field basis law of a matroid $M$.
For each connected component $C$, put
$s_C=\min(\operatorname{rk}_M(C),|C|-\operatorname{rk}_M(C))$.
Suppose that every basis exchange inside $C$ changes a Hermitian
observable $A$ in operator norm by at most $a_C$. Then
\begin{equation}\label{eq:component-gaussian}
 \Pr\{\lambda_{\max}(A-\E A)\geq t\}
 \leq d\exp\left(-\frac{t^2}{2\sum_Cs_Ca_C^2}\right)
\end{equation}
when the denominator is positive. If it is zero, $A$ is constant
on the reference support.
\end{corollary}
\begin{proof}
In each component of rank larger than its corank, replace the
basis by its complement and the field by its reciprocal. This
pushes $\mu$ to a product-field basis law of a direct sum of
matroids of total rank $s=\sum_Cs_C$. If $s=0$, the support is
a singleton. Otherwise components of rank zero may be removed. For its deletion--replacement generator, a step
chooses component $C$ with probability $s_C/s$ and exchanges
only within that component. Thus
$\Gamma_L(A)\preceq(2s)^{-1}\sum_Cs_Ca_C^2 I$.
The complete modified logarithmic Sobolev constant is at least
$1/s$. Apply \eqref{eq:nonlinear-tail} and pull the observable
back through the complementation bijection.
\end{proof}

\section{Matrix concentration for arbitrary-order subconfiguration averages}
\label{sec:ustatistics}

A subconfiguration average inherits the kernel's matrix Laplace bound,
with reciprocal deletion coefficient playing the role of sample size.

Let $\alpha$ have a degree-$r$ Lorentzian occupation law, and fix
$1\leq\ell\leq r$. For a Hermitian matrix kernel $B_\beta$
on the retained $\ell$-occupations, define
\begin{equation}\label{eq:ustatistic}
 U(\alpha)=\sum_\beta D_{r\to\ell}(\beta\mid\alpha)B_\beta,
 \qquad \overline B=\E_{\mu_\ell}B_\beta=\E_\mu U.
\end{equation}
For a matroid basis $B$ this is
\[
 U(B)=\binom r\ell^{-1}\sum_{\substack{I\subseteq B\\|I|=\ell}}B_I.
\]
Thus the kernel may depend on an entire $\ell$-subset. Throughout
this section one may take $c=r/\ell$, or the larger value
$c=a_r/a_\ell$ under the hypothesis of
Theorem~\ref{thm:all-level-gap}.

\begin{theorem}[Matrix Laplace comparison]\label{thm:u-laplace}
For every Hermitian $K$ and every real $\theta$,
\begin{equation}\label{eq:u-laplace}
 \E\Tr e^{K+\theta U}
 \leq\Tr\exp\!\left(K+c\log\E_{\mu_\ell}
                                        e^{\theta B_\beta/c}\right).
\end{equation}
In particular, with $M_\theta=\E_{\mu_\ell}e^{\theta B_\beta/c}$,
\[
 \E\Tr e^{\theta U}\leq\Tr M_\theta^c\leq d\|M_\theta\|^c.
\]
\end{theorem}
\begin{proof}
In Corollary~\ref{cor:all-order-trace}, put
$A_\beta=e^{\theta B_\beta/c}$. Then
$cD^*\log A=\theta U$. The last comparisons follow from spectral
calculus; the real exponent $c$ need not be an integer.
\end{proof}

\begin{theorem}[Chernoff and binomial bounds]\label{thm:u-chernoff}
Suppose $0\preceq B_\beta\preceq RI_d$, where $R>0$, and put
$m_+=\lambda_{\max}(\overline B)$ and
$m_-=\lambda_{\min}(\overline B)$. For $\delta\geq0$ and $m_+>0$,
\begin{equation}\label{eq:u-chernoff-upper}
 \Pr\{\lambda_{\max}(U)\geq(1+\delta)m_+\}
 \leq d\left(\frac{e^\delta}{(1+\delta)^{1+\delta}}\right)^{cm_+/R}.
\end{equation}
For $0<\delta<1$ and $m_->0$,
\begin{equation}\label{eq:u-chernoff-lower}
 \Pr\{\lambda_{\min}(U)\leq(1-\delta)m_-\}
 \leq d\left(\frac{e^{-\delta}}{(1-\delta)^{1-\delta}}\right)^{cm_-/R}.
\end{equation}
For $s_+=m_+/R<u\leq1$ and $0\leq u<s_-=m_-/R$,
the respective refinements are
\[
 \Pr\{\lambda_{\max}(U)\geq Ru\}\leq d e^{-c\operatorname{kl}(u\Vert s_+)},
 \qquad
 \Pr\{\lambda_{\min}(U)\leq Ru\}\leq d e^{-c\operatorname{kl}(u\Vert s_-)},
\]
where $\operatorname{kl}(u\Vert s)=u\log(u/s)+(1-u)\log((1-u)/(1-s))$.
Endpoint cases follow by limits, and all bounds may be capped at one.
\end{theorem}
\begin{proof}
For every real $t$, spectral calculus gives
\[
 e^{tB_\beta/R}\preceq I+(e^t-1)B_\beta/R.
\]
If $t>0$, the largest eigenvalue of the expectation on the right
is at most $1+s_+(e^t-1)$. Theorem~\ref{thm:u-laplace} and
the trace exponential Markov bound, with $\theta=ct/R$, give
\[
 \Pr\{\lambda_{\max}(U)\geq Ru\}
 \leq d\inf_{t>0}
       \exp\{c[-tu+\log(1-s_++s_+e^t)]\}.
\]
Optimizing gives the binary relative-entropy bound. For $t<0$,
the largest eigenvalue is controlled by $s_-$, and the same argument
with $\lambda_{\min}(U)$ gives the lower bound. Finally use
$\log(1+s(e^t-1))\leq s(e^t-1)$ and take
$t=\log(1+\delta)$ or $t=\log(1-\delta)$ to obtain the displayed
Chernoff estimates.
\end{proof}

\begin{theorem}[Centered Bennett and Bernstein bounds]
\label{thm:u-bernstein}
Put $C_\beta=B_\beta-\overline B$ and suppose
$\|C_\beta\|\leq R$ on the support, where $R>0$. Define
\[
 V=\E_{\mu_\ell}C_\beta^2,\qquad \sigma^2=\|V\|,
 \qquad X=U-\overline B.
\]
If $\sigma^2>0$, then for every $t>0$,
\begin{align}
 \Pr\{\lambda_{\max}(X)\geq t\}
 &\leq d\exp\!\left[-\frac{c\sigma^2}{R^2}
                 h\!\left(\frac{Rt}{\sigma^2}\right)\right]
                                                        \label{eq:u-bennett}\\
 &\leq d\exp\!\left[-\frac{ct^2}{2\sigma^2+2Rt/3}\right],
                                                        \label{eq:u-bernstein}
\end{align}
where $h(v)=(1+v)\log(1+v)-v$. The same bounds hold for
$\lambda_{\max}(-X)$, and a union bound gives a two-sided norm
bound with prefactor $2d$. If $V=0$, then $X=0$ almost surely.
\end{theorem}
\begin{proof}
For $s\geq0$ and $-R\leq x\leq R$,
\[
 e^{sx}\leq1+sx+\frac{e^{sR}-1-sR}{R^2}x^2.
\]
Apply this by spectral calculus and use $\E C_\beta=0$. Since
$\log$ is operator monotone and $\log(I+A)\preceq A$ for
$A\succeq0$,
\[
 c\log\E e^{\theta C_\beta/c}
 \preceq\frac{c}{R^2}(e^{\theta R/c}-1-\theta R/c)V.
\]
The Laplace comparison applies to the centered kernel itself, so
\[
 \E\Tr e^{\theta X}
 \leq\Tr\exp\!\left[
 \frac{c}{R^2}(e^{\theta R/c}-1-\theta R/c)V\right].
\]
The exponential trace is monotone on Hermitian matrices.
Bounding the last expression by its largest eigenvalue and optimizing
at $\theta=(c/R)\log(1+Rt/\sigma^2)$ proves
\eqref{eq:u-bennett}. The scalar inequality
$h(v)\geq v^2/(2+2v/3)$ gives \eqref{eq:u-bernstein}.
If $V=0$, positivity gives $C_\beta=0$ on the reference support.
\end{proof}

\begin{corollary}[Intrinsic dimension]\label{cor:u-intrinsic}
Under Theorem~\ref{thm:u-bernstein}, suppose $\sigma^2>0$ and put
$d_{\rm eff}=\Tr V/\|V\|$. For $t>0$,
\begin{align*}
 \Pr\{\lambda_{\max}(X)\geq t\}
 &\leq \min\!\left\{1,\inf_{\theta>0}
 \frac{\Tr(e^{g_c(\theta)V}-I)}
      {e^{\theta t}-\theta t-1}\right\}\\
 &\leq e\,d_{\rm eff}
       \exp\!\left[-\frac{ct^2}{2\sigma^2+2Rt/3}\right],
 \qquad
 g_c(\theta)=\frac{c}{R^2}
        (e^{\theta R/c}-1-\theta R/c).
\end{align*}
The same upper bound holds for the opposite tail, and the norm
bound has prefactor $2e\,d_{\rm eff}$.
\end{corollary}
\begin{proof}
Use the intrinsic-dimension argument of
Corollary~\ref{cor:nonlinear-gaussian}, replacing
$\theta^2V/\alpha$ by $g_c(\theta)V$ in its spectral test.
The Laplace bound in Theorem~\ref{thm:u-bernstein} gives the first
inequality and bounds its numerator by
$d_{\rm eff}(e^{g_c(\theta)\sigma^2}-1)$.
Set
\[
 E=\frac{ct^2}{2\sigma^2+2Rt/3},\qquad
 \theta=\frac{ct}{\sigma^2+Rt/3}.
\]
Then $\theta t=2E$, $\theta R/c<3$, and
$e^u-1-u\le u^2/[2(1-u/3)]$ gives
$g_c(\theta)\sigma^2\le E$.
The same scalar ratio bound and small-exponent case from that
corollary give $e\,d_{\rm eff}e^{-E}$. Apply this to $-X$ for the
opposite tail.
\end{proof}

\begin{corollary}[Matrix sums on every matroid]
\label{cor:matroid-matrix-sums}
Let $B$ be a product-field random basis of any rank-$r$ matroid, $r\ge1$,
and put $Y=\sum_{i\in B}A_i$, $Q=\E Y$, for $d\times d$ matrices.
If $0\preceq A_i\preceq RI$ with $R>0$, then, for $\delta\ge0$ and
$q_+=\lambda_{\max}(Q)>0$,
\[
 \Pr\{\lambda_{\max}(Y)\ge(1+\delta)q_+\}
 \le d\left(\frac{e^\delta}{(1+\delta)^{1+\delta}}\right)^{q_+/R}.
\]
For $0<\delta<1$ and $q_-=\lambda_{\min}(Q)>0$,
\[
 \Pr\{\lambda_{\min}(Y)\le(1-\delta)q_-\}
 \le d\left(\frac{e^{-\delta}}{(1-\delta)^{1-\delta}}\right)^{q_-/R}.
\]
For arbitrary Hermitian $A_i$, write $p_i=\Pr(i\in B)$ and
\[
 V_0=\sum_i p_i A_i^2-\frac{Q^2}{r},\qquad
 R_0=\max_{i:p_i>0}\|A_i-Q/r\|.
\]
If $\|V_0\|>0$, then for every $t>0$,
\[
 \Pr\{\lambda_{\max}(Y-Q)\ge t\}
 \le d\exp\!\left(-\frac{t^2}{2\|V_0\|+2R_0t/3}\right).
\]
The same bound holds for the opposite tail.
\end{corollary}
\begin{proof}
Use the one-particle kernel $B_i=rA_i$ in the preceding theorems.
Then $U=Y$, $c=r$, its bound is $rR$, and its centered variance
and radius are $rV_0$ and $rR_0$. The rank cancels in each exponent.
Positivity of $V_0$ follows from operator Jensen for the
probabilities $p_i/r$.
\end{proof}

\begin{corollary}[Spectral approximation from arbitrary matroid bases]
\label{cor:matroid-spectral}
Let $B$ be a product-field random basis of an arbitrary matroid,
let $P_i\succeq0$, and suppose
$Q=\E\sum_{i\in B}P_i$ has rank $d_0\ge1$. Interpret
$Q^{-1/2}$ on its range. If
$\|Q^{-1/2}P_iQ^{-1/2}\|\leq\rho$ on supported elements,
then for independent bases $B_1,\ldots,B_m$ and $0<\varepsilon<1$,
\[
 \Pr\!\left\{\frac1m\sum_{j=1}^m\sum_{i\in B_j}P_i
           \notin[(1-\varepsilon)Q,(1+\varepsilon)Q]\right\}
 \leq2d_0\exp\!\left(-\frac{m\varepsilon^2}{3\rho}\right).
\]
Here interval membership is in the positive semidefinite order.
If $Q=0$, every supported summand vanishes and the approximation is exact.
\end{corollary}
\begin{proof}
Write $r$ for the rank of the matroid. The joint basis law is the
law of a direct sum of $m$ matroids, of total degree $mr$. Give
its one-particle label $(j,i)$ the kernel
$B_{(j,i)}=rQ^{-1/2}P_iQ^{-1/2}$. Its normalized average is the
whitened matrix in the event, its mean is $I$, its kernel bound is
$R=r\rho$, and $c=mr$. Thus $c/R=m/\rho$ in
Theorem~\ref{thm:u-chernoff}. The standard scalar estimates for
the two Chernoff functions give the claim.
\end{proof}

Theorems~\ref{thm:u-laplace}--\ref{thm:u-bernstein} apply to
every matroid, without representability assumptions, and to
all homogeneous Lorentzian occupation laws.
The effective sample size for a kernel on $\ell$ particles is
$r/\ell$. The quadratic gap gives the larger value $a_r/a_\ell$
for simple matroids and any other law satisfying
\eqref{eq:hodge-gap}.

\section{Exact marginal and matrix H\"older regions}
\label{sec:weighted-marginals}

Selected-particle entropy controls the event that a coordinate is
present. Applying it also to the complement controls the event
that it is absent. Together they give the half-weight inequality
for full Bernoulli coordinates. Component factorization then gives
the whole admissible region, and a field concentrating on one
basis exchange proves every necessary pair constraint.
Section~\ref{sec:joint-entropy-rigidity} then determines fixed-field
sharpness and every equality density for the half vector.

Let $M$ be a matroid on $[n]$, of rank $r$, and let
\[
 f(z)=\sum_{B\in\mathcal B(M)}c_Bz^B,\qquad c_B>0.
\]
Assume that both $f$ and its complement polynomial
$f^*(z)=z^{[n]}f(z_1^{-1},\ldots,z_n^{-1})$ are Lorentzian.
For a positive field put
\[
 \mu_\lambda(B)=\frac{c_B\lambda^B}{f(\lambda)},
 \qquad \lambda_i>0.
\]
This includes product-field basis measures of every matroid and all
homogeneous strongly Rayleigh laws: a multiaffine stable polynomial
and its complement are stable, hence Lorentzian. The coefficient
weights are otherwise arbitrary subject to the two stated hypotheses.
Write $X_i=\mathbf1_{i\in B}$ and let $\mu_{\lambda,i}$ be its
Bernoulli reference law.

\begin{proposition}[Full Bernoulli marginals]\label{r4q:prop:half}
For every positive semidefinite matrix function $F$ on $\mathcal B(M)$,
\begin{equation}\label{r4q:eq:half}
 \frac12\sum_{i=1}^n
 \mathcal H_{\mu_{\lambda,i}}\bigl(\E[F\mid X_i]\bigr)
 \leq\mathcal H_{\mu_\lambda}(F).
\end{equation}
\end{proposition}
\begin{proof}
If $r=0$ or $r=n$, the law is a point mass and the assertion is
immediate. Otherwise put $p_i=\Pr(X_i=1)$, $G_i=\E[F\mid X_i=1]$,
$K_i=\E[F\mid X_i=0]$, and $\sigma=\E F$; omit missing conditioning
events. Applying Theorem~\ref{r4l:main} with $\ell=1$ and
multiplying by $r$ gives
\[
 \sum_ip_i\Tr[G_i\log G_i]-r\Tr[\sigma\log\sigma]
 \leq\mathcal H_\mu(F).
\]
Complementary bases have the law of $f^*$ at field $\lambda^{-1}$.
The same argument therefore gives
\[
 \sum_i(1-p_i)\Tr[K_i\log K_i]-(n-r)\Tr[\sigma\log\sigma]
 \leq\mathcal H_\mu(F).
\]
Their sum is the left side of \eqref{r4q:eq:half} before multiplication
by $1/2$. Constant coordinates cause no problem in these identities.
\end{proof}

\begin{lemma}[Factorization over matroid components]
\label{lem:lorentzian-component-factorization}
The coefficient law of a multiaffine Lorentzian polynomial factors
over the connected components of its support matroid, at every
positive field.
\end{lemma}
\begin{proof}
Choose adjacent bases in each of two distinct components and fix
bases in the other components. Differentiate by the common elements
and set all remaining coordinates except the four exchanged elements
to zero. The resulting quadratic is Lorentzian and has the form
$\sum_{i,j=1}^2 c_{ij}x_i y_j$, with $c_{ij}>0$.
Its Hessian has off-diagonal blocks $C$ and $C^{\mathsf T}$;
its number of positive eigenvalues is $\operatorname{rank}C$.
Hence $c_{11}c_{22}=c_{12}c_{21}$.
These rectangle identities say that the coefficient ratio along
a basis exchange in one component is unchanged by an exchange in
another. Each component's basis graph is connected. Fixing reference
bases and following exchange paths therefore factors the coefficients,
and positive fields preserve the factorization. Each factor inherits
Lorentzianity by differentiating in fixed bases of the other
components. When $f^*$ is also Lorentzian, the same argument gives
Lorentzianity of the complementary factors.
\end{proof}

Let $C_1,\dots,C_m$ be the nontrivial connected components of $M$;
loops and coloops are constant coordinates and may be omitted. Define
\[
 \mathcal W(M)=\{w\in\mathbb R_{\geq0}^n:
    w_i+w_j\leq1\text{ whenever }i\ne j\text{ lie in the same }C_a\}.
\]

\begin{theorem}[Sharp all-field quantum coefficient region]
\label{r4q:thm:weighted}
Assume $f$ and $f^*$ are Lorentzian as above. For $w\geq0$,
the following statements are equivalent.
\begin{enumerate}
\item $w\in\mathcal W(M)$.
\item For every positive external field, every finite quantum register
  $Q$, and every classical--quantum state supported on
  $\mathcal B(M)\otimes Q$,
  \begin{equation}\label{r4q:eq:weighted}
   \sum_iw_iD(\rho_{X_iQ}\Vert\mu_{\lambda,i}\otimes\rho_Q)
   \leq D(\rho_{BQ}\Vert\mu_\lambda\otimes\rho_Q).
  \end{equation}
\item Statement (2) holds with a one-dimensional quantum register.
\end{enumerate}
The same equivalence applies to every nonzero minor law, with the
connected components of its support matroid.
\end{theorem}
\begin{proof}
First suppose $M$ is connected with at least two elements. The half
vector satisfies \eqref{r4q:eq:weighted} by Proposition~\ref{r4q:prop:half},
and every coordinate unit vector does by data processing. If
$h=\max_iw_i\leq1/2$, the weight vector is dominated coordinatewise
by the half vector. Otherwise choose $j$ with $w_j=h>1/2$.
All other weights are at most $1-h$, and $h\leq1$. Therefore
\[
 w\leq(2h-1)e_j+(2-2h)\frac12\one.
\]
This is a convex combination of two admissible coefficient vectors.
Nonnegativity of relative entropy proves the claim by domination.

For several components, Lemma~\ref{lem:lorentzian-component-factorization}
makes $\mu_\lambda$ a product of its component laws. If
$B_a=B\cap C_a$, the complete product inequality
\begin{equation}\label{r4q:eq:ssa-product}
 \sum_aD(\rho_{B_aQ}\Vert\mu_{\lambda,a}\otimes\rho_Q)
 \leq D(\rho_{BQ}\Vert\mu_\lambda\otimes\rho_Q)
\end{equation}
follows from strong subadditivity~\cite{LiebRuskai}. Indeed, the difference of the right
and left sides is
$\sum_a H(B_a\mid Q)-H(B_1\cdots B_m\mid Q)\geq0$.
Apply the connected-component inequality to each $B_a$ and use
\eqref{r4q:eq:ssa-product}. This proves (1)$\Rightarrow$(2), and
(2)$\Rightarrow$(3) is immediate.

For necessity, two distinct elements $i,j$ belong to the same connected
component precisely when some circuit contains both. Extending that
circuit minus $j$ to a basis gives a basis $B$ such that
$B'=B\setminus\{i\}\cup\{j\}$ is also a basis. The segment joining
their indicator vectors is an edge of the matroid base polytope.
Expose this edge by a linear functional $q$ and take
$\lambda_k(t)=\exp(tq_k)$. The laws converge to a two-point law
on $\{B,B'\}$ assigning mass $\pi=c_B/(c_B+c_{B'})\in(0,1)$
to $B$. Test with $\nu=\delta_B$. The joint divergence and
the two changing coordinate divergences tend to $-\log\pi>0$;
all other coordinate divergences tend to zero. Hence
$w_i+w_j\leq1$. Differentiation and setting variables to zero
preserve both hypotheses, since complementation interchanges these
operations. The same proof therefore applies to every nonzero minor.
\end{proof}

\begin{remark}[The all-field matrix entropy ribbon]
After removing deterministic coordinates, define the matrix
$x\log x$ entropy ribbon of a law $\mu$ by the weights
$w\in[0,1]^n$ satisfying
$\sum_iw_i\mathcal H_{\mu_i}(\E[F\mid X_i])
\leq\mathcal H_\mu(F)$ in every finite matrix dimension.
In the terminology of Wang--Gohari~\cite{WangGohari2026},
Theorem~\ref{r4q:thm:weighted} determines its all-field intersection:
\[
 \bigcap_{\lambda>0}\mathcal R_{x\log x}(\mu_\lambda)
 =\mathcal W(M).
\]
The component partition alone determines this intersection.
For a single fixed field the ribbon can be larger.
\end{remark}

\begin{corollary}[The exact matrix H\"older region]
\label{cor:matrix-holder-region}
Under the generating-and-complement hypotheses of
Theorem~\ref{r4q:thm:weighted}, fix $w_i\ge0$.
For positive definite $A_{i,0},A_{i,1}$ put
$M_i=\E_{\mu_i}A_{i,X_i}$.
The inequality
\[
 \E_\mu\Tr\exp\!\left(K+\sum_iw_i\log A_{i,X_i}\right)
 \le\Tr\exp\!\left(K+\sum_iw_i\log M_i\right)
\]
holds for all such matrices, in every matrix dimension, for every
Hermitian $K$ and every positive field, if and only if
$w\in\mathcal W(M)$.
\end{corollary}
\begin{proof}
In the Gibbs variational formula, let $F$ have normalized trace
mean $N$ and put $J_i=\E[F\mid X_i]$.
The weighted marginal entropy bound gives
$\mathcal H_\mu(F)\ge\sum_iw_i\mathcal H_{\mu_i}(J_i)$.
Joint convexity gives
\[
 \E_{\mu_i}\Tr[J_i(\log J_i-\log A_i)]
 \ge\Tr[N(\log N-\log M_i)].
\]
Substitution in the variational formula and optimization over $N$
prove the trace comparison. For the converse take $F>0$, normalized
by $\Tr N=1$, and set $A_i=J_i$ and
$K=-(\sum_iw_i-1)\log N$. The resulting normalized trace bound,
followed by Klein's inequality exactly as in
Theorem~\ref{thm:channel-duality}, gives the weighted marginal
entropy bound. Regularization includes semidefinite inputs.
Its exact region is Theorem~\ref{r4q:thm:weighted}.
\end{proof}

\section{Marginal equality, reciprocity, and sharp fields}
\label{sec:joint-entropy-rigidity}

The half vector of Section~\ref{sec:weighted-marginals} is sharp at a fixed
field precisely when the law and its complement have nonzero compatible
translation directions. On connected support, exact integration gives two Bernoulli product marginals,
the sharp-field locus, and every matrix extremizer.

\subsection{A joint kernel and quantitative strictness}

Absorb a fixed positive external field and write
\[
 f(z)=\sum_{B\in\mathcal B(M)}\mu(B)z^B,\qquad f(\mathbf1)=1.
\]
Assume that \(f\) and \(f^*(z)=z^{[n]}f(z^{-1})\) are Lorentzian.
Remove deterministic coordinates; the point-mass case is trivial.
Thus \(r=\operatorname{rank}M\), \(s=n-r\) are positive and
\(p_i=\Pr(i\in B)\), \(q_i=1-p_i\) belong to \((0,1)\).
The selected-particle distribution assigns probability \(p_i/r\)
to label \(i\). Put
\[
 A=\nabla^2f(\mathbf1),\quad A^*=\nabla^2f^*(\mathbf1),\quad
 \mathsf D=\operatorname{diag}(p_i/q_i),\qquad
 \mathcal J=\{h:D_hf=0,\ D_{-\mathsf D h}f^*=0\}.
\]
Lemma~\ref{lem:translation-kernel} identifies the translation space
with the Hessian kernel in degree at least two.  In degree one it is
the mean-zero hyperplane.  Consequently, uniformly in \(r,s\ge1\),
\[
 \mathcal J=\{h:p^{\mathsf T}h=0,\ Ah=0,\ A^*\mathsf D h=0\}.
\]
For \(p^{\mathsf T}h=0\) define
\[
 Q(h)=-\frac{h^{\mathsf T}Ah}{r}
       -\frac{(\mathsf D h)^{\mathsf T}A^*(\mathsf D h)}s,\qquad
 V(h)=\sum_i\frac{p_i}{q_i}h_i^2,\qquad
 \gamma=\min_{h\ne0,\ p^{\mathsf T}h=0}\frac{Q(h)}{V(h)}.
\]
Both summands of \(Q\) are nonnegative by the Lorentzian signature,
since \(q^{\mathsf T}\mathsf D h=p^{\mathsf T}h=0\).
Its kernel is \(\mathcal J\); hence \(\gamma>0\) exactly when
\(\mathcal J=0\).

For a positive semidefinite matrix function \(F\) on a register
\(\mathcal Q\), put \(\sigma=\E_\mu F\) and write
\[
 \mathcal M_\mu(F)=\sum_i
 \mathcal H_{\operatorname{Bern}(p_i)}(\E[F\mid X_i]),
 \qquad X_i=\one_{i\in B}.
\]

\begin{theorem}[Exact sharpness criterion and quantitative strictness]
\label{thm:joint-kernel-coefficient}
In every finite matrix dimension,
\[
 \frac12\mathcal M_\mu(F)\le\frac1{1+\gamma}\mathcal H_\mu(F).
\]
The optimal complete coefficient \(\beta_\mu\) in
\(\mathcal M_\mu(F)/2\le\beta_\mu\mathcal H_\mu(F)\), taken over all
finite matrix dimensions, equals one
if and only if \(\mathcal J\ne0\).  When \(\mathcal J=0\), this gives the
dimension-independent estimate \(\beta_\mu\le(1+\gamma)^{-1}<1\).
The same bound and coefficient-one criterion hold when every entropy
is replaced by $\mathcal H_{\nu,a}(F)=\Tr\E_\nu F^a-\Tr(\E_\nu F)^a$,
$1<a\le2$, also for pointwise trace-one inputs.
\end{theorem}

\begin{proof}
Compress to \(\operatorname{supp}\sigma\).  For \(\sigma>0\) use the
Bogoliubov--Kubo--Mori (BKM) inner product
\[
 \langle U,V\rangle_\sigma
 =\operatorname{Tr}\!\left[
 U\int_0^\infty(\sigma+tI)^{-1}V(\sigma+tI)^{-1}\,dt\right].
\]
Put \(G_i=\mathbb E[F\mid X_i=1]\), \(K_i=\mathbb E[F\mid X_i=0]\),
\(Z_i=G_i-\sigma\), and \(Y_i=K_i-\sigma=-(p_i/q_i)Z_i\).
Write \(D_g(U\Vert V)=\operatorname{Tr}[U(\log U-\log V)-U+V]\).
The profile identity of Theorem~\ref{thm:entropy-profile} and the quadratic remainder in
Theorem~\ref{thm:all-level-gap} give
\begin{align*}
 \delta_+&:=\mathcal H_\mu(F)-\sum_i p_iD_g(G_i\Vert\sigma)
       \ge-\frac1r\sum_{ij}A_{ij}\langle Z_i,Z_j\rangle_\sigma,\\
 \delta_-&:=\mathcal H_\mu(F)-\sum_i q_iD_g(K_i\Vert\sigma)
       \ge-\frac1s\sum_{ij}A^*_{ij}\langle Y_i,Y_j\rangle_\sigma .
\end{align*}
Indeed \(H_r-rH_1\ge(r-1)(H_2-2H_1)\), and the ordered pair
matrix is \(A/[r(r-1)]\).  For \(r=1\), the first selected channel
is the identity and its remainder is zero; likewise when \(s=1\).
Since \(\sum_i p_iZ_i=0\), expansion in a BKM orthonormal basis yields
\[
 2\mathcal H_\mu(F)-\mathcal M_\mu(F)
 =\delta_++\delta_-\ge
 \gamma\sum_i\frac{p_i}{q_i}\|Z_i\|_\sigma^2.
\]
The tangent inequality for \(\log\), used in Theorem~\ref{thm:all-level-gap},
gives \(D_g(U\Vert\sigma)\le\|U-\sigma\|_\sigma^2\).  Therefore
\[
 \mathcal M_\mu(F)\le
 \sum_i\bigl(p_i\|Z_i\|_\sigma^2+q_i\|Y_i\|_\sigma^2\bigr)
 =\sum_i\frac{p_i}{q_i}\|Z_i\|_\sigma^2.
\]
This proves the bound; regularization includes singular inputs.
The scalar equality family constructed next proves logarithmic sharpness.
For trace powers, apply this bound to $F+tI$ and integrate the shifted-log
formula from Theorem~\ref{thm:complete-power-contraction}; continuity gives
$a=2$. If $0\ne h\in\mathcal J$, the tangent
$u(B)=\sum_{i\in B}h_i$ of that equality family has half-marginal
variance ratio one. Thus the trace-one qubit perturbation
$F_\varepsilon=\frac12\diag(1+\varepsilon u,1-\varepsilon u)$ has
power-entropy ratio tending to one, with exact equality for $a=2$.
If $\mathcal J=0$, the strict bound applies.
\end{proof}

\subsection{Exact integration of the joint directions}

\begin{theorem}[A pointwise primal--dual identity]
\label{thm:joint-tilt-identity}
The normalized scalar equality densities for the logarithmic half-weight inequality
are exactly
\[
 F_B=\prod_{i\in B}(1+h_i)
     =\prod_{i\notin B}\left(1-\frac{p_i}{q_i}h_i\right),
 \qquad h\in\mathcal J,\quad -1\le h_i\le q_i/p_i .
\]
In particular the two products agree on every supported basis.
\end{theorem}

\begin{proof}
First let \(a_i=1+h_i>0\), \(b_i=1-(p_i/q_i)h_i>0\).
Theorem~\ref{thm:lorentzian-equality} gives
normalized densities \(F_B^+=a^B\), \(F_B^-=b^{B^c}\), with entropies
\[
 S_+=\sum_i p_i a_i\log a_i,\qquad
 S_-=\sum_i q_i b_i\log b_i.
\]
Their present and absent conditional means are \(a_i,b_i\) in both
cases, because \(p_i a_i+q_i b_i=1\).
Applying complementary deletion to \(F^+\) and original deletion
to \(F^-\) gives \(S_-\le S_+\le S_-\).  Moreover,
\[
 D(\mu F^+\Vert\mu F^-)
 =S_+-\mathbb E_\mu[F^+\log F^-]=S_+-S_-=0.
\]
Thus the densities agree pointwise.  Replace \(h\) by \(th\) and
let \(t\uparrow1\) for boundary parameters.
Conversely, half equality forces both nonnegative selected deficits
to vanish.  Their scalar equality classifications and the common
Bernoulli mean give precisely \(a,b,h\) above.
In degree one, the same description is just the mean-zero condition.

For \(0\ne h\in\mathcal J\), small nonzero multiples give nonconstant
positive densities: the derivative \(u(B)=\sum_{i\in B}h_i\) has
\(\operatorname{Var}u=\sum_i p_i h_i^2>0\), since \(Ah=0\).
This completes Theorem~\ref{thm:joint-kernel-coefficient}.
\end{proof}

\subsection{Automatic reciprocity and connected rigidity}

The next algebraic statement uses only nonnegative coefficients and
the translation space \(\mathcal N_f=\{h:D_hf=0\}\).

\begin{theorem}[Automatic reciprocity and independent halves]
\label{efg:reciprocity}
Let \(f\ne0\) be homogeneous and multiaffine of degree \(m\ge1\)
on \(2m\) variables, with nonnegative coefficients.  Suppose
\(h\in\mathcal N_f\) and \(h_i\ne0\) for every \(i\).
Then its two sign sets \(A,A'\) each have size \(m\), and both are
supported monomials.  With \(\lambda_i=|h_i|\), the normalized polynomial
\(g(z)=f(\lambda z)/f(\lambda)\) satisfies
\[
 g^*=g,\qquad
 g(z_A,\mathbf1_{A'})=2^{-m}\prod_{i\in A}(1+z_i),\qquad
 g(\mathbf1_A,z_{A'})=2^{-m}\prod_{j\in A'}(1+z_j).
\]
Equivalently,
\begin{equation}\label{efg:reciprocity-formula}
 f^*(z)=\left(\prod_i|h_i|\right)^{-1}
                 f(h_1^2z_1,\ldots,h_{2m}^2z_{2m}).
\end{equation}
In particular, if \(f\) is Lorentzian, so is \(f^*\).
Each half of the balanced reference law consists of independent
Bernoulli\((1/2)\) coordinates.
\end{theorem}

\begin{proof}
Put \(\epsilon_i=\operatorname{sign}(h_i)\).
Then \(D_\epsilon g=0\) and
\(g(\mathbf1+\epsilon)=g(\mathbf1-\epsilon)=g(\mathbf1)=1\).
Nonnegativity therefore makes each sign set contain a supported
\(m\)-subset.  There are only \(2m\) labels, so each set has size \(m\)
and its coefficient in \(g\) is \(2^{-m}\).
Translation now gives
\[
 g(z_A,\mathbf1_{A'})
 =g(z_A+\mathbf1_A,0_{A'})
 =2^{-m}\prod_{i\in A}(1+z_i),
\]
and similarly for \(A'\).

For reciprocity, the sign change places $g$ at middle weight in the
Boolean $\mathfrak{sl}_2$ representation on $2m$ variables.
A middle-weight vector killed by lowering is invariant; the Weyl
element then acts by signed complementation.
Put \(\widehat g(w)=g((\epsilon_iw_i)_i)\).
On the multiaffine vector space consider
\[
 \mathsf F=\sum_i\partial_{w_i},\qquad
 \mathsf E=\sum_iw_i(1-w_i\partial_{w_i}),\qquad
 [\mathsf E,\mathsf F]=2\sum_iw_i\partial_{w_i}-2m.
\]
In the inner product making monomials orthonormal,
\(\mathsf E^*=\mathsf F\).  Since \(\widehat g\) has degree \(m\)
and \(\mathsf F\widehat g=0\), the commutator gives
\(\|\mathsf E\widehat g\|^2
=\langle\mathsf F\mathsf E\widehat g,\widehat g\rangle=0\).
Thus both raising and lowering annihilate \(\widehat g\).
On one variable,
\(\exp(\mathsf F)\exp(-\mathsf E)\exp(\mathsf F)\) sends
\(1\mapsto-w\) and \(w\mapsto1\).
On degree \(m\) in \(2m\) labels its tensor action is
\((-1)^m\) times coefficient complementation.  Hence
\(\widehat g^*=(-1)^m\widehat g\).
Exactly \(m\) signs \(\epsilon_i\) are negative, so
\(\epsilon^{S^c}=(-1)^m\epsilon^S\) for every \(m\)-subset \(S\).
The signs cancel, proving \(g^*=g\).
Substituting in this identity gives \eqref{efg:reciprocity-formula};
positive diagonal substitutions preserve Lorentzianity.
\end{proof}

The size \(2m\) is minimal: the same sign argument shows that a
translation with no zero coordinate in degree \(m\) requires at least
\(m\) labels of each sign.

The coefficient law factors over matroid components by Lemma~\ref{lem:lorentzian-component-factorization}.  Translation spaces, and hence \(\mathcal J\),
split over these components: each logarithmic directional derivative
depends on its own variables and has homogeneity \(-1\).

\begin{theorem}[Two bases and independent coordinate halves]
\label{thm:joint-connected-rigidity}
For a nontrivial connected component \(C\), of rank \(m\),
\(\dim\mathcal J_C\le1\).  If it is nonzero, there are two bases
\(C=A\sqcup A'\), each of size \(m\), and \(\pi\in(0,1)\), such that
\[
 p_i=\pi\ (i\in A),\qquad p_i=1-\pi\ (i\in A'),\qquad
 \mathcal J_C=\mathbb R h,\quad
 h_i=\begin{cases}\pi^{-1},&i\in A,\\-(1-\pi)^{-1},&i\in A'.\end{cases}
\]
Within \(A\) the coordinates are independent Bernoulli\((\pi)\);
within \(A'\) they are independent Bernoulli\((1-\pi)\).
In particular \(|B\cap A|\) is \(\operatorname{Binomial}(m,\pi)\).
For \(m\ge2\), the component is simple and cosimple.
\end{theorem}

\begin{proof}
Choose \(0\ne h\in\mathcal J_C\) and put \(k_i=-p_i h_i/q_i\).
Any two elements \(i,j\) of a connected matroid can be exchanged
in some basis.  Compare the pointwise identities of
Theorem~\ref{thm:joint-tilt-identity} on those bases, at parameter \(th\):
\[
 (1+th_i)(1+tk_i)=(1+th_j)(1+tk_j).
\]
All pairs \((h_i,k_i)\) have the same sum and product.  A zero coordinate
would force every coordinate to vanish.  Otherwise their common
unordered pair is \(\{u,v\}\), with \(u>0>v\), and the two sign classes
have probabilities \(\pi=-v/(u-v)\) and \(1-\pi\).
The identities \(\sum_i p_i h_i=0\) and \(\sum_i p_i=m\) show that
both classes have size \(m\).  Evaluation at any coordinate is
injective on \(\mathcal J_C\), proving its dimension bound.

After the positive diagonal rescaling to the translation
\((\one_A,-\one_{A'})\), Theorem~\ref{efg:reciprocity} makes both
halves bases and their coordinates independent. Rescaling back gives
the stated probabilities and binomial count.
A parallel pair must cross the halves.  Its probabilities sum to one,
so every basis contains exactly one member.  That pair and its
complement have ranks \(1\) and \(m-1\), making the pair a direct
component.  Connectedness excludes this for \(m\ge2\).
The same argument in the dual proves cosimplicity.
\end{proof}

\subsection{Sharp fields and product components}

For a positive field \(\lambda\), write \(\beta_\lambda\) for the optimal
complete coefficient in \(\mathcal M_{\mu_\lambda}(F)/2\le
\beta_\lambda\mathcal H_{\mu_\lambda}(F)\).

\begin{theorem}[The complete sharp-field locus]
\label{efg:sharp-fields}
For a connected rank-\(m\) law with \(f,f^*\) Lorentzian, \(\beta_\lambda=1\) at some
positive field if and only if
\[
 n=2m\quad\text{and}\quad
 \mathcal N_f\text{ contains a vector }h\text{ with every }h_i\ne0.
\]
When this holds, \(\mathcal N_f=\mathbb Rh\), and all sharp fields are
exactly
\[
 \lambda_i=a|h_i|\quad(i\in A),\qquad
 \lambda_i=b|h_i|\quad(i\in A'),\qquad a,b>0.
\]
The partition is unique up to exchanging its halves.
Modulo common scaling of the field, the sharp locus is one line in
logarithmic field coordinates.  The sufficient condition itself makes
the hypothesis that \(f^*\) is Lorentzian automatic.
\end{theorem}

\begin{proof}
For necessity use the joint-kernel criterion and connected rigidity
of Theorems~\ref{thm:joint-kernel-coefficient} and~\ref{thm:joint-connected-rigidity}: a sharp field forces
\(n=2m\) and a joint translation with no zero coordinate.  Undoing
the field gives the stated \(h\).
Conversely, Theorem~\ref{efg:reciprocity} gives a self-complementary
balanced law with all inclusion probabilities \(1/2\).
Its translations \(\epsilon,-\epsilon\) satisfy the joint criterion,
so this field is sharp.

To prove uniqueness, let \(v\in\mathcal N_f\).  For small real \(t\),
\(h+tv\) has the same signs as \(h\), and the fields \(|h+tv|\)
all have inclusion vector \((1/2,\ldots,1/2)\).
The exponential-family mean map is injective modulo common scaling
on connected matroid support.  Indeed, if two fields have the same
inclusion vector, then
\[
 D(\mu_\lambda\Vert\mu_\nu)+D(\mu_\nu\Vert\mu_\lambda)
 =\sum_i(p_i(\lambda)-p_i(\nu))\log(\lambda_i/\nu_i)=0.
\]
Their laws coincide, so \((\lambda/\nu)^B\) is constant on all bases.
Connected basis exchanges force all \(\lambda_i/\nu_i\) to agree.
Consequently \(|h+tv|\) is proportional to \(|h|\), and \(v\in\mathbb Rh\).

At an arbitrary sharp field, the primal translation is proportional
to \((h_i/\lambda_i)_i\).  Connected rigidity requires one common
positive value and one common negative value, giving the displayed
field locus.
Conversely, at such a field the present probabilities are
\(a/(a+b)\) on \(A\) and \(b/(a+b)\) on \(A'\).
The primal direction is \((1/a,-1/b)\); the complementary direction
required by the joint criterion is
\[
 (-1/b,1/a)=-(ab)^{-1}(a,-b).
\]
Reciprocity gives \((a,-b)\) as a complementary translation.
Thus every displayed field is sharp.
\end{proof}

Write \(G=2^m g\), \(p=a/(a+b)\), and \(q=1-p\).
We call \(G\) balanced: it is multiaffine and homogeneous of degree
\(m\) on two \(m\)-element halves, has nonnegative coefficients,
pure-half coefficients one, and \(D_{(\mathbf1,-\mathbf1)}G=0\).
The whole sharp family is \(f_p(z)=G(pz_A,qz_{A'})\), with
\[
 f_p(z_A,\mathbf1)=\prod_{i\in A}(pz_i+q),\qquad
 f_p(\mathbf1,z_{A'})=\prod_{j\in A'}(qz_j+p).
\]
The count \(|B\cap A|\) is \(\operatorname{Binomial}(m,p)\);
the conditional law given that count is independent of \(p\).
In particular
\[
 D(\mu_p\Vert\mu_t)
   =mD(\operatorname{Bern}(p)\Vert\operatorname{Bern}(t)).
\]

These product marginals also explain the half bound directly.
The complete product inequality and data processing give
\[
 \sum_{i\in A}\mathcal H_{\mu_i}(\E[F\mid X_i])
 \le\mathcal H_{\mu_A}(\E[F\mid X_A])
 \le\mathcal H_\mu(F).
\]
Apply the same argument to \(A'\) and average.
This uses only the product strong subadditivity
\eqref{r4q:eq:ssa-product} and the two product marginals.

\begin{theorem}[Exact component tensorization and generic strictness]
\label{thm:marginal-component-coefficient}
Let \(\mu=\bigotimes_C\mu_C\) be the factorization over the nontrivial
matroid components. Then
\[
 \beta_\mu=\max_C\beta_{\mu_C}.
\]
The same value is obtained when inputs are restricted to pointwise
trace-one matrices.

Consequently the global set of sharp positive fields is the union
of the component cylinders specified by Theorem~\ref{efg:sharp-fields},
with the fields outside the selected component arbitrary.
If a binary direct component exists, every positive field is sharp.
Otherwise each nonempty cylinder has codimension \(|C|-2\ge2\)
in logarithmic field coordinates. The complete half bound is therefore
strict on an open dense set of full measure, uniformly in matrix
dimension; the improvement is uniform on compact subsets of that set.
\end{theorem}

\begin{proof}
For \(F_C=\E[F\mid B\cap C]\), apply the component inequalities and
\eqref{r4q:eq:ssa-product}:
\[
 \frac12\mathcal M_\mu(F)
 \le\sum_C\beta_{\mu_C}\mathcal H_{\mu_C}(F_C)
 \le\left(\max_C\beta_{\mu_C}\right)\mathcal H_\mu(F).
\]
Inputs depending on one component prove the reverse inequality for the
optimal coefficients.
To compare unrestricted inputs with trace-one encodings, enlarge the
register by one dimension and use the flag extension
\[
 \widehat F_\varepsilon(B)
   =\operatorname{diag}\bigl(\varepsilon F(B),
                       1-\varepsilon\operatorname{Tr}F(B)\bigr).
\]
For small positive \(\varepsilon\) this is a pointwise density matrix.
The first block contributes exactly
\(\varepsilon\mathcal H_\mu(F)\), and the scalar block contributes
\(O(\varepsilon^2)\); the same expansion holds after every marginal
channel. Thus every nonzero entropy ratio is approached by trace-one
ratios. The reverse comparison is immediate by restriction.
The suprema agree; exact attainment is treated below.

Lemma~\ref{lem:lorentzian-component-factorization} and
Theorem~\ref{efg:sharp-fields} give the stated cylinders.
Binary cylinders fill their field space; other eligible components
have at least four labels. A finite union of proper affine subspaces
in logarithmic coordinates is closed and null.
Continuity of \(\gamma>0\) off this union gives compact uniformity.
\end{proof}

\Needspace{16\baselineskip}
\subsection{Matrix equality and exact attainment}

\begin{theorem}[Only binary components carry noncommuting factors]
\label{thm:half-matrix-form}
Let \(\mathcal S\) be the nonbinary connected components with
\(\mathcal J_C\ne0\), and let \(\mathcal U\) be the binary
\(U_{1,2}\) direct components.
All logarithmic half-weight equality densities with
\(\operatorname{Tr}\mathbb E_\mu F=1\) have exactly the form
\begin{equation}\label{eq:half-matrix-normal-form}
 F_B=\bigoplus_z w_z
       \left(\prod_{C\in\mathcal S}\prod_{i\in B\cap C}(1+t_{C,z}h_i^C)\right)
       \tau_0^z\otimes
       \bigotimes_{D\in\mathcal U}A^{D,z}_{Y_D(B)}.
\end{equation}
Here \(w_z>0\), \(\sum_z w_z=1\); \(Y_D\) is the binary choice with
reference law \(\pi_D\);
\(\sum_y\pi_D(y)A_y^{D,z}=\tau_D^z\);
the \(A_y^{D,z}\) are positive semidefinite and the \(\tau_0^z,\tau_D^z\)
are faithful density matrices on their respective tensor factors.
For each \(C\in\mathcal S\), take \(h^C\) and \(\pi_C\) from
Theorem~\ref{thm:joint-connected-rigidity}; then
\(-\pi_C\le t_{C,z}\le1-\pi_C\).  Empty products are
permitted, with \(0^0=1\).
\end{theorem}

\begin{proof}
Suppose first that the global rank is at least two.
Primal and dual equality give
\((G_i-\sigma)_i\in\mathcal J\otimes\operatorname{Herm}(\mathcal Q)\),
by \eqref{eq:quantum-null-directions}; for a degree-one complement,
its null identity is the centered mean relation.
Apply Theorem~\ref{thm:universal-lorentzian-matrix-equality}.
Outside surviving components \(G_i=\sigma\).
On each surviving nonbinary component, Theorem~
\ref{thm:joint-connected-rigidity} makes every parallel class a
singleton.  Its matrices are therefore scalar multiples of \(\sigma\)
in each central block, with parameters on the one-dimensional joint
line.  On a binary component the class mean is \(\sigma\), leaving
exactly its two matrix factors.  Every remaining class factor is
constant. This is \eqref{eq:half-matrix-normal-form}.
Conversely, entropy is additive over normalized independent tensor
factors and affine over these orthogonal blocks.  Each scalar factor
has half equality, and the two coordinates of a binary component
both reveal its entire choice. Hence \eqref{eq:half-matrix-normal-form}
has equality.
For global rank one, \(G_i=F_{\{i\}}\).
Joint nullity therefore makes \(F\) constant unless the component is
a binary pair, where the assertion allows an arbitrary two-state ensemble.
\end{proof}

\begin{corollary}[Binary rigidity and unattained Holevo sharpness]
\label{cor:half-fixed-trace}
For positive product-field matroid laws, \(\beta_\mu=1\) exactly
when a binary direct component exists.
For arbitrary coefficients under the two Lorentzian hypotheses,
an informative equality encoding satisfying
\(\operatorname{Tr}F_B=1\) exists exactly when a binary component exists.
Nevertheless the optimal coefficient restricted to such encodings
equals one exactly when \(\mathcal J\ne0\).
Thus, if \(\mathcal J\ne0\) and there is no binary component,
scalar densities attain the bound while informative trace-one
encodings approach it without attaining it.
\end{corollary}

\begin{proof}
For product-field laws, the translation classification in Theorem~\ref{thm:matroid-equality}, using the Hessian strictness of~\cite{MNY}, gives only
weighted parallel-class contrasts.  Its dual
gives series-class contrasts.  A circuit \(\{i,j\}\) and cocircuit
\(\{i,k\}\) cannot meet in exactly one element, so a common nonzero
coordinate must belong to a binary direct component.
Conversely such a component has a joint contrast.

Under pointwise trace one, Corollary~\ref{cor:fixed-marginal-no-quotient-tilt} removes every
quotient scalar tilt in \eqref{eq:half-matrix-normal-form};
only binary factors remain.
Any nonconstant trace-one ensemble on a binary component suffices.
For \(0\ne h\in\mathcal J\), put \(u(B)=\sum_{i\in B}h_i\) and
\[
 F_\varepsilon(B)=\tfrac12
   \operatorname{diag}(1+\varepsilon u(B),1-\varepsilon u(B)).
\]
Differentiating the pointwise identity gives
\(u(B)=\sum_{i\notin B}(-\mathsf D h)_i\).
Thus \(u\) is simultaneously a primal and dual variance equality
mode.  The second-order entropy expansion makes the half-marginal
ratio tend to one.  When \(\mathcal J=0\), the strict complete bound
already applies to these encodings.
\end{proof}

For product-field matroid laws, marginal entropy equality and the
unit-edge eigenvalue-two space thus single out the same binary direct
components (Theorem~\ref{thm:discovery-spectral-separation}).

\subsection{Determinantal models and the cubic classification}

For any \(m\times m\) unitary \(U\), put
\[
 G_U(x,y)=\det\!\bigl(\diag(x)+U\diag(y)U^*\bigr),\qquad
 f_p(z)=G_U(pz_A,(1-p)z_{A'}),\quad 0<p<1.
\]
Cauchy--Binet gives nonnegative coefficients. Upper-half-plane variables
make the matrix's imaginary part positive definite, so \(G_U\) is
stable and hence Lorentzian~\cite{BrandenHuh}.
Its full-support translation is \((\one_A,-\one_{A'})\), and
Theorem~\ref{efg:reciprocity} gives complement reciprocity.
The law of \(f_p\) is the projection determinantal law associated with
the row space of \([\sqrt p\,I,\sqrt{1-p}\,U]\).
For generic \(U\), every square minor is nonzero, giving support
\(U_{m,2m}\), simple and cosimple for \(m\ge2\).
Corollary~\ref{cor:half-fixed-trace} gives scalar half equality and
sharp but unattained informative trace-one ratios.

For \(m=3\), self-complementarity and translation force
\begin{equation}\label{efg:cubic-polynomial}
 G_M(x,y)=x_1x_2x_3+y_1y_2y_3+
 \sum_{i,j=1}^3M_{ij}
   \bigl(x_{[3]\setminus\{i\}}y_j+x_i y_{[3]\setminus\{j\}}\bigr),
\end{equation}
where \(M\) is doubly stochastic.
Indeed complementation pairs the mixed coefficients; the coefficients
of \(x_{[3]\setminus\{i\}}\) and \(y_{[3]\setminus\{j\}}\) in
\((\sum\partial_x-\sum\partial_y)G_M=0\) give its row and column sums.
Conversely those sums give the translation identity.

\begin{theorem}[Complete cubic classification]
\label{efg:cubic-classification}
For a \(3\times3\) doubly stochastic matrix \(M\), the following
conditions are equivalent:
\begin{enumerate}
\item \(G_M\) is Lorentzian.
\item The Hessian of the single quadratic
  \(\partial_{x_1}G_M\) has at most one positive eigenvalue.
\item The three numbers \(\sqrt{M_{2j}M_{3j}}\), \(1\le j\le3\),
  are the side lengths of a possibly degenerate triangle.
\item \(M_{ij}=|U_{ij}|^2\) for a unitary matrix \(U\).
\item \(G_M=\det(\diag(x)+U\diag(y)U^*)\) for a unitary \(U\).
\item \(G_M\) is real stable.
\end{enumerate}
Thus every connected rank-three sharp law is a two-orthonormal-bases
projection determinantal law after positive rescaling.
\end{theorem}

\begin{proof}
Lorentzianity implies (2), and (5) implies stability by positivity of
the imaginary part of the determinant matrix.
Stability implies Lorentzianity by \cite{BrandenHuh}.
We prove the remaining implications explicitly.

Set \(a_j=M_{2j}\), \(b_j=M_{3j}\), and
\(\ell_j=\sqrt{a_jb_j}\).
In the Hessian of \(\partial_{x_1}G_M\), omitting its zero \(x_1\) row and ordered by
\(x_2,x_3,y_1,y_2,y_3\), the first block is
\(\left(\begin{smallmatrix}0&1\\1&0\end{smallmatrix}\right)\).
Its Schur complement \(S\) satisfies
\[
 S\mathbf1=0,\qquad S_{jj}=-2\ell_j^2,\qquad
 S_{ij}=\ell_i^2+\ell_j^2-\ell_k^2
       \quad(\{i,j,k\}=\{1,2,3\}).
\]
Every principal \(2\times2\) minor of \(S\) equals
\begin{align*}
 \Delta
 &=2\sum_{i<j}\ell_i^2\ell_j^2-\sum_i\ell_i^4\\
 &=(\ell_1+\ell_2+\ell_3)(-\ell_1+\ell_2+\ell_3)
   (\ell_1-\ell_2+\ell_3)(\ell_1+\ell_2-\ell_3).
\end{align*}
Since \(\det S=0\), the principal-minor criterion gives
\(S\preceq0\) exactly when \(\Delta\ge0\).
The saddle block has inertia \((+,-)\), so Schur-complement inertia
proves (2)\(\Leftrightarrow\)(3).
When the triangle exists, \(\Delta\) is sixteen times its squared area.

Choose phases satisfying
\(\sum_j\sqrt{a_jb_j}\,e^{i\theta_j}=0\).
The unit rows \((\sqrt{a_j})_j\) and
\((\sqrt{b_j}e^{i\theta_j})_j\) are orthogonal.
Complete them to a unitary \(U\), assigning them rows two and three.
Column norms force \(|U_{1j}|^2=1-a_j-b_j=M_{1j}\).
This is the classical unistochastic triangle construction
\cite[Section~3, equations~(21)--(25)]{Bengtsson}; it also covers
degenerate triangles and zero entries.

Finally Cauchy--Binet gives coefficient \(|U_{ij}|^2\) for
\(x_{[3]\setminus\{i\}}y_j\).  Complementary minors of a unitary give
the same coefficient for \(x_i y_{[3]\setminus\{j\}}\).
The pure coefficients are one, proving (4)\(\Rightarrow\)(5).
The converse follows from the first coefficients.
\end{proof}

Real orthogonal representations correspond exactly to degenerate
triangles. In rank two, the matrix
\(\left(\begin{smallmatrix}t&1-t\\1-t&t\end{smallmatrix}\right)\)
is realized by a rotation; rank one is immediate. The determinantal
classification therefore holds through rank three.

For \(\sigma\in S_m\), put
\(P_\sigma(x,y)=\prod_i(x_i+y_{\sigma(i)})\).
Through rank three, every nonnegative balanced polynomial is a convex
combination of these matching products: its coefficient matrix is
doubly stochastic, Birkhoff's theorem expresses that matrix as a convex
combination of permutation matrices, and the polynomial depends
affinely on it \cite[Section~2]{Bengtsson}.
Br\"and\'en obtained obstructions to determinantal representability
from matroid representability \cite{BrandenDeterminantal}.
For the symmetric family below, two-by-two minors force an impossible
Gram matrix.

\begin{proposition}[Two incomparable models from rank four]
\label{thm:balanced-incomparability}
For every \(m\ge4\), connected real stable sharp laws on \(U_{m,2m}\)
that admit positive matching decompositions and those that admit
homogeneous positive-semidefinite determinantal representations form
incomparable classes. In particular, the determinantal classification
ends exactly at rank three.
\end{proposition}

\begin{proof}
One direction has the explicit witness
\[
 G_{\mathrm{sym}}(x,y)
 =\frac1{m!}\operatorname{per}\bigl[(x_i+y_j)_{i,j=1}^m\bigr]
 =\sum_{k=0}^m\frac{e_{m-k}(x)e_k(y)}{\binom mk}.
\]
Here \(\operatorname{per}\) denotes the permanent, equal to
\(\sum_{\sigma\in S_m}P_\sigma(x,y)\). Counting permutations sending a
fixed \(k\)-subset to another gives the second identity.
This polynomial is the separate polarization of
\((x+y)^m\), so it is real stable
\cite[Section~2.2, Proposition~2.4]{BorceaBrandenPolarization}.
It has every degree-\(m\) squarefree monomial with positive coefficient
and translation \((\one,-\one)\), hence gives connected sharp laws.
To sample \(G_{\mathrm{sym}}(px,(1-p)y)\), draw
\(K\sim\operatorname{Binomial}(m,p)\), then independently choose uniform
subsets of sizes \(K\) and \(m-K\) from the two halves, conditional on \(K\).

A two-basis representation
\(G_{\mathrm{sym}}=\det(\diag(x)+U\diag(y)U^*)\) forces
\(|U_{ij}|^2=1/m\). For the unit phase ratios \(z_1,\ldots,z_m\)
of two rows, the two-by-two minor coefficients give
\[
 \frac{2-2\operatorname{Re}(z_j\overline z_k)}{m^2}
   =\binom m2^{-1},
 \qquad
 \operatorname{Re}(z_j\overline z_k)=-\frac1{m-1}\quad(j\ne k).
\]
Thus the \(z_j\) would form a regular \((m-1)\)-simplex in
\(\mathbb R^2\), with Gram matrix of rank \(m-1\). This forces \(m\le3\).

This also rules out every homogeneous PSD determinant representation
\(G=\det(\sum_i x_i A_i+\sum_j y_j B_j)\), including a positive
scalar prefactor. Indeed, a nonzero such determinant has size \(m\).
At the positive matrix \(P=\sum_iA_i+\sum_jB_j\), multiaffinity
forces each \(A_i,B_j\) to have rank one, since
\(\det(P+(t-1)A_i)\) has degree \(\operatorname{rank}A_i\).
Translation makes
\(\det(P+t(\sum_iA_i-\sum_jB_j))\) constant.
The normalized difference is Hermitian and all its eigenvalues vanish,
so \(\sum_iA_i=\sum_jB_j=:R>0\).
After conjugation by \(R^{-1/2}\), each group consists of \(m\)
rank-one matrices summing to \(I_m\), hence is an orthonormal basis.
The pure-half coefficient one fixes the scalar factor and reduces the
representation to the two-basis form already excluded.

Appendix~\ref{app:matching-obstruction} gives the reverse obstruction
in every rank \(m\ge4\), persisting under perturbation to uniform support.
\end{proof}

\subsection{The exact transverse infinitesimal gap}
\label{sec:entropy-local-spectrum}

Let \(f_p=G(pz_A,qz_{A'})\) be a connected sharp law of rank
\(m\ge2\), and let \(M_{ij}\) be the coefficient of
\(x_{A\setminus\{i\}}y_j\) in \(G\).
Translation and complementation make \(M\) doubly stochastic.
On centered \(L^2(\mu_p)\) define
\[
 \mathcal L u=\frac12\sum_{i\in A\cup A'}\E[u\mid X_i],
 \qquad K=|B\cap A|.
\]

\begin{proposition}[A field-independent local spectrum]
\label{efg:local-spectrum}
Write the singular values of \(M\) as
\(1=s_1>s_2\ge\cdots\ge s_m\ge0\).
The nonzero spectrum of \(\mathcal L\), with multiplicities, is
\[
 1,\qquad \frac{1+s_j}{2},\ \frac{1-s_j}{2}\quad(2\le j\le m).
\]
Its eigenvalue-one space is \(\mathbb R(K-mp)\), and it vanishes
on the orthogonal complement of the coordinate span.
Thus the optimal infinitesimal half-marginal entropy coefficient
on \((K-mp)^\perp\) is \((1+s_2)/2\), with gap
\((1-s_2)/2\), independently of \(p\).
For matrix perturbations $F_\varepsilon=\sigma+\varepsilon H$ about a
faithful constant state, the same coefficient holds on the transverse space
\[
 \E H=0,\qquad \E[(K-mp)H]=0.
\]
It is attained even with pointwise trace-one perturbations.
\end{proposition}

\begin{proof}
The independent halves give orthonormal families
\[
 Z_i=\frac{X_i-p}{\sqrt{pq}}\quad(i\in A),\qquad
 W_j=\frac{1-X_j-p}{\sqrt{pq}}\quad(j\in A').
\]
Translate \(\partial_{x_i}\partial_{y_j}G\) from
\((p\mathbf1,q\mathbf1)\) to \((\mathbf1,0)\).
There it equals \(1-M_{ij}\), so
\(\operatorname{Cov}(X_i,1-X_j)=pqM_{ij}\).
Consequently the Gram matrix of \((Z,W)\) is
\(\left(\begin{smallmatrix}I&M\\M^{\mathsf T}&I\end{smallmatrix}\right)\).
The operator \(\mathcal L\) is half its frame operator; its nonzero
eigenvalues are therefore half the nonzero Gram eigenvalues.

For simplicity of \(s_1\), put \(R=\diag(pI,qI)\).
The same covariance calculation gives
\[
 R^{-1}\nabla^2f_p(\mathbf1)R^{-1}
 =\mathbf1\mathbf1^{\mathsf T}
       -\begin{pmatrix}I&M\\M^{\mathsf T}&I\end{pmatrix}.
\]
The Hessian kernel equals the translation space
by Lemma~\ref{lem:translation-kernel}, which is one-dimensional by
Theorem~\ref{efg:sharp-fields}.
Since \(m\ge2\), the constant vector has nonzero Hessian eigenvalue
in this congruence; the singular-vector decomposition then makes
\(s_1=1\) simple and gives the stated spectrum.
The common coordinate sum is
\(\sum_i Z_i=\sum_jW_j=(K-mp)/\sqrt{pq}\).

Finally the second entropy variation at a faithful constant
\(\sigma\) is the BKM quadratic form
\(\frac12\E\Tr[H\,D\log(\sigma)[H]]\).
Conditional expectations act only on the classical variable, so the
same frame calculation is tensored with the matrix tangent space.
Scalar perturbations, or \(H(B)=u(B)T\) with a traceless qubit matrix
\(T\), attain the corresponding Rayleigh quotients.
\end{proof}

The value \((1+s_2)/2\) is the optimal coefficient for the second
variation on the stated transverse space. On the full tangent space,
the directions \((K-mp)T\) give coefficient one, including under
pointwise trace one when \(T\) is traceless.

For \(G_{\mathrm{sym}}\), \(M=m^{-1}\one\one^{\mathsf T}\), so for \(m\ge2\)
the coordinate spectrum is \(1\) once and \(1/2\) with multiplicity \(2m-2\).
Its maximal transverse gap \(1/2\) coexists, for \(m\ge4\), with
the absence of a homogeneous PSD determinant representation.

\appendix
\section{Exact realization of entropy profiles}
\label{app:profile-cone}
The profile inequalities in Section~\ref{sec:entropy-profiles} are
complete when the reference law, finite ground set, and encoding
vary. The following construction realizes every allowed profile
using only a classical register and a partition matroid.

\begin{theorem}[Exact realization of every entropy profile]
\label{thm:exact-profile-cone}
Fix $r\ge1$. As the finite ground set, Lorentzian reference law, and
finite-dimensional encoding vary, the profiles
\[
 (I(\mathsf A_0:Q),\ldots,I(\mathsf A_r:Q))
\]
with classical marginal equal to the reference law are exactly the
sequences $h=(h_0,\ldots,h_r)$ satisfying
\begin{equation}\label{eq:exact-profile-cone}
 h_0=0,\qquad h_1\ge0,\qquad
 h_{j+1}-2h_j+h_{j-1}\ge0\quad(1\le j<r).
\end{equation}
Every such profile has an exact realization with a classical register
and the uniform basis law of a partition matroid.
\end{theorem}
\begin{proof}
Necessity is Theorem~\ref{thm:entropy-profile}. For sufficiency, decompose
$h$ into ramps:
\[
 h_k=\sum_{j=0}^{r-1}c_j(k-j)_+,
 \qquad c_0=h_1,\qquad
 c_j=h_{j+1}-2h_j+h_{j-1}\quad(j\ge1).
\]
Choose a prime power $q\ge r$ with $\log q\ge\sum_jc_j$, and let
$X_1,\ldots,X_r$ be independent uniform elements of $\mathbb F_q$.
The associated partition matroid chooses one label $(a,X_a)$ from
each of $r$ classes; its normalized basis polynomial is
\[
 f(z)=q^{-r}\prod_{a=1}^r\left(\sum_{v\in\mathbb F_q}z_{a,v}\right).
\]

For $0\le j<r$, choose a linear map
$A_j:\mathbb F_q^r\to\mathbb F_q^{r-j}$ whose kernel is an
$[r,j]$ MDS code, and put $Q_j=A_jX$. Such codes are obtained by
evaluating polynomials of degree less than $j$ at $r$ distinct field
points; for $j=0$ use the zero code. Conditional on $Q_j$, the vector
$X$ is uniform on a coset of this code. Its projection onto any
$k$ coordinates has dimension $\min(k,j)$. Hence, for every
$S\subseteq[r]$ with $|S|=k$,
\[
 I(X_S:Q_j)=\bigl(k-\min(k,j)\bigr)\log q=(k-j)_+\log q.
\]
Uniform retention reveals $S$ and $X_S$, and $S$ is independent of
$(X,Q_j)$. Thus the same formula holds for $I(\mathsf A_k:Q_j)$.
This is the MDS ramp-secret-sharing construction
\cite{KuriharaUyematsuMatsumoto} in the language of retained bases.

Choose an independent visible flag $J$ with
$\Pr(J=j)=c_j/\log q$, and assign its remaining probability to an
uninformative encoding. More explicitly, on this extra flag value
take $Q_J$ to be constant. Since $J$ is independent of the input,
the chain rule gives
\[
 I(\mathsf A_k:J,Q_J)
 =\sum_{j=0}^{r-1}\frac{c_j}{\log q}
                  I(\mathsf A_k:Q_j)=h_k.
\]
All alphabets are finite; arbitrary real values of the $c_j$ are
carried by the flag probabilities. This proves exact realization.
\end{proof}

The theorem also identifies the cone of all matrix-entropy profiles
when the classical marginal may vary: every such profile
satisfies \eqref{eq:exact-profile-cone}, and every point of that cone
already has the fixed-marginal classical realization above.
The cone is simplicial, with extreme rays
$((k-j)_+)_{k=0}^r$. Quantum side information realizes no additional
profiles. In particular, a linear inequality
\[
 \sum_{k=1}^r a_kH_k\ge0
\]
holds for all these profiles if and only if
\[
 \sum_{k=j+1}^r a_k(k-j)\ge0\qquad(0\le j<r).
\]
This exact finite realization complements the Shannon descriptions of
partition-symmetric entropy regions in \cite{ChenYeung}.

\begin{corollary}[Sharpness of every chord]
\label{cor:sharp-profile-chords}
For $0\le j<\ell<r$, the coefficient $(\ell-j)/(r-j)$ in
Corollary~\ref{cor:conditional-quantum-information} is attained.
More generally, for $j>0$ and every feasible pair
$H_j=a\ge0$, $H_r=b\ge ra/j$, equality in
\eqref{eq:entropy-profile-chord} is simultaneously attained at all
$j<\ell<r$. Every coefficient of the Green-kernel deficit formula
\eqref{eq:profile-green-deficit} is separately attained by a profile
with just that second difference nonzero.
\end{corollary}
\begin{proof}
Use the $j$th ramp for the first assertion. For the second, realize
\[
 h_k=\begin{cases}
 ak/j,&k\le j,\\[1mm]
 a+(b-a)(k-j)/(r-j),&k\ge j.
 \end{cases}
\]
This sequence is convex precisely under the stated endpoint condition.
A single ramp isolates one second difference and proves the last claim.
\end{proof}

\section{Determinantal laws without positive matching decompositions}
\label{app:matching-obstruction}

Here is the determinantal obstruction used in
Proposition~\ref{thm:balanced-incomparability}, with the matching
products \(P_\sigma\) defined there.

\begin{proposition}[An exact obstruction to matching mixtures]
\label{efg:matching-obstruction}
For every \(m\ge4\), there are connected sharp laws with uniform support
\(U_{m,2m}\) whose balanced polynomials are not nonnegative sums of the
\(P_\sigma\).
\end{proposition}

\begin{proof}
Start with the rational orthogonal matrix
\[
 U=\frac1{45}
 \begin{pmatrix}
 32&-6&-2&31\\
 -30&0&-15&30\\
 -10&-15&40&10\\
 -1&-42&-14&-8
 \end{pmatrix},
 \qquad UU^{\mathsf T}=I,
\]
and put \(G_U=\det(\diag(x)+U\diag(y)U^{\mathsf T})\).
If \(G_U=\sum_\sigma w_\sigma P_\sigma\), \(w_\sigma\ge0\),
the pure \(x\) coefficient gives \(\sum_\sigma w_\sigma=1\).
Comparing coefficients of \(x_{I^c}y_J\) gives
\[
 |\det U_{I,J}|^2=\Pr(\sigma(I)=J),\qquad |I|=|J|.
\]
Every permutation of four labels satisfies
\[
 \{\sigma\{1,3\}=\{1,4\}\}\subseteq
 \{\sigma(2)=2\}\cup\{\sigma(3)=4\}
                   \cup\{\sigma\{1,2\}=\{3,4\}\}.
\]
Indeed, outside the middle event the left event forces
\(\sigma(3)=1,\sigma(1)=4\); then \(\sigma(2)\) is either \(2\) or \(3\).
Its required probability inequality is violated by the exact values
\[
 |\det U_{13,14}|^2=\frac{196}{2025}
 >|U_{22}|^2+|U_{34}|^2+|\det U_{12,34}|^2
 =0+\frac{100}{2025}+\frac{81}{2025}.
\]
The gap is \(1/135\).
It persists under sufficiently small orthogonal perturbations.
A generic such perturbation has all square minors nonzero, so its
two-basis law has uniform support \(U_{4,8}\), and is connected.
It remains a sharp determinantal law by the preceding field construction.
For \(m>4\), embed the displayed matrix as \(U\oplus I_{m-4}\).
The same permutation inclusion holds with the additional event
\(\{\sigma(2)>4\}\) on the right. Its determinantal probability is
\(\sum_{j=5}^m|U_{2j}|^2=0\) at this embedding, so the strict violation
again persists under generic orthogonal perturbation to uniform support.
\end{proof}

This excludes positive matching mixtures with the prescribed common
Bernoulli weights.

\section*{Acknowledgment}
The author used ChatGPT for mathematical exploration and proof development, literature searches, computational experiments, and code for numerical certificates. Proofs developed with this assistance were carefully checked, adapted, and rewritten by the author, who also rewrote the exposition. The author takes full responsibility for the mathematical content, computations, and any errors.

\end{document}